\documentclass[11pt]{article}
\usepackage{amsmath,amsbsy,amsfonts}
\usepackage[frenchb]{babel}
\usepackage{makeidx}
 \newtheorem{thm}{Th\'eor\`eme}[subsection]
 
 \newtheorem{lem}[thm]{Lemme}
 \newtheorem{prop}[thm]{Proposition}
\title{Holomorphie des op\'erateurs d'entrelacement normalis\'es \`a l'aide des param\`etres d'Arthur}
\author{C. M{\oe}glin\\ CNRS, Institut de Math\'ematiques de Jussieu\\ moeglin@math.jussieu.fr}
\date{}
\begin{document}
\maketitle
 
\section*{Introduction}
Le but de cet article est de prouver que les op\'erateurs d'entrelacement qui interviennent dans les s\'eries d'Eisenstein construites avec des formes automorphes de carr\'e int\'egrable et des paraboliques maximaux sont holomorphes au voisinage de l'axe r\'eel positif. Ici on ne travaille que localement donc on r\'eexprime le r\'esultat diff\'eremment; on renvoie le lecteur \`a \cite{seriedeisenstein} pour les rapports locaux/globaux. On fixe donc un gros groupe $\mathcal{G}$ classique est un parabolique maximal de ce groupe donc de la forme $GL(D)\times G$ o\`u $G$ est un groupe classique de m\^eme type que ${\mathcal{G}}$; c'est $G$ qui intervient dans cet article. On fixe une repr\'esentation irr\'eductible $\pi$ de $G$; on suppose qu'elle est composante locale d'une forme automorphe de carr\'e int\'egrable; on oublie le fait que la conjecture de Ramanujan n'est pas connue pour les groupes $GL$, cela n'a pas d'importance s\'erieuse pour nos probl\`emes d'holomorphie et on suppose donc que $\pi$ est dans un paquet d'Arthur; on rappelle les d\'efinitions et constructions ci-dessous. On fixe une repr\'esentation de Steinberg de $GL(D)$, not\'e ici $\sigma$ et on d\'efinit une fonction m\'eromorphe $r(\sigma,\psi,s)$ qui d\'epend de $\sigma$ et de $\psi$ (ceci est fait en \ref{normalisation}) mais on peut dire simplement que $r(\sigma,\psi,s)$ est le facteur de normalisation d\'efini au moins th\'eoriquement par Langlands et Shahidi pour l'entrelacement $\sigma\vert\,\vert^{s}\times \pi_{L} \rightarrow \sigma^*\vert\,\vert^{-s}\times \pi_{L}$ o\`u $\pi_{L}$ est une des repr\'esentations du paquet de Langlands inclus dans le paquet d'Arthur d\'efini par $\psi$. Donc en g\'en\'eral la fonction $r(\sigma,\pi,s)$ n'est pas le facteur de normalisation de Langlands-Shahidi pour l'entrelacement pr\'ec\'edent quand on remplace $\pi_{L}$ par $\pi$ sauf si $\psi$ est temp\'er\'e c'est-\`a-dire trivial sur la 2e copie de $SL(2,{\mathbb C})$; c'est le seul cas o\`u paquet de Langlands et paquet d'Arthur co\"{\i}ncident. On note $M(\sigma,\pi,s)$ l'op\'erateur d'entrelacement standard associ\'e comme pr\'ec\'edemment \`a l'induite $\sigma\vert\,\vert^s\times \pi$ et on pose $N(\sigma,\pi,s):=r(\sigma,\psi,s)^{-1}M(\sigma,\pi,s)$ et on montre dans cet article que $N(\sigma,\pi,s)$ est holomorphe au voisinage de l'axe r\'eel positif.

L'article commence par un rappel des constructions des repr\'esentations dans un paquet d'Arthur;  on simplifie la construction dans le cas d'un morphisme g\'en\'eral en \ref{assouplissement} ci-dessous; cette simplification a son int\'er\^et en soi. L'id\'ee est toujours la m\^eme on sait d\'efinir avec une grande pr\'ecision les repr\'esentations dans un paquet d'Arthur associ\'e \`a un morphisme $\psi$ si la restriction de $\psi$ \`a $W_{F}\times SL(2,{\mathbb C})$ o\`u $SL(2,{\mathbb C})$ est plong\'e diagonalement dans le produit $SL(2,{\mathbb C})\times SL(2,{\mathbb C})$ est sans multiplicit\'e. Le point est d'obtenir le cas g\'en\'eral \`a partir de ce cas particulier en utilisant des modules de Jacquet. Dans la d\'efinition orginale on avait rigidifi\'e les choix de fa\c{c}on \`a avoir une param\'etrisation la plus canonique possible; en fait cette param\'etrisation est surtout beaucoup trop rigide et peu compatible aux proc\'edures standard d'induction restriction. On s'affranchit donc de toute rigidit\'e, on perd encore plus de renseignements sur la param\'etrisation mais on peut alors faire des choix dans chaque situation et les rendre compatibles \`a l'induction ou la restriction que l'on consid\`ere. Dans cette introduction, on ne peut pas expliquer toutes les notations; mais il est certainement assez clair que $\psi$ peut \^etre vu comme une repr\'esentation de $W_{F}\times SL(2,{\mathbb C})\times SL(2,{\mathbb C})$ dans un groupe d'automorphismes d'une forme bilin\'eaire d'un type donn\'e (\'eventuellement ne respectant la forme bilin\'eaire qu'\`a un scalaire pr\`es). On dit que $\psi$ a bonne parit\'e si toutes les sous-repr\'esentations irr\'eductibles incluses dans $\psi$ sont aussi \`a valeurs dans un groupe d'automorphismes d'une forme de m\^eme type. On avait montr\'e en \cite{general} repris en \cite{pourshahidi} que par une induction qui pr\'eserve l'irr\'eductibilit\'e on se ram\`ene \`a d\'ecrire les paquets de repr\'esentations associ\'es \`a un $\psi$ g\'en\'eral au cas des $\psi$ ayant bonne parit\'e. Ensuite on construit les paquets associ\'es \`a un morphisme $\psi$ ayant bonne parit\'e en construisant un morphisme $\psi_{>>}$ ``dominant'' $\psi$ et tel que $\psi_{>>}$ soit de restriction discr\`ete \`a la diagonale et en trouvant les repr\'esentations associ\'ees \`a $\psi$ par module de Jacquet \`a partir de celles associ\'ees \`a $\psi_{>>}$ (cf. ci-dessous \ref{assouplissement}); c'est la notion de ``dominant'' que l'on a simplifi\'ee. On dit que $\psi_{>>}$ domine $\psi$ s'il existe un ordre total sur l'ensemble des sous-repr\'esentations irr\'eductibles incluses dans $\psi_{>>}$ et sur l'ensemble des sous-repr\'esentations irr\'eductibles incluses dans $\psi$ (il y a une ambiguit\'e si $\psi$ a de la multiplicit\'e) et une bijection respectant les ordres, du premier ensemble sur le second compatible \`a l'action de $W_{F}$; de fa\c{c}on explicite en identifiant une repr\'esentation irr\'eductible de $SL(2,{\mathbb C})$ \`a sa dimension, on voit les sous-repr\'esentations irr\'eductibles incluses dans $\psi$ comme un triplet $(\rho,a_{>>},b_{>>})$ et la bijection envoie un tel triplet sur un triplet $(\rho,a,b)$ tel que $inf(a,b)=inf(a_{>>},b_{>>})$, $sup(a_{>>},b_{>>})\geq sup(a,b)$ et $(a_{>>}-b_{>>})(a-b)\geq 0$. Remarquons que l'on ne peut pas trouver de morphisme dominant $\psi$ et de restriction discr\`ete \`a la diagonale si  $\psi$ n'a pas bonne parit\'e. 

La preuve de l'holomorphie des op\'erateurs d'entrelacement se fait de la fa\c{c}on la plus standard qui soit: on s'arrange pour inclure $\pi$ dans une induite $\tau\times \pi'$ o\`u $\pi'$ sera encore dans un paquet d'Arthur et on d\'ecompose l'op\'erateur d'entrelacement standard en le produit des 3 op\'erateurs d'entrelacements standard \'evidents et on montre que ce produit n'a pas de p\^ole en utilisant une hypoth\`ese de r\'ecurrence pour $\tau'$. On explique en d\'etail les choix en \ref{shemadelapreuve}; le point de d\'epart est pour $s\in \mathbb{R}_{>0}$, le r\'esultat d'Harish-Chandra qui montre que si $\pi$ est une s\'erie discr\`ete, alors l'op\'erateur d'entrelacement standard est holomorphe pour $Re\, s>0$ (cf. par exemple \cite{waldsJIMJ});  en $s=0$, c'est uniquement le fait que l'op\'erateur d'entrelacement normalis\'e est autodual et les applications sont pour $s>0$.
\tableofcontents
\section{Quelques notations g\'en\'erales\label{notations}}
Dans tout le travail $F$ est un corps $p$-adique.
On consid\`ere ici un groupe $G$ classique; le prototype sont les groupes orthgonaux symplectiques ou unitaires mais gr\^ace \`a une remarque d'Arthur (\cite{arthursp4}) on peut leur adjoindre les groupes $Gspin(2n+1)$ et les groupes non connexes $Gspin(2n)$ et pour les m\^emes raisons on peut obtenir les groupes $GU(m,F'/F)$ o\`u $F'$ est une extension quadratique de $F$. Chacun de ces groupes (hormi le cas des groupes unitaires) a un groupe dual de la forme $G^*\times W_{F}$ o\`u $W_{F}$ est le groupe de Weyl de $F$ et o\`u $G^*$ est un groupe classique. Si $G$ est un groupe unitaire, $G^*$ sera plut\^ot un groupe lin\'eaire, ceci est largement \'ecrit dans la litt\'erature. On consid\`ere la repr\'esentation naturelle de $G^*$ dans un groupe $GL(m^*_{G},{\mathbb C})$ (cela d\'efinit $m^*_{G}$); si $G$ est un groupe $Gspin$ ou $GU$ il faut prendre $GL^(m^*_{G},{\mathbb C})\times {\mathbb C}^*$ comme expliqu\'e par exemple dans \cite{arthursp4}. On note $\theta^*$ l'automorphisme $g\mapsto \, ^tg^{-1}$ (ou $(g,\lambda)\mapsto (\, ^tg^{-1})\lambda,\lambda)$) de $GL(m^*_{G},{\mathbb C})$ (ou $GL(m^*_{G},{\mathbb C})\times {\mathbb C}^*)$) et on note $\theta$ l'automorphisme ''dual'' et on impose \`a $\theta$ de respecter un \'epinglage.

On consid\'erera des morphismes de $W_{F}\times SL(2,{\mathbb C})\times SL(2,{\mathbb C})$ dans $G^*$ continus-alg\'ebriques, born\'es et semi-simples; en composant avec la repr\'esentation naturelle cela donne une repr\'esentation de $W_{F}\times SL(2,{\mathbb C})\times SL(2,{\mathbb C})$. On note g\'en\'eriquement $\psi$ un tel morphisme et on note $Jord(\psi)$ l'ensemble des sous-repr\'esentations irr\'eductibles incluses dans $\psi$; toute repr\'esentation irr\'eductible de $SL(2,{\mathbb C})$ est uniquement d\'etermin\'ee par sa dimension et on identifie donc $Jord(\psi)$ \`a un ensemble de triplets $(\rho,a,b)$ o\`u $\rho$ est une repr\'esentation irr\'eductible de $W_{F}$ et $a,b$ sont des entiers. On tient compte des multiplicit\'es de $\psi$ dans $Jord(\psi)$. On appelle caract\`ere du centralisateur de $\psi$ une application, $\epsilon$, de $Jord(\psi)$ dans $\{\pm 1\}$ et on appelle restriction de $\epsilon$ au centre de $G^*$, le signe $\epsilon(z_{G^*}):=\times_{(\rho,a,b)\in Jord(\psi)}\epsilon (\rho,a,b)$ o\`u les multiplicit\'es de $Jord(\psi)$ sont prises en compte. Et on dit que la restriction de $\epsilon$ au centre de $G^*$ est d\'etermin\'ee par le type de $G$ si $\epsilon(z_{G^*})=+1$ exactement quand $G$ est d\'etermin\'e par une forme bilin\'eaire dont l'invariant de Hasse est $+1$. C'est une version combinatoire des notions usuelles.

Gr\^ace aux travaux de Zelevinsky on sait donc associer une repr\'esentation irr\'eductible de $GL(m^*_{G},F)$ not\'ee $\pi^{GL}(\psi)$. Cette repr\'esentation est stable par l'action de l'automorphisme ext\'erieur de $GL(m^*_{G},F)$, not\'e $\theta$ ci-dessus et on prolonge cette repr\'esentation au produit semi-direct de $GL(m^*_{G},F)\times\{1,\theta\}$; ici il n'est pas important de fixer le prolongement. On appelle $\theta$-trace la trace de cette repr\'esentation limit\'ee \`a la composante non connexe de ce produit semi-direct; on identifie une telle trace \`a une fonction sur $GL(n,F)$ invariante par $\theta$-conjugaison.

\

On utilisera aussi la notation suivante bien commode: soit $\pi$ une repr\'esentation de $G$ et $\rho$ une repr\'esentation cuspidale unitaire d'un groupe lin\'eaire $GL(d_{\rho},F)$, ce qui d\'efinit $d_{\rho}$. On suppose qu'il existe un groupe $G'$ de m\^eme type que $G$ et de  rang $d_{\rho}$ plus petit tel que $GL(d_{\rho},F)\times G'$ soit un sous-groupe de Levi de $G$. On note $x$ un r\'eel et $Jac_{\rho\vert\,\vert^x}\pi$ ou simplement $Jac_{x}\pi$ (quand $\rho$ est fix\'e) l'\'el\'ement du groupe de Grothendieck de $G'$, encore vu comme une repr\'esentation semi-simple de $G'$ tel que le module de Jacquet de $\pi$ relativement \`a un parabolique de Levi $GL(d_{\rho},F)\times G'$ soit de la forme $\rho\vert\,\vert^{x}\otimes Jac_{x}\pi \oplus _{\sigma'\not\equiv \rho\vert\,\vert^{x}, \pi'}\sigma'\otimes \pi'$, \'egalit\'e dans le groupe de Grothendieck. Soit $\pi^{GL}$ une repr\'esentation de $GL(m^*_{G},F)$; on d\'efinit de fa\c{c}on analogue $Jac^\theta_{x}\pi^{GL}$ comme repr\'esentation semi-simple de $GL(m^*_{G}-2d_{\rho},F)$ telle que le module de Jacquet de $\pi^{GL}$ relativement \`a un parabolique de Levi $GL(d_{\rho},F)\times GL(m^*_{G}-2d_{\rho},F)\times GL(d_{\rho},F)$ soit, dans le groupe de Grothendieck, de la forme $\rho\vert\,\vert^{x}\otimes Jac^\theta_{x}\pi^{GL}\otimes \theta(\rho\vert\,\vert^x)\oplus_{\sigma',\sigma'',\tau}\sigma'\otimes \tau\otimes \sigma''$ o\`u $(\sigma',\sigma'')$ parcourt l'ensemble des couples de repr\'esentations irr\'eductibles de $GL(d_{\rho},F)$ sauf le couple $(\rho\vert\,\vert^{x},\theta(\rho\vert\,\vert^x))$.

\

L'induction pour les repr\'esentations est not\'ee par $x$; soit $\rho$ une repr\'esentation cuspidale d'un groupe lin\'eaire $GL(d_{\rho},F)$ et soit $[x,y]$ un segment croissant ou d\'ecroissant. Gr\^ace aux travaux de Bernstein et de Zelevinsky, on sait que l'induite $\rho\vert\,\vert^{x}\times \cdots \times \rho\vert\,\vert^{y}$ a un unique sous-module irr\'eductible. Ce sous-module est not\'e, soit $<\rho\vert\,\vert^x, \cdots, \rho\vert\,\vert^y>$ soit de fa\c{c}on plus traditionnelle $<x, \cdots, y>_{\rho}$. Gr\^ace aux travaux de Zelevinsky, on peut g\'en\'eraliser cette construction a un ensemble de multisegments dans les cas suivants; consid\'erons une matrice, $\mathcal A$, (non n\'ecessairement rectangulaire), dont les lignes et les colonnes sont des segments de croissance diff\'erentes$$
\begin{matrix} &B_{1}&\cdots &\cdots &A_{1}\\
&\vdots &\vdots &\vdots
\\
B_{t}&\cdots&\cdots &A_{t}
\end{matrix}
$$
o\`u pour tout $i\in [1,t]$, $[B_{i},A_{i}]$ est un segment li\'e au sens de Zelevinsky \`a celui qui le pr\'ec\`ede et \`a celui qui le suit, soit croissant avec alors $B_{1}>\cdots >B_{t}$ et soit d\'ecroissant avec $B_{1}< \cdots <B_{t}$. Dans l'\'ecriture, l'\'el\'ement en dessous de $B_{1}$ est $B_{1}-1$, donc $B_{t}$ n'est pas n\'ecessairement d\'ecal\'e par rapport \`a $B_{1}$. Dans ce cas l'induite, pour le $GL$ convenable $<B_{1}, \cdots, A_{1}>_{\rho}\times \cdots <B_{t}, \cdots, A_{t}>_{\rho}$ a un unique sous-module irr\'eductible que l'on note $<\mathcal{A}>_{\rho}$. Le type m\^eme de repr\'esentations qui interviennent d\`es que l'on \'etudie les paquets d'Arthur sont les repr\'esentations $Speh(St(\rho,a),b)$ pour des entiers $a,b$ qui, dans cette \'ecriture sont les repr\'esentations associ\'ees \`a la matrice rectangulaire ci-dessous et aussi \`a la transpos\'ee de cette matrice.
$$\begin{matrix}
(a-b)/2 &\cdots &(a+b)/2-1\\
\vdots &\vdots &\vdots\\
-(a+b)/2+1 &\cdots &-(a-b)/2
\end{matrix}
$$

\bf Une convention importante\rm

On s'int\'eresse aux op\'erateurs d'entrelacement associ\'es \`a l'\'el\'em\'ent du groupe de Weyl de longueur maximale, pour des induites de la forme $\sigma\vert\,\vert^s\times \pi$ o\`u $\sigma$ est une repr\'esentation de Steinberg (unitaire) et $s\in {\mathbb R}_{\geq 0}$. En g\'e\'erale l'image est \`a valeurs dans une induite de la forme $\tilde{\sigma}\vert\,\vert^{-s}\times \pi$ o\`u $\tilde{\sigma}$ est la contragr\'ediente de $\sigma$ parfois tordu par un caract\`ere d\'ependant de $\pi$ (cas des groupes de similitudes); pour simplifier la notation on suppose que $\tilde{\sigma}\simeq \sigma$. C'est le seul cas important.

\section{Rappel sur les paquets d'Arthur}
\subsection{D\'efinition des param\`etres}

On fixe $\psi$ comme dans \ref{notations}.
A un tel morphisme Arthur associe un ensemble fini de repr\'esentations irr\'eductibles de $G$, ensemble not\'e $\Pi(\psi)$ uniquement d\'etermin\'e par le fait qu'une combinaison lin\'eaire convenable \`a coefficient tous non nuls des traces des \'el\'ements de $\Pi(\psi)$ se transf\`ere en la $\theta$-trace de $\pi^{GL}(\psi)$.

Le seul cas particulier de ce r\'esultat dont nous avons besoin est celui beaucoup plus simple des paquets de s\'eries discr\`etes; on suppose que $\psi$ est trivial sur la 2e copie de $SL(2,{\mathbb C})$ et que la repr\'esentation d\'efinie par $\psi $ n'a pas de multiplicit\'e. Dans ce cas $\Pi(\psi)$ est form\'ee de s\'eries discr\`etes de $G$; leur nombre est exactement $2^{long(\psi)-1}$ o\`u $long(\psi)$ est la longueur de $\psi$ en tant que repr\'esentation de $W_{F}\times SL(2,{\mathbb C})$. Avec cela on a d\'eduit l'existence d'une classification des repr\'esentations cuspidales de $G$ sous la forme suivante (les d\'efinitions de sans trou et altern\'e sont apr\`es l'\'enonc\'e):
\begin{prop}Il existe une bijection entre l'ensemble des classes d'isomorphie de repr\'esentations cuspidales irr\'eductibles de $G$ et l'ensemble des couples $\psi,\epsilon$ o\`u $\psi$ est un morphisme discret et sans trou et $\epsilon$ est un caract\`ere alt\'ern\'e du centralisateur de $\psi$ de restriction au centre de $G^*$ impos\'ee par $G$.
\end{prop}
On dit que $\psi$ est discret si la repr\'esentation de $W_{F}\times SL(2,{\mathbb C})$ est sans multiplicit\'e; on dit que $\psi$ est sans trou si quand $Jord(\psi)$ contient $(\rho,a)$ avec $a>2$, il contient aussi $(\rho,a-2)$ . Le caract\`ere du centralisateur de $\psi$ est dit altern\'e si, vu comme application de $Jord(\psi)$ dans $\{\pm 1\}$ il prend des valeurs diff\'erentes sur $(\rho,a)$ et $(\rho,a-2)$ en prenant la valeur $-$ sur tout couple de la forme $(\rho,2)$.

\

L'application qui \`a une repr\'esentation cuspidale $\pi$ associe le morphisme $\psi$ est uniquement d\'etermin\'ee par la propri\'et\'e suivante: soit $\rho$ une repr\'esentation cuspidale de $GL(d_{\rho},F)$ autoduale. On sait alors d'apr\`es Silberger qu'il existe un unique r\'eel $x_{\rho,\pi}\in {\mathbb R}_{\geq 0}$ tel que l'induite $\rho\vert\,\vert^{x_{\rho,\pi}}\times \pi$ soit r\'eductible. On pose $Red(\pi):=\{(\rho,x_{\rho,\pi}); x_{\rho,\pi}>1/2\}$. Alors $Jord(\psi)=\{(\rho,2x_{\rho,\pi}-1); (\rho,x_{\rho,\pi})\in Red(\pi)\}$.

\

La seule d\'emonstration que je connais de ce r\'esultat utilise l'existence de $\Pi(\psi)$ dans le cas particulier des morphismes triviaux sur la 2e copie de $SL(2,{\mathbb C})$ et ceci est d\^u \`a Arthur et n\'ecessite le lemme fondamental ordinaire et tordu, ainsi que le transfert pour toutes les fonctions localement constantes \`a support compact.

\

Avec cela on peut utiliser \cite{europe}, \cite{ams} pour construire toutes les s\'eries discr\`etes de $G$, ou on peut utiliser les r\'esultats d'Arthur pour cette m\^eme construction; on a d\'ej\`a montr\'e que les constructions co\"{\i}ncident. On ne fait alors plus d'hypoth\`eses sur $\psi$ et on construit un paquet de repr\'esentations $\Pi(\psi)$ en \cite{paquetdiscret} repris en \cite{pourshahidi}. Et la construction est faite de telle sorte que la propri\'et\'e de transfert pour les $\psi$ discrets donne la propri\'et\'e de transfert pour les $\psi$ g\'en\'eraux. Ici on n'a pas besoin de savoir de quelle combinaison lin\'eaire il s'agit. Les paquets que nous avons construits sont donc les paquets qu'Arthur trouvent (une fois tous les lemmes fondamentaux d\'emontr\'es) par son \'etude des formules des traces. Arthur annonce que la somme des s\'eries discr\`etes dans un paquet $\Pi(\psi)$ avec $\psi$ de restriction triviale \`a la 2e copie de $SL(2,{\mathbb C})$ est la (\`a un scalaire pr\`es) combinaison lin\'eaire stable dans le paquet. On a alors montr\'e en \cite{paquetdiscret} quelle est la combinaison lin\'eaire stable dans le paquet $\Pi(\psi)$; les coefficients sont des signes explicites. Ici on a besoin d'une description/construction des \'el\'ements de $\Pi(\psi)$ pour un $\psi$ g\'en\'eral que l'on va rappeler et que l'on simplifiera un peu ci-dessous. Dans cet article, on va utiliser le calcul des coefficients  mais on pourrait tr\`es bien garder des coefficients non calcul\'es, cela ne changerait rien.

\

Fixons $\psi$; on note {$\psi_{\vert \Delta}$} la restriction de $\psi$  \`a $W_{F}\times SL(2,{\mathbb C})$ o\`u $SL(2,{\mathbb C})$ est plong\'e diagonalement dans $SL(2,{\mathbb C})\times SL(2,{\mathbb C})$.

Pour construire les \'el\'ements de $\Pi(\psi)$ on commence par le faire dans des cas o\`u $\psi_{\vert \Delta}$ est sans multiplicit\'e; on appelle ce cas le cas o\`u  $\psi$ est de restriction discr\`ete \`a la diagonale. Ceci a \'et\'e trait\'e en d\'etail dans \cite{paquetdiscret} avec les r\'esultats de \cite{selecta}.

Sous cette hypoth\`ese $\Pi(\psi)$ est en bijection avec les couples $\underline{t}, \underline{\eta}$ d'applications d\'efinies sur $Jord(\psi)$ o\`u pour tout $(\rho,a,b)\in Jord(\psi)$, $\underline{t}(\rho,a,b)\in [0,[inf(a,b)/2]]$ et $\underline{\eta}(\rho,a,b)\in \{\pm 1\}$ avec $\underline{\eta}(\rho,a,b)=+$ si $\underline{t}(\rho,a,b)=inf(a,b)/2$; de plus \`a un tel couple $\underline{t},\underline{\eta}$ on associe (cf. ci-dessous) un caract\`ere du centralisateur de $\psi$, pour nous une application, $\epsilon_{\underline{t},\underline{\eta}}$ de $Jord(\psi)$ dans $\{\pm 1\}$ et ce caract\`ere \`a sa restriction au centre de $G^*$ d\'etermin\'ee par la forme de $G$, c'est-\`a-dire pour nous $\times_{(\rho,a,b)\in Jord(\psi)}\epsilon_{\underline{t},\underline{\eta}}(\rho,a,b)=\epsilon_{G}$ o\`u $\epsilon_{G}$ est l'invariant de Hasse intervenant dans la d\'efinition de la forme bilin\'eaire que $G$ respecte (\'eventuellement \`a un scalaire pr\`es). Par d\'efinition:
$$\forall (\rho,a,b)\in Jord(\psi),
\epsilon_{\underline{t}, \underline{\eta}}(\rho,a,b)=\underline{\eta}(\rho,a,b)^{inf(a,b)}(-1)^{[inf(a,b)/2]+\underline{t}(\rho,a,b)}.
$$
On note $\pi(\psi,\underline{t},\underline{\eta})$ ou plus simplement $\pi(\underline{t},\underline{\eta})$ la repr\'esentation correspondant \`a $\underline{t},\underline{\eta}$. On pose aussi $\epsilon_{\underline{t},\underline{\eta}}(s_{\psi}):=\times_{(\rho,a,b)\in Jord(\psi)}\epsilon_{\underline{t},\underline{\eta}}(\rho,a,b)^{b-1}$ et on a  a alors montr\'e en loc.cite que $\sum_{\underline{t},\underline{\eta}}\epsilon_{\underline{t},\underline{\eta}}(s_{\psi}) \pi(\underline{t},\underline{\eta})$ est une distribution stable qui a le bon transfert pour un bon choix d'action de $\theta$ qui ne nous importe pas ici puisque cela ne change la combinaison lin\'eaire que par un signe d\'ependant uniquement de $\psi$ et de l'action de $\theta$ fix\'ee. Ce qui est important pour ce travail est d'avoir une description en termes de repr\'esentations des $\pi(\psi,\underline{t},\underline{\eta})$. On rappelle cette description telle qu'elle a \'et\'e faite en loc.cite; elle se fait par r\'ecurrence sur $\ell(\psi):=\sum_{(\rho,a,b)\in Jord(\psi)}(inf(a,b)-1)$. 

\subsection{Le cas des morphismes \'el\'ementaires de restriction discr\`ete \`a la diagonale}
Le premier cas est donc celui des morphismes v\'erifiant $\ell(\psi)=0$ que l'on a appel\'e morphismes \'el\'ementaires; on garde l'hypoth\`ese que la restriction $\psi_{\vert \Delta}$ est sans multiplicit\'e. Dans ce cas n\'eces\-sairement $\underline{t}\equiv 0$ et $\epsilon_{\underline{t},\underline{\eta}}\equiv \underline{\eta}$. On note $\psi_{\vert \Delta}$ la restriciton de $\psi$ \`a $W_{F}$ fois la diagonale de $SL(2,{\mathbb C})\times SL(2,{\mathbb C})$. On remarque que l'application $(\rho,a,b)\in Jord(\psi) \rightarrow (\rho,sup(a,b))\in Jord(\psi_{\vert \Delta})$ est une bijection; on sait donc d\'efinir une s\'erie discr\`ete associ\'ee \`a $\psi_{\vert \Delta}$ et au carac\`etre $\underline{\eta}$ que l'on note $\pi(\psi_{\vert \Delta},\underline{\eta})$ et la repr\'esentation $\pi(\psi,\underline{t}, \underline{\eta})$ est l'image de $\pi(\psi_{\vert \Delta}, \underline{\eta})$ par une application qui g\'en\'eralise l'involution d'Iwahori-Matsumoto. Mais en \cite{elementaire} on en a donn\'e une description simple qui distingue 3 cas, le 2e et le 3e ne sont pas exclusifs l'un de l'autre. On  fait ici une r\'ecurrence sur $\sum_{(\rho,a,b)\in Jord(\psi)}sup(a,b)$.

1e cas: $\psi_{\vert \Delta}$ est sans trou et $\underline{\eta}$ est altern\'e; la repr\'esentation associ\'ee est alors cuspidale et son existence r\'esulte de la classification des repr\'esentations cuspidales

2e cas: le morphisme $\psi_{\vert \Delta}$ a un trou; il existe donc $(\rho,a,b)\in Jord(\psi)$ avec $sup(a,b)>2$ et $(\rho,sup(a,b)-2)\notin Jord(\psi_{\vert \Delta})$; on inclut ici aussi le cas o\`u $sup(a,b)=2$ et $\underline{\eta}(\rho,a,b)=+$. On note $(\rho,a',b')$ le triplet tel que $(a-b)(a'-b')\geq 0$ et $sup(a,b)-2=sup(a',b')$ et on note $\psi'$ le morphisme qui se d\'eduit de $\psi$ en rempla\c{c}ant $(\rho,a,b)$ par $(\rho,a',b')$. On identifie naturellement $\underline{\eta}$ \`a un morphisme de $Jord(\psi')$ dans $\{\pm 1\}$. On sait donc d\'efinir $\pi(\psi',\underline{t}\equiv 0, \underline{\eta})$ et on a montr\'e que l'induite
$
\rho\vert\,\vert^{(a-b)/2)}\times \pi(\psi',\underline{t},\underline{\eta})$ a un unique sous-module irr\'eductible et on a d\'efini $\pi(\psi,\underline{t}\equiv 0, \underline{\eta})$ comme ce sous-module irr\'eductible.

3e cas: il existe $(\rho,a,b)\in Jord(\psi)$ avec $sup(a,b)> 2$ et $(\rho,sup(a,b)-2x)\in Jord(\psi_{\vert \Delta})$ pour tout $x\in [0, [inf(a,b)/2]]$ et $\underline{\eta}$ n'est pas altern\'e sur ces \'el\'ements mais l'est sur les \'el\'ement $(\rho,sup(a,b)-2x)$ pour tout $x\in [1,[inf(a,b)/2]]$. On pose $\zeta:=+$ si $a>1$ et $\zeta=-$ si $b>1$. On note $\psi'$ le morphisme qui se d\'eduit de $\psi$ en enlevant  de $Jord(\psi)$ les \'el\'ements $(\rho,a,b)$ et $(\rho,a',b')$ o\`u $sup(a',b')=sup(a,b)-2$. On restreint $\underline{\eta}$ en une application de $Jord(\psi')$ dans $\{\pm 1\}$. On sait d\'efinir, par r\'ecurrence, la repr\'esentation $ \pi(\psi',\underline{t})$.  On a montr\'e que la repr\'esentation induite
$$
<\rho\vert\,\vert^{ (a-b)/2}, \cdots, \rho\vert\,\vert^{- (a'-b')/2}>\times \pi(\psi',\underline{\eta})
$$
a exactement 2 sous-modules irr\'eductibles. Et la repr\'esentation cherch\'ee $\pi(\psi,\underline{\eta})$ est l'un de ces sous-modules; on a pr\'ecis\'e lequel mais cela n'a pas d'importance ici.

\subsection{Le cas g\'en\'eral des morphismes de restriction discr\`ete \`a la diagonale}
Ici on suppose que $\psi$ est tel que $\psi_{\vert \Delta}$ est une repr\'esentation sans multiplicit\'e de $W_{F}\times SL(2,{\mathbb C})$ et on suppose que $\ell(\psi)>0$. On fixe $(\rho,a,b)$ avec $inf(a,b)>1$. On note  $\zeta$ le signe de $a-b$ et on pose $\zeta=+$ si $a=b$. On fixe encore $\underline{t}$ et $\underline{\eta}$ satisfaisant aux conditions g\'en\'erales et on d\'ecrit $\pi(\psi,\underline{t},\underline{\eta})$ suivant les cas.

1e cas: $\underline{t}(\rho,a,b)=0$; on note $\psi'$ le morphisme qui se d\'eduit de $\psi$ en changeant uniquement $(\rho,a,b)$ en $\cup_{c\in [\vert a-b\vert +1,a+b-1]; c\equiv a+b-1[2]}(\rho,sup(1,\zeta c),sup (1,-\zeta c))$; on d\'efinit $\underline{t}'$ sur $Jord(\psi')$ en restreignant $\underline{t}$ \`a $Jord(\psi')\cap Jord(\psi)$ et en \'etendant ensuite par $0$ et on d\'efinit $\underline{\eta}'$ sur $Jord(\psi')$ en restreignant $\underline{\eta}$ \`a $Jord(\psi')\cap Jord(\psi)$ et en \'etendant par:
$$\forall c\in  [\vert a-b\vert +1,a+b-1]; c\equiv a+b-1[2],
\underline{\eta}'(\rho,sup(1,\zeta c), sup(1,-\zeta c))=\underline{\eta}(\rho,a,b) ^{(c-\vert a-b\vert-1)/2}.$$
On a clairement $\psi_{\vert \Delta}=\psi'_{\vert \Delta}$ et
$$
\epsilon_{\underline{t},\underline{\eta}}(\rho,a,b)=\underline{\eta}(\rho,a,b)^{inf(a,b)}(-1)^{[inf(a,b)/2]}=\times_{c\in  [\vert a-b\vert +1,a+b-1]; c\equiv a+b-1[2]}\underline{\eta}'(\rho,sup(1,\zeta c), sup(1,-\zeta c)).
$$
Donc on a bien la relation $$\times_{(\rho',a',b')\in Jord(\psi')}\epsilon_{\underline{t}',\underline{\eta}'}(\rho',a',b')= \times_{(\rho'',a'',b'')\in Jord(\psi)}\epsilon_{\underline{t}, \underline{\eta}}(\rho'',a'',b'')=\epsilon_{G}.$$
De plus $\ell(\psi')=\ell(\psi)-inf(a,b)+1<\ell(\psi)$ on sait d\'efinir $\pi(\psi',\underline{t'},\underline{\eta}')$ et on pose tout simplement:
$$
\pi(\psi,\underline{t},\underline{\eta})=\pi(\psi',\underline{t}',\underline{\eta}').
$$
2e cas: $\underline{t}(\rho,a,b)>0$; en particulier $inf(a,b)\geq 2$. On note $\psi'$ le morphisme qui se d\'eduit de $\psi$ en changeant simplement $(\rho,a,b)$ en $(\rho,a',b')$ o\`u $(a',b')$ est uniquement d\'etermin\'e par le fait que $sup(a',b')=sup(a,b)$, $inf(a',b')=inf(a,b)-2$ et $\zeta (a'-b')>0$; si $inf(a,b)=2$, on supprime simplement $(\rho,a,b)$ et on se rappelle qu'alors $\underline{\eta}(\rho,a,b)=+$. On d\'efinit $\underline{t}',\underline{\eta}'$ en prenant la restriction de $\underline{t}$ et $\underline{\eta}$ sur $Jord(\psi)\cap Jord(\psi')$ et en posant $\underline{\eta}'(\rho,a',b')=\underline{\eta}(\rho,a,b)$, $\underline{t}'(\rho,a',b')=\underline{t}(\rho,a,b)-1$. On a clairement $\epsilon_{\underline{t},\underline{\eta}}(\rho,a,b)=\epsilon_{\underline{t}',\underline{\eta}'}(\rho,a',b')$. On consid\`ere le groupe $G'$ de m\^eme type que $G$ mais de rang celui de $G$ moins $d_{\rho}$ ($d_{\rho}$ est la dimension de $\rho$ vue comme repr\'esentation de $W_{F}$). Comme $\ell(\psi')=\ell(\psi)-2$, on sait d\'efinir $\pi(\psi',\underline{t}',\underline{\eta}')$ avec des propri\'et\'es sur les modules de Jacquet que l'on r\'ecrira en \ref{proprietejac}. On pose $A:=(a+b)/2-1$ et $B=(a-b)/2$. On a montr\'e en \cite{paquetdiscret} que l'induite
$$
<\rho\vert\,\vert^{\zeta B}, \cdots, \rho\vert\,\vert^{-\zeta A}>\times \pi(\psi',\underline{t}',\underline{\eta}')$$
a un unique sous-module irr\'eductible et  $\pi(\psi,\underline{t},\underline{\eta})$ est ce sous-module irr\'eductible, par d\'efinition. D'o\`u, en revenant aux d\'efinitions de $B,A$
$$
\pi(\psi,\underline{t},\underline{\eta})\hookrightarrow
<\rho\vert\,\vert^{(a-b)/2}, \cdots, \rho\vert\,\vert^{-\zeta ((a+b)/2-1)}>\times \pi(\psi',\underline{t}',\underline{\eta}') .
$$
\subsection{Cas d'un morphisme de bonne parit\'e}
On fixe $\psi$ et $(\rho,a,b)\in Jord(\psi)$; on dit que $(\rho,a,b)$ a bonne parit\'e si la repr\'esentation de $W_{F}\times SL(2,{\mathbb C})\times SL(2,{\mathbb C})$ d\'efinie par ce triplet est \`a valeurs dans un groupe de m\^eme type que $G^*$. On dit que $\psi$ a bonne parit\'e si tous les \'el\'ements de $Jord(\psi)$ ont bonne parit\'e. 

Fixons un ordre total sur $Jord(\psi)$; comme $Jord(\psi)$ a en g\'en\'eral de la multiplicit\'e, on pr\'ecise que 2 \'el\'ements $(\rho,a,b),(\rho',a',b')$ de $Jord(\psi)$ sont ordonn\'es m\^eme si $\rho=\rho', a=a', b=b'$; on note $\geq_{Jord(\psi)}$ cet ordre. On impose toujours \`a l'ordre mis de v\'erifier la propri\'et\'e {$(\mathcal{P})$}:

\begin{center}$(\rho',a',b')\geq_{Jord(\psi)} (\rho,a,b)$ 
si $\rho=\rho'$, $(a-b)(a'-b')>0$, $\vert a'-b'\vert>\vert a-b\vert$ et $a'+b'>a+b$. \end{center}

Soit $G_{>>}$ un groupe de m\^eme type que $G$ mais de rang plus grand et soit $\psi_{>>}$ un morphisme pour $G_{>>}$. On suppose que $Jord(\psi_{>>})$ est aussi muni d'un ordre total avec la propri\'et\'e ci-dessus et on dit que $\psi_{>>}$ domine $\psi$ si $\vert Jord(\psi_{>>})\vert=\vert Jord(\psi)\vert=:L$ et pour tout $i\in [1,L]$, en notant $(\rho_{>>,i},a_{>>,i},b_{>>,i})$ et $(\rho_{i},a_{i},b_{i})$ le i\`eme \'el\'ement de $Jord(\psi_{>>})$ et $Jord(\psi)$ respectivement, on a $\rho_{>>,i}=\rho_{i}$, $inf(a_{i},b_{i})=inf(a_{>>,i},b_{>>,i})$ et $(a_{>>,i}-b_{>>,i})(a_{i}-b_{i})\geq 0$.

Si $\psi$ a bonne parit\'e et seulement dans ce cas, il existe des morphismes $\psi_{>>}$ de restriction discr\`ete \`a la diagonale et dominant $\psi$. Pour d\'efinir $\Pi(\psi)$ le paquet de repr\'esentations associ\'ees \`a $\psi$ l'id\'ee est de l'obtenir par module de Jacquet \`a partir de $\Pi(\psi_{>>})$ o\`u $\psi_{>>}$ est un morphisme dominant $\psi$ et de restriction discr\`ete \`a la diagonale. L'application de $\Pi(\psi_{>>})$ sur $\Pi(\psi)$ est celle qui est expliqu\'ee en \ref{assouplissement} ci-dessous; comme on va le montrer en \ref{assouplissement} on obtient bien ainsi $\Pi(\psi)$ mais \'evidemment la param\'etrisation \`a l'int\'erieur de $\Pi(\psi)$ d\'epend du choix de $\psi_{>>}$; on a fait une construction qui ne d\'epend pas du choix tout simplement en imposant \`a l'ordre mis sur $Jord(\psi)$ plus de condition que  {$(\mathcal{P})$}. Pr\'ecis\'ement, on a demand\'e que 
\begin{center}
$\forall (\rho',a',b'), (\rho,a,b)\in Jord(\psi)
(\rho,a',b')\geq_{Jord(\psi)}(\rho,a,b)$ si $\rho=\rho'$ et $\vert (a'-b')\vert>\vert (a-b)\vert$ ou $\vert (a'-b')\vert = \vert (a-b)\vert$ et $a'+b'>a+b$ ou $(a',b')=(a,b)$ en tant qu'ensemble non ordonn\'e et $a'\geq b'$.
\end{center}
On renvoie \`a \cite{pourshahidi}.
On a alors montr\'e que la param\'etrisation de $\Pi(\psi)$ obtenue gr\^ace \`a celle de $\Pi(\psi_{>>})$ est ind\'ependante du choix de $\psi_{>>}$ dominant $\psi$ pour un ordre v\'erifiant cette propri\'et\'e; en \cite{general} on a m\^eme montr\'e que la param\'etrisation de $\Pi(\psi)$ qui en r\'esulte \`a l'aide de couples $(\underline{t},\underline{\eta})$ d\'efinis sur l'ensemble $Jord(\psi)\simeq Jord(\psi_{>>})$ peut se d\'efinir directement sur l'ensemble $Jord(\psi)$ vu comme ensemble ordinaire avec multiplicit\'e. Ce raffinement ne sert pas  pour ce travail et au contraire on pr\'ef\`erera garder une grande libert\'e sur le choix de l'ordre mis sur $Jord(\psi)$.

\subsection{D\'efinition g\'en\'erale de $\Pi(\psi)$}
Il reste encore \`a enlever l'hypoth\`ese que $\psi$ a bonne parit\'e. De l'alg\`ebre lin\'eaire \'el\'ementaire
montre alors que $\psi$ est la somme de fa\c{c}on unique $\psi_{mp}\oplus \psi_{bp}$ de 2 morphismes o\`u $\psi_{bp}$ a bonne parit\'e et aucun \'el\'ement de $Jord(\psi_{mp})$ a bonne parit\'e. On \'ecrit alors $\psi_{mp}=\psi_{1/2,mp}\oplus \psi_{-1/2,mp}$ de telle sorte que les repr\'esentations des groupes lin\'eaires associ\'ee $\pi^{GL}(\psi_{1/2,mp})$ et $\pi^{GL}(\psi_{-1/2,mp})$ se d\'eduisent l'une de l'autre par application de $\theta$; ici il n'y a pas unicit\'e mais cela n'a pas d'importance. On a alors montr\'e en \cite{pourshahidi} 3.2 que pour tout $
\pi\in \Pi(\psi_{bp})$ l'induite $\pi^{GL}(\psi_{1/2,mp})\times \pi$ est irr\'eductible et que cette tensorisation d\'efinit une bijection de $\Pi(\psi_{bp})$ sur $\Pi(\psi)$. Cela termine la description des paquets.
\subsection{Paquet d'Arthur et induction\label{induction}}
On fixe $\psi=:\psi_{mp}\oplus \psi_{bp}$ (cf. ci-dessus) et on suppose  que $Jord(\psi_{bp})$ a de la multiplicit\'e; on fixe $(\rho,a,b)\in Jord(\psi_{bp})$ intervenant avec au moins multiplicit\'e 2 et on note $\psi'$ le morphisme qui se d\'eduit de $\psi$ en enlevant 2 copies de $(\rho,a,b)$. L'hypoth\`ese de bonne parit\'e assure que $\psi'$ est encore un morphisme pour un groupe $G'$ de m\^eme type que $G$. On note $\pi(\rho,a,b)$ la repr\'esentation du groupe lin\'eaire $Speh( St(\rho,a),b)$ c'est-\`a-dire le sous-module de Langlands de l'induite $St(\rho,a)\vert\,\vert^{-(b-1)/2)}\times \cdots \times St(\rho,a)\vert\,\vert^{(b-1)/2)}$. 
\begin{prop}L'ensemble $\Pi(\psi)$ est l'ensemble des sous-modules irr\'eductibles des induites de la forme $\pi(\rho,a,b)\times \pi'$ o\`u $\pi'$ parcourt $\Pi(\psi')$.
\end{prop}
Cette proposition a \'et\'e d\'emontr\'ee en \cite{general} sans utiliser le fait que $\Pi(\psi')$ est form\'e de repr\'esentation unitaires. Ici, on en donne une d\'emonstration tr\`es rapide qui utilise l'unitarit\'e des \'el\'ements de $\Pi(\psi')$; on aura cette unitarit\'e d\`es que la totalit\'e des r\'esultats d'Arthur seront disponibles c'est-\`a-dire d\`es que l'on aura tous les lemmes fondamentaux y compris pond\'er\'es; sans l'unitarit\'e et avec la d\'emonstration ci-dessous on obtient la proposition en rempla\c{c}ant sous-modules par sous-quotients. On sait que $\pi^{GL}(\psi)=\pi(\rho,a,b)\times \pi^{GL}(\psi')\times \theta(\pi(\rho,a,b))$. Soit $T(\psi')$ la combinaison lin\'eaire stable de repr\'esentations de $G'$ dont la trace se transf\`ere en la $\theta$-trace de $\pi^{GL}(\psi')$.  La $\theta$-trace de $\pi^{GL}(\psi)$ est alors un transfert de la trace de l'induite $\pi(\rho,a,b)\times T(\psi')$ (d\'efinie dans le groupe de Grothendieck de $G$). Ainsi $\Pi(\psi)$ est form\'e des sous-quotients irr\'eductibles des induites de la forme $\pi(\rho,a,b)\times \pi'$ o\`u $\pi'\in \Pi(\psi')$; l'unitarit\'e admise de $\pi'$ assure qu'une telle induite est semi-simple d'o\`u le r\'esultat. Remarquons qu'en g\'en\'eral une telle induite est r\'eductible; on a montr\'e en \cite{general} qu'elle est sans multiplicit\'e de longueur inf\'erieure ou \'egale \`a $inf(a,b)+1$. Le fait qu'il n'y ait pas de multiplicit\'e se retrouve ici  parce que l'on sait que $\Pi(\psi)$ est sans multiplicit\'e (cf. \cite{pourshahidi}) mais cela utilise toute la force des r\'esultats d'Arthur alors que \cite{general} est \'el\'ementaire.
\subsection{Propri\'et\'es des modules de Jacquet\label{proprietejac}}
 On fixe $\psi$ d'o\`u $Jord(\psi)$ et $\pi\in \Pi(\psi)$. Soient $\rho$ une repr\'esentation cuspidale unitaire irr\'eductible, $\zeta=\pm$, $A,B$ des demi-entiers tels que $A-B\in \mathbb{Z}_{\geq 0}$. 
Le lemme ci-dessous est d\'ej\`a dans \cite{general} mais on pr\'ef\`ere en redonner la d\'emonstration pour pouvoir l'utiliser dans \ref{assouplissement} et rendre ce travail  ind\'ependant  de \cite{general}.
 \begin{lem}  On suppose que $B\neq 0$ et que  $Jac_{\zeta B, \cdots, \zeta A}\pi\neq 0$ alors il existe un sous ensemble de $Jord(\psi)$ de la forme $(\rho,A_{i},B_{i},\zeta)$ pour $i\in [1,v]$ avec $v$ un entier convenable tel que $B_{1}=B$, $A_{v}\geq A$ et pour tout $i\in ]1,v]$, $B_{i}\in ]B_{i-1},A_{i-1}+1]$.
 \end{lem}
On d\'emontre le lemme par la r\'ecurrence suivante. On fixe $\psi$; on a mis un ordre sur $Jord(\psi)$, pour d\'efinir $\pi$; on consid\`ere un \'el\'ement  $(\rho,A',B',\zeta')\in Jord(\psi)$ tel que pour tout $(\rho,A'',B'',\zeta'') \in Jord(\psi)$ plus grand que cet \'el\'ement, on a $B''>>A'$ et pour un 3e \'el\'ement de $Jord(\psi)$, $(\rho,A''',B''',\zeta''')$ strictement plus grand que $(\rho,A'',B'',\zeta'')$, on a $B'''>>A''$;  la notion de $>>$ est aussi relative aux $A,B$ du l'\'enonc\'e. Le plus grand \'el\'ement de $Jord(\psi)$ convient toujours;  si $(\rho,A',B',\zeta')$ peut \^etre le plus petit \'el\'ement de $Jord(\psi)$ alors $\psi$ est de restriction discr\`ete \`a la diagonale et tous les blocs de Jordan sauf \'eventuellement le plus petit sont de la forme $(\rho',A',B',\zeta')$ avec $B'>>A$. Il n'y a aucune difficult\'e \`a montrer que dans ce cas $Jac_{\zeta B, \cdots, \zeta A}\pi=0$ sauf \'eventuellement si pour $(\rho',A',B',\zeta')$ comme ci-dessus, $\rho'=\rho$, $\zeta'=\zeta$, $B'=B$ et $A'\geq A$.

On raisonne donc par r\'ecurrence sur la place d'un tel \'el\'ement $(\rho',A',B',\zeta')$ dans $Jord(\psi)$. On fixe un tel \'el\'ement dans $\psi$ et on suppose que ce n'est pas le plus petit \'el\'ement puisque ce cas a d\'ej\`a \'et\'e vu; on note $\psi'$ le morphisme qui se d\'eduit de $\psi$ en rempla\c{c}ant $(\rho,A',B',\zeta')$ par $(\rho,A'+T',B'+T',\zeta')$ avec $T'$ grand mais tel que si $(\rho,A'',B'',\zeta'')>_{Jord(\psi)}(\rho,A',B',\zeta')$, on a encore $B''>>A'+T'$. On peut alors appliquer le lemme par r\'ecurrence aux repr\'esentations dans le paquet associ\'e \`a $\psi'$. 

On note $S(\rho,A',B',T',\zeta')$ la repr\'esentation irr\'eductible associ\'ee aux multi-segments
$$
\begin{matrix}
\zeta' (B'+T' )&\cdots & \zeta'(A'+T')\\
\vdots &\cdots &\vdots\\
\zeta'(B'+1) &\cdots & \zeta'(A'+1)
\end{matrix}
$$
Par d\'efinition, pour toute repr\'esentation irr\'eductible $\pi$ du paquet associ\'e \`a $\psi$, il existe une repr\'esenta\-tion irr\'eductible $\pi'$ du paquet associ\'e \`a $\psi'$ et une inclusion
$$
\pi'\hookrightarrow S(\rho,A',B',T',\zeta')\times \pi.
$$
On suppose que $Jac_{\zeta B, \cdots, \zeta A}\pi\neq 0$. On fixe une repr\'esentation irr\'eductible $\sigma$ et une inclusion
$$
\pi \hookrightarrow \rho\vert\,\vert^{\zeta B} \times \cdots \times \rho\vert\,\vert^{\zeta A}\times \sigma.
$$
Avec une r\'ecurrence facile sur $A-B$, on peut supposer que cette inclusion se factorise en une inclusion 
$$
\pi\hookrightarrow <\rho\vert\,\vert^{\zeta B}, \cdots, \rho\vert\,\vert^{\zeta A}>\times \sigma.
$$
En remontant, on a une inclusion
$$
\pi'\hookrightarrow S(\rho,A',B',T',\zeta')\times <\rho\vert\,\vert^{\zeta B}, \cdots, \rho\vert\,\vert^{\zeta A}>\times \sigma.
$$
Supposons d'abord que l'induite $S(\rho,A',B',T',\zeta')\times <\rho\vert\,\vert^{\zeta B}, \cdots, \rho\vert\,\vert^{\zeta A}>$ dans le GL convenable est irr\'eductible. On peut \'echanger les 2 facteurs et on obtient donc par r\'eciprocit\'e de Frobenius que $Jac_{\zeta B, \cdots, \zeta A}\pi'\neq 0$. On applique le lemme par r\'ecurrence \`a $\pi'$ mais les blocs de Jordan qui interviennent sont n\'ecessairement comparables \`a $(\rho,A,B,\zeta)$ (qui n'est pas un bloc de Jordan, en g\'en\'eral) et sont donc des blocs de Jordan de $\psi$ et on obtient le lemme pour $\pi$. On sait que l'induite est irr\'eductible si $\zeta'\neq \zeta$: c'est une application imm\'ediate du fait que $\vert \zeta' (B'+1)-\zeta B\vert$ est le plus petit \'ecart entre une repr\'esentation intervenant dans $S(\rho,A',B',T',\zeta')$ et une repr\'esentation intervenant dans $\rho\vert\,\vert^{\zeta B}\times  \cdots \times  \rho\vert\,\vert^{\zeta A}$ et que ce nombre vaut $\vert (B'+1)-\zeta'\zeta B\vert$ ou encore $B'+1+B\geq 1+B>1$ car on a suppos\'e que $B\neq 0$.
On suppose donc que $\zeta'=\zeta$; si $\zeta=\zeta'=-$ on a encore l'irr\'eductibilit\'e d'apr\`es par exemple \cite{mw} 1.9 si $B\leq B'+1$ et $A\geq A'+1$. En utilisant l'involution de Zelevinsky, on a aussi le m\^eme r\'esultat si $\zeta=\zeta'=+$. En fait on sait que $Jac_{\zeta B}\pi=0$ s'il n'existe pas $(\rho,\tilde{A},B,\zeta)\in Jord(\psi)$ donc on peut se limiter au cas o\`u $B\leq B'$. Le seul cas qui reste est donc $B\leq B'$ et $A\leq A'$. Si $B=B'$, on prend $v=1$ et $(\rho,A_{1},B_{1},\zeta_{1})=(\rho,A',B',\zeta')$. Si $B<B'$, on factorise 
$$
<\rho\vert\,\vert^{\zeta B}, \cdots, \rho\vert\,\vert^{\zeta A}>\hookrightarrow <\rho\vert\,\vert^{\zeta B}, \cdots, \rho\vert\,\vert^{\zeta (B'-1)}>\times <\rho\vert\,\vert^{\zeta B'}, \cdots,\rho\vert\,\vert^{\zeta A'}>$$
et on remarque que l'induite $S(\rho,A',B',T',\zeta')\times <\rho\vert\,\vert^{\zeta B}, \cdots, \rho\vert\,\vert^{\zeta (B'-1)}>$ est irr\'eductible. On a donc $Jac_{\zeta B, \cdots, \zeta (B'-1)}\pi'\neq 0$. On applique alors le lemme pour $\pi'$ et on obtient le lemme pour $\pi$ en ajoutant l'\'el\'ement $(\rho,A',B',\zeta')$. Cela termine la preuve.

\subsection{Assouplissement des choix\label{assouplissement}}
Ici, on fixe $\psi$ et on veut d\'ecrire les repr\'esentations de $\Pi(\psi)$ avec des choix aussi peu contraignants que possible. On fixe un ordre total sur $Jord(\psi)$ ayant uniquement la propri\'et\'e:

\

soit  $(\rho,A,B,\zeta),(\rho',A',B',\zeta')\in Jord(\psi)$ tels que 
$\rho=\rho'$, $\zeta=\zeta'$,  $A>A'$ et $B>B'$ alors $(\rho,A,B,\zeta)>(\rho',A',B',\zeta')$.

\

Cette propri\'et\'e est la propri\'et\'e not\'ee $\mathcal P$ ci-dessus sauf qu'elle est exprim\'ee plus simplement quand on a remplac\'e $(\rho,a,b)$ par $(\rho,A,B,\zeta)$ puisque l'on a d\'ej\`a mis un choix en fixant $\zeta$ si $a=b$.

On construit $\psi_{>>}$ dominant $\psi$ pour cet ordre c'est-\`a-dire que pour tout $(\rho,A,B,\zeta) \in Jord(\psi)$, on fixe $T_{\rho,A,B,\zeta}$ de telle sorte que si $(\rho,A,B,\zeta)>(\rho',A',B',\zeta')$ alors $T_{\rho,A,B,\zeta}>> T_{\rho',A',B',\zeta'}$. Soient $\underline{t}, \underline{\eta}$ v\'erifiant \ref{notations}; ici on voit ces applications comme d\'efinies sur $Jord(\psi)\simeq Jord(\psi_{>>})$, c'est-\`a-dire qu'elles ne sont pas d\'efinies sur $Jord(\psi)$ vu comme ensemble sans multiplicit\'e. On construit $\pi(\psi_{>>},\underline{t},\underline{\eta})$. On calcule ensuite les modules de Jacquet de cette repr\'esentation $Jac_{\psi} \pi(\psi_{>>},\underline{t},\underline{\eta}):=$
$$
\circ_{(\rho,A,B,\zeta)\in Jord(\psi))} \circ_{j\in [1,T_{\rho,A,B,\zeta}]}Jac_{\zeta (B+j), \cdots, \zeta (A+j)}\pi(\psi_{>>},\underline{t},\underline{\eta}),
$$
o\`u les $(\rho,A,B,\zeta)$ sont pris dans l'ordre d\'ecroissant, c'est-\`a-dire, avec l'inversion provoqu\'ee par l'\'ecritu\-re de $\circ$, que l'on commence par ''baisser'' le plus petit \'el\'ement de $Jord(\psi_{>>})$ pour le faire devenir le plus petit \'el\'ement de $Jord(\psi)$ et et caetera.

\begin{prop} Avec les notations pr\'ec\'edentes, $Jac_{\psi} \pi(\psi_{>>},\underline{t},\underline{\eta})$ est nulle ou irr\'eductible; si cette repr\'esentation est non nulle elle est dans $\Pi(\psi)$ et toute repr\'esentation de $\Pi(\psi)$ est obtenue par cette proc\'edure pour un unique choix de $\underline{t}$ et $\underline{\eta}$. 
\end{prop}

On fixe $\psi$ et $\psi_{>>}$ et pour tout $k\in [1, \vert Jord(\psi)\vert]$, on note $(\rho_{k},A_{k},B_{k},\zeta_{k})$ le $k$ i\`eme \'el\'ement de $Jord(\psi)$, ces \'el\'ements \'etant ordonn\'es de telle sorte que le premier soit le plus petit. Et on note $(\rho_{k},A_{>>,k},B_{>>,k},\zeta_{k})$ le $k$ i\`eme bloc de Jordan de $\psi_{>>}$ avec le m\^eme principe pour l'ordre. Par hypoth\`ese $(\rho_{k},A_{>>,k},B_{>>,k},\zeta_{k})$ ''domine'' $(\rho_{k},A_{k},B_{k},\zeta_{k})$ ce qui veut dire que $T_{k}:=A_{>>,k}-A_{k}=B_{>>,k}-B_{k}>>0$. On note $\psi_{\leq k}$ le morphisme dont les blocs de Jordan sont tous les $(\rho_{i},A_{i},B_{i},\zeta_{i})$ pour $i\leq k$ et les $(\rho_{j},A_{>>,j},B_{>>,j},\zeta_{j})$ pour touts les $j>k$. On d\'emontre la proposition par r\'ecurrence sur $k$; pour s'autoriser $T_{1}\neq 0$, on commence \`a $k=0$, o\`u il n'y a rien \`a d\'emontrer. On le d\'emontre donc pour $k$ en l'admettant pour $k-1$. On note $S(\rho,A_{k},B_{k},T_{k},\zeta)$ la repr\'esentation associ\'ee aux multi-segments
$$
\begin{matrix}
B_{k}+T_{k} &\cdots &A_{k}+T_{k}\\
\vdots &\cdots &\vdots\\
B_{k}+1 &\cdots & A_{k}+1
\end{matrix}
$$
Le premier point \`a d\'emontrer est que si $\pi_{k}:=Jac_{\zeta (B_{k}+1),\cdots, \zeta (A_{k}+1)} \cdots Jac_{\zeta (B_{k}+T_{k}), \cdots, \zeta (A_{k}+T_{k})}\pi_{k-1}\neq 0$, alors il existe une repr\'esentation irr\'eductible $\sigma$, une inclusion
$$
\pi_{k-1}\hookrightarrow S(\rho,A_{k},B_{k},T_{k},\zeta)\times \sigma
$$et que $\pi_{k}=\sigma$. 

Pour tout $\ell\in [0,T_{k}[$, on d\'efinit la repr\'esentation irr\'eductible $S(\rho,A_{k}+\ell, B_{k}+\ell,T_{k},\zeta)$ en ne consid\'erant que les $T_{k}-\ell$ premi\`eres lignes du tableau ci-dessus. Et on d\'emontre  par r\'ecurrence descendante sur $\ell$ l'assertion suivante: il existe une repr\'esentation  $\sigma_{\ell}$ v\'erifiant $Jac_{\zeta C, \cdots, \zeta A'}\sigma_{\ell}=0$ pour tout $C \in [B_{k}+\ell+1, A_{k}+T_{k}]$ et $A\in [ A_{k}+\ell+1, A_{k}+T_{k}]$, avec $A\geq C$ et tel que l'on ait une inclusion
$$
\pi_{k-1}\hookrightarrow S(\rho,A_{k}+\ell,B_{k}+\ell, T_{k},\zeta)\times \sigma_{\ell}.$$
Remarquons tout de suite que comme $\pi_{k-1}$ est irr\'eductible, on peut supposer (comme nous le ferons) que $\sigma_{\ell}$ est irr\'eductible, on garde la propri\'et\'e de nullit\'e de certains modules de Jacquet de $\sigma_{\ell}$ et cette propri\'et\'e de nullit\'e assure que
$$
\circ_{i\in [T_{k},\ell+1]}Jac_{\zeta (B_{k}+i), \cdots, \zeta (A_{k}+i)}S(\rho,A_{k}+\ell,B_{k}+\ell, T_{k},\zeta)\times \sigma_{\ell}=\sigma_{\ell}.
$$
En particulier pour $\ell=0$, on aura le r\'esultat cherch\'e, soit $\pi_{k}=0$ soit $\pi_{k}=\sigma_{0}$. D\'emontrons donc l'assertion.

 Pour $\ell=T_{k}$, il n'y a rien presque rien \`a d\'emontrer en prenant $\sigma_{T_{k}}=\pi_{k-1}$; il suffit d'utiliser le fait que $Jac_{\zeta C}\pi_{k-1}=0$ s'il n'existe pas d'\'el\'ement de $Jord(\psi_{k-1})$ de la forme $(\rho,\tilde{A}, C,\zeta)$. Or cette non existence est vraie par hypoth\`ese puisque $C\in [B_{k}+T_{k}+1,A_{k}+T_{k}]$ . On suppose l'assertion vraie pour $\ell+1\leq T_{k}$ et on va la d\'emontrer pour $\ell$. On a une inclusion
$$
\pi_{k-1}\hookrightarrow S(\rho,A_{k}+\ell+1, B_{k}+\ell+1,T_{k},\zeta)\times \sigma_{\ell+1}.
$$
On v\'erifie alors
$$
Jac_{\zeta (B_{k}+\ell), \cdots, \zeta (A_{k}+\ell)}\circ_{i\in [T_{k},\ell+1]}Jac_{\zeta (B_{k}+i), \cdots, \zeta (A_{k}+i)} \pi_{k-1}\hookrightarrow Jac_{\zeta (B_{k}+\ell), \cdots, \zeta (A_{k}+\ell)} \sigma_{\ell+1}.
$$
Comme $\pi_{k}\neq 0$, on a une non nullit\'e du membre de droite. Ainsi il existe une repr\'esentation irr\'eductible $\sigma_{\ell}$ et une inclusion
$$
\sigma_{\ell+1}\hookrightarrow \rho\vert\,\vert^{\zeta (B_{k}+\ell)}\times \cdots \times \rho\vert\,\vert^{\zeta (A_{k}+\ell)}\times \sigma_{\ell}.
$$
Comme $Jac_{\zeta C, \cdots, \zeta (A_{k}+\ell)}\sigma_{\ell+1}=0$ par hypoth\`ese, cette inclusion se factorise par l'unique sous-module irr\'eductible (pour le GL convenable) de l'induite $\rho\vert\,\vert^{\zeta (B_{k}+\ell)}\times \cdots \times \rho\vert\,\vert^{\zeta (A_{k}+\ell)}$, d'o\`u
$$
\sigma_{\ell+1}\hookrightarrow <\rho\vert\,\vert^{\zeta (B_{k}+1)}, \cdots, \rho\vert\,\vert^{\zeta (A_{k}+1)}>\times \sigma_{\ell}.
$$
En remontant, on obtient une inclusion
$$
\pi_{k-1}\hookrightarrow S(\rho,A_{k}+\ell+1, B_{k}+\ell+1,T_{k},\zeta)\times <\rho\vert\,\vert^{\zeta (B_{k}+\ell)}, \cdots, \rho\vert\,\vert^{\zeta (A_{k}+\ell)}>\times \sigma_{\ell}.
$$
On veut encore d\'emontrer que cette inclusion se factorise par l'unique sous-module irr\'eductible (pour un GL convenable) de $S(\rho,A_{k}+\ell+1, B_{k}+\ell+1,T_{k},\zeta)\times <\rho\vert\,\vert^{\zeta (B_{k}+\ell)}, \cdots, \rho\vert\,\vert^{\zeta (A_{k}+\ell)}>$ qui est pr\'ecis\'ement $S(\rho,A_{k}+\ell,B_{k}+\ell,T_{k},\zeta)$. S'il n'en est pas ainsi, l'inclusion se factorise par 
$$
\rho\vert\,\vert^{\zeta (B_{k}+\ell)}\times S(\rho,A_{k}+\ell+1, B_{k}+\ell+1,T_{k},\zeta)\times <\rho\vert\,\vert^{\zeta (B_{k}+\ell+1)}, \cdots, \rho\vert\,\vert^{\zeta (A_{k}+\ell)}>\times \sigma_{\ell}.
$$
Or l'induite $S(\rho,A_{k}+\ell+1, B_{k}+\ell+1,T_{k},\zeta)\times <\rho\vert\,\vert^{\zeta (B_{k}+\ell+1)}, \cdots, \rho\vert\,\vert^{\zeta (A_{k}+\ell)}>$ est irr\'eductible (cf. \cite{mw} 1.9); on peut donc encore \'echanger les 2 facteurs; on obtient alors une inclusion pour $\tau$ une repr\'esentation convenable:
$$
\pi_{k-1}\hookrightarrow \rho\vert\,\vert^{\zeta (B_{k}+\ell)}\times <\rho\vert\,\vert^{(B_{k}+\ell+1)}, \cdots, \rho\vert\,\vert^{A_{k}+\ell}>\times \tau.
$$
On applique le lemme \ref{proprietejac} \`a $\pi_{k-1}$ et $\psi_{k-1}$ pour $A=A_{k}+\ell$ et $B=B_{k}+\ell$. Comme $\ell <T_{k}$, l'\'el\'ement de $Jord(\psi_{k-1})$ not\'e $(\rho,A_{1},B_{1},\zeta_{1})$ de loc.cit est de la forme $(\rho,A',B',\zeta)$ v\'erifiant $(\rho,A',B',\zeta)<_{Jord(\psi_{k-1}} (\rho,A_{k}+T_{k},B_{k}+T_{k},\zeta)$; on a donc aussi $(\rho,A',B',\zeta)<_{Jord(\psi)}(\rho,A_{k},B_{k},\zeta)$ par choix de l'ordre. Ceci reste vrai pour tout $i\in [1,v]$ puisque les blocs se chevauchent et en particulier $B_{v},A_{v}$ (avec les notations de loc. cit) v\'erifie soit $B_{v}\leq B_{k}$ soit $A_{v}\leq A_{k}$; or on a $B_{1}=B_{k}+\ell>B_{k}$ d'o\`u aussi $B_{v}>B_{k}$ donc $A_{v}\leq A_{k}$ ce qui est contradictoire avec $A_{v}\geq A_{k}+\ell$. D'o\`u l'assertion. On a donc montr\'e l'inclusion
$$
\pi_{k-1}\hookrightarrow S(\rho,A_{k}+\ell,B_{k}+\ell,T_{k},\zeta)\times \sigma_{\ell}.
$$
Il reste \`a montrer que pour tout $C\geq B_{k}+\ell$ et tout $A\in [A_{k}+\ell,A_{k}+T_{k}]$, avec $A\geq C$, on a
$$
Jac_{\zeta C, \cdots, \zeta A}\sigma_{\ell}=0.
$$
Raisonnons par l'absurde et on fixe $C$ maximum avec une telle propri\'et\'e, d'o\`u l'existence d'une repr\'esentation irr\'eductible $\sigma'$ et une inclusion
$$
\sigma_{\ell}\hookrightarrow <\rho\vert\,\vert^{\zeta C}, \cdots, \rho\vert\,\vert^{\zeta A}>\times \sigma'.
$$
La repr\'esentation $S(\rho,A_{k}+\ell,B_{k}+\ell,T_{k},\zeta)\times <\rho\vert\,\vert^{\zeta C}, \cdots, \rho\vert\,\vert^{\zeta A}>$ est irr\'eductible \`a cause des hypoth\`eses sur $C\geq B_{k}+\ell$ et $A\leq A_{k}+T_{k}$. En remontant on trouve encore une inclusion de $\pi_{k-1}$ dans une induite $<\rho\vert\,\vert^{\zeta C}, \cdots, \rho\vert\,\vert^{\zeta A}>\times \sigma''$ avec $\sigma''$ convenable et on peut conclure comme ci-dessus sauf dans le cas o\`u $C=B_{k}+T_{k}$ qui est tout \`a fait possible. Ici on aurait que $Jac_{\zeta (B_{k}+T_{k}), \zeta (B_{k}+T_{k})}\pi_{k-1}\neq 0$ ce qui est exclu par exemple par \cite{paquetdiscret} 5.2 (c'est facile). 

On a donc ainsi d\'efini $\pi_{k}$ comme repr\'esentation irr\'eductible. On remarque aussi que $\pi_{k}$ d\'etermine uniquement $\pi_{k-1}$ car l'induite $S(\rho,A_{k},B_{k},T_{k},\zeta)\times \pi_{k}$ a un unique sous-module irr\'eductible par r\'eciprocit\'e de Frobenius. On a donc d\'emontr\'e le premier point de la proposition.

\vskip 0.5cm

Montrons que les repr\'esentations obtenues sont dans $\Pi(\psi)$ encore par r\'ecurrence sur $k$ puisque l'on conna\^{\i}t le r\'esultat dans le cas de $\psi_{>>}$.
On sait donc qu'une bonne combinaison lin\'eaire des $\pi_{k-1}$ \`a coefficients tous non nuls se tranf\`ere en la trace tordue de $\pi^{GL}(\psi_{k})$; cela donne une identit\'e de caract\`ere. Etant donn\'e la compatibilit\'e du transfert \`a la restriction partielle (\cite{selecta} 4.2) la m\^eme combinaison lin\'eaire de $\pi_{k}$ se transf\`ere en $Jac_{\psi}^\theta\pi^{GL}(\psi_{k-1})$.  L'ind\'ependance lin\'eaire des caract\`eres du groupe $G$ donnera le r\'esultat cherch\'e d\`es que l'on aura montr\'e que $Jac_{\psi}^\theta \pi^{GL} (\psi_{k-1})=\pi^{GL}(\psi_{k})$.

Conceptuellement, on d\'efinit $Jac_{\psi}^\theta\pi^{GL}(\psi_{>>})$ comme on a d\'efini $Jac_{\psi}$ et le point est que$$
Jac_{\psi}^\theta\pi^{GL}(\psi_{>>})=\pi^{GL}(\psi).
$$
En effet, soit $(\rho',A',B',\zeta')$ un quadruplet comme ceux que l'on consid\`ere ici. On consid\`ere le tableau suivant:
$$\begin{array}{ccc} \zeta' B' & \cdots & -\zeta'  A'\\ \vdots & \vdots & \vdots \\ \zeta ' A' & \cdots & -\zeta ' B'\end{array}
$$
on peut lui associer gr\^ace aux travaux de Zelevinsky une repr\'esentation irr\'eductible bas\'ee sur la cuspidale $\rho'$; on note $Z(\rho',A',B',\zeta')$ cette repr\'esentation.

Pour une repr\'esentation $\sigma$ d'un groupe lin\'eaire, on note $Jac^{g}_{\rho\vert\,\vert^x}\sigma$ l'analogue de $Jac_{\rho\vert\,\vert^x}$ d\'efini en \ref{notations} et $Jac^{d}_{\rho\vert\,\vert^{x}}\sigma$ l'object de m\^eme nature quand on \'echange la gauche et la droite, ce qui a un sens pour les groupes lin\'eaires et n'en avait pas pour les groupes classiques. Evidemment cela suppose que l'on a fix\'e les matrices triangulaires comme sous-groupe de Borel  et que tous les paraboliques contiennent ce sous-groupe de  Borel. Avec ces notations $Jac^\theta_{\rho\vert\,\vert^x} \sigma= Jac^{g}_{\rho\vert\,\vert^{x}}Jac^{d}_{\rho\vert\,\vert^{-x}}\sigma=Jac^{d}_{\rho\vert\,\vert^{-x}}Jac^{g}_{\rho\vert\,\vert^x}\sigma$.

Fixons $\rho$; soit $[x,y]$ un segment croissant ou d\'ecroissant $x,y$ \'etant des demi-entiers. Les repr\'esentations $Z(\rho',A',B',\zeta')$ v\'erifient $Jac_{\rho\vert\,\vert^{x}, \cdots, \rho\vert\,\vert^{y}}Z(\rho',A',B',\zeta')\neq 0$ exactement quand $\rho'=\rho$, $x=\zeta' B'$ et $y\in [\zeta' B', -\zeta' A']$ ou $y\in [\zeta' B', \zeta' A']$, les 2 possibilit\'es sont exclusives. Dans le premier cas le r\'esultat est l'unique repr\'esentation irr\'eductible associ\'ee au tableau
$$\begin{array}{ccccc}& &y-\zeta'1& \cdots & -\zeta'  A'\\ \zeta ' (B'+1)&\cdots & \cdots &\cdots & -\zeta' (A'-1)\\ \vdots &
\vdots&\vdots & \vdots& \vdots \\ \zeta ' A' &\cdots & \cdots & \cdots & -\zeta ' B'\end{array}
$$
Et dans le deuxi\`eme cas, c'est l'unique repr\'esentation irr\'eductible associ\'ee au tableau
$$
\begin{array}{cccc}& \zeta' (B'-1)& \cdots & -\zeta'  A'\\ & \vdots &\vdots & \vdots \\
y+\zeta' 1&\vdots & \vdots& \vdots \\ \vdots & \vdots & \vdots & \vdots \\
\zeta' A' &\cdots &\cdots & -\zeta ' B'\end{array}
$$
On montre le r\'esultat suivant: soient $\psi'$ un morphisme et $(\rho,A,B,\zeta)\in Jord(\psi')$; on suppose que $Jord(\psi')$ est muni d'un ordre total v\'erifiant:

 pour tout $(\rho',A',B',\zeta')$ tel que $(\rho',A',B',\zeta)> (\rho,A,B,\zeta)$, $B'>>A$;
 
 pour tout $(\rho',A',B',\zeta')<(\rho,A,B,\zeta)$ soit $\rho'\not\simeq \rho$ soit $\zeta'\neq \zeta$ avec $B'\neq 0$ soit $B'\leq B$ soit $A'\leq A$.
 
 On note $\psi'_{>}$ le morphisme qui se d\'eduit de $\psi$ en changeant $(\rho,A,B,\zeta)$ en $(\rho, A+1,B+1,\zeta)$ et on va montrer que $Jac^\theta_{\zeta (B+1), \cdots, \zeta (A+1)}\pi^{GL}(\psi'_{>})=\pi^{GL}(\psi')$. De proche en proche, cela donnera l'assertion cherch\'ee.
 
 On \'ecrit
 $
 \pi^{GL}(\psi'_{>})=\times _{(\rho',A',B',\zeta') \in Jord(\psi)-\{(\rho,A,B,\zeta)\}}Z(\rho',A',B',\zeta') \times Z(\rho,A+1,B+1,\zeta)$.
 On calcule $Jac^g_{\zeta (B+1), \cdots, \zeta (A+1)}\pi^{GL}(\psi'_{>})$. Les formules standard de Bernstein Zelevinsky, disent que le r\'esultat a une filtration dont les sous-quotients sont index\'es par les d\'ecoupages de l'ensemble $\mathcal{E}:=\{\zeta (B+1), \cdots, \zeta (A+1)\}$ en $\vert Jord(\psi'_{>})\vert$ sous-ensembles, $\mathcal{E}_{(\rho',A',B',\zeta')}$, le sous-quotient correspondant \`a ce d\'ecoupage \'etant isomorphe \`a 
 $
 \times_{(\rho',A',B',\zeta')\in Jord(\psi'_{>})} Jac_{x\in \mathcal{E}_{(\rho',A',B',\zeta'}}Z(\rho',A',B',\zeta')$.
 D'apr\`es ce que l'on a vu, $\mathcal{E}
 _{(\rho',A',B',\zeta')}=\emptyset$ si $\zeta'B'\notin [\zeta(B+1),\cdots,\zeta (A+1)]$. Ainsi $\mathcal{E}_{(\rho',A',B',\zeta')}\neq\emptyset$ entra\^{\i}ne que n\'ecessairement $(\rho',A',B',\zeta')\leq (\rho,A,B,\zeta)$ et que $\rho=\rho'$, $B'>B$ et $\zeta'=\zeta$; supposons que $(\rho',A',B',\zeta')\neq (\rho,A+1,B+1,\zeta)$,  les propri\'et\'es de l'ordre assurent alors que  $A'\leq A$ et cela entra\^{\i}ne que $\zeta (A+1) \notin \mathcal{E}_{(\rho',A',B',\zeta')}$. Ainsi $\zeta (A+1)\in \mathcal{E}_{(\rho,A+1,B+1,\zeta)}$ mais cela force $\zeta (B+1)\in \mathcal{E}_{(\rho,A+1,B+1,\zeta)}$ et $\mathcal{E}=\mathcal{E}_{(\rho,A+1,B+1,\zeta)}$. Ainsi qu'il n'y a qu'un seul d\'ecoupage possible; le m\^eme argument s'applique pour $Jac^{d}_{-\zeta (A+1), \cdots, -\zeta (B+1)}$ et on trouve
 $$
 Jac^\theta_{\zeta (B+1), \cdots, \zeta (A+1)}\pi^{GL}(\psi'_{>})=\pi^{GL}(\psi').
 $$
 On note $\sum_{\pi_{k-1}\in \Pi(\psi_{k-1})}a(\pi_{k-1})\pi_{k-1}$ la combinaison lin\'eaire stable qui se transf\`ere en la trace tordue de $\pi^{GL}(\psi_{k-1})$ pour des bons coefficients $a(\pi_{k-1})\neq 0$. On a alors montr\'e que la combinaison $\sum_{\pi_{k-1}\in \Pi(\psi_{k-1}} a(\pi_{k-1})Jac_{\psi_{k}}\pi_{k-1}$ est stable et se tranf\`ere en la trace tordue de $\pi^{GL}(\psi_{k})$. On rappelle que les $Jac_{\psi_{k}}\pi_{k-1}$ sont soit nuls soit irr\'eductibles et ceux qui sont non nuls sont tous distincts. C'est alors l'ind\'ependance lin\'eaire des caract\`eres qui permet de conclure.

\section{D\'efinition de la normalisation et \'enonc\'e du th\'eor\`eme}
\subsection{Normalisation\label{normalisation}}
On fixe $\Psi$ et $\pi\in \Pi(\Psi)$; on fixe aussi un entiers $a_{0}$ et une repr\'esentation cuspidale autoduale $\rho$. On consid\`ere l'op\'erateur d'entrelacement standard
$$M(s):=St(\rho,a_{0})\vert\,\vert^{s}\times \pi \rightarrow St(\rho,a_{0})\vert\,\vert^{-s}\times \pi$$ au voisinage de $s \in \mathbb{R}_{\geq 0}$. On d\'efinit la fonction m\'eromorphe de $s$ suivante
$$
r(s,\psi):={ \times}_{(\tau,a,b)\in Jord(\Psi)}\frac{L(St(\rho,a_{0})\times St(\tau,a),s-(b-1)/2)}{L(St(\rho,a_{0})\times St(\tau,a),s+(b+1)/2)}$$
$$ \times \frac{L(\rho,r_{G},2s)}{L(\rho,r_{G},2s+1)},
$$o\`u $r_{G}$ est d\'efini par la repr\'esentation naturelle de $G^*$ qui sert \`a d\'efinir le transfert.
Et on pose $N_{\psi}(s, \pi):= M(s,\pi)r(s,\psi)^{-1}$ et en g\'en\'eral $\psi$ est fix\'e et on pose donc $N(s,\pi):=N_{\psi}(s,\pi)$. On conna\^{\i}t compl\`etement explicitement les fonctions $L$ qui interviennent ci-dessus. On rappelle les r\'esultats de Shahidi \cite{shahidi} que $$L(St(\rho,a_{0})\times St(\rho',a),s)=\times_{k\in [\vert (a-a')/2\vert, (a+a')/2[}L(\rho\times \rho',s+k).$$
On sait aussi qu'il existe une factorisation $L(\rho\times \rho,s)=L(\rho,r_{G},s)L(\rho,r'_{G},s)$ o\`u $r_{G}\oplus r_{G}'$ est la d\'ecomposition de la repr\'esentation $\otimes^2 {\mathbb C}^n$ o\`u ${\mathbb C}^n$ est la repr\'esentation naturelle de $GL(n,{\mathbb C})$ ($n=m^*_{G}$) la d\'ecomposition d\'epend du plongement de $G^*$ dans $GL(n,{\mathbb C})$.

On en d\'eduit facilement que les d\'enominateurs ne peuvent avoir de p\^oles au voisinage d'un nombre r\'eel positif ou nul et les num\'erateurs n'ont pas de z\'eros. Ainsi au voisinage d'un nombre r\'eel positif ou nul, la fonction $r(s,\psi)$ est d'ordre inf\'erieur ou \'egal \`a 0.

On a  vu en \cite{seriedeisenstein} 2.1  que $r(s,\psi)r(-s,\psi)$ diff\`ere du produit $M(-s)\circ M(s)$ vu comme op\'erateur scalaire par une fonction holomorphe inversible au voisinage de tout point r\'eel. En loc. cite on a bien exliqu\'e que cette normalisation n'est pas en g\'en\'eral la normalisation de Langlands-Shahidi; c'est la normalisation de Langlands-Shahidi uniquement si $\pi$ est dans le paquet de Langlands associ\'e au paquet d'Arthur. 
\subsection{Enonc\'e du th\'eor\`eme et sh\'ema de la preuve\label{shemadelapreuve}}
\begin{thm} L'op\'erateur $N_{\psi}(s,\pi)$ est holomorphe en tout point $s$ r\'eel positif ou nul.
\end{thm}
Le choix de la normalisation est de nature globale et les cons\'equences attendues comme en \cite{seriedeisenstein} sont aussi de nature globale, mais je n'ai pas de d\'emonstration globale de ce r\'esultat car tout essai s'est heurt\'e au fait qu'un tel op\'erateur peut \^etre identiquement z\'ero et que pour obtenir des r\'esultats en la place $v$ il faut au moins emp\^echer des 0 aux autres places; or on a peu de libert\'e car la normalisation d\'epend du paquet, $\psi$. Certes il y a plusieurs fa\c{c}on de globaliser ce paquet mais la repr\'esentation de la deuxi\`eme copie de $SL(2,{\mathbb C})$ est elle de nature globale et est donc fix\'ee par la place $v$ et ne donne lieu \`a aucune libert\'e.

Le cas $s=0$: ici on sait que $N(-s,\pi)\circ N(s,\pi)$ est une fonction holomorphe inversible, en particulier n'a pas de p\^ole en $s=0$. Donc \`a une fonction holomorphe inversible pr\`es $N(-s,\pi)$ est dual de $N(s,\pi)$ pour la dualit\'e naturelle entre $St(\rho,a_{0})\vert\,\vert^{s}\times \pi$ et $St(\rho,a_{0})\vert\,\vert^{-s}\times \pi$ construite gr\^ace \`a l'unitarit\'e de $\pi$. Cela suffit pour avoir l'holomorphie en $s=0$ et m\^eme le fait que $N(0,\pi)$ est bijectif.

On se place donc au voisinage d'un r\'eel $s=s_{0}>0$. On pose $b_{0}=2s_{0}-1$. Dans ce paragraphe, on explique la d\'emonstration. Cette d\'emonstration se  fait par r\'eduction au cas o\`u $\pi$ est une repr\'esentation temp\'er\'ee. Supposons momentan\'ement que $\pi$ soit temp\'er\'ee; sous cette hypoth\`ese, on sait gr\^ace aux travaux d'Harish-Chandra que l'op\'erateur d'entrelacement non normalis\'e est holomorphe en $s=s_{0}$ par positivit\'e stricte. Il suffit alors de v\'erifier que $r(s,\psi)$ n'a pas de z\'ero en $s=s_{0}$; on a d\'ej\`a vu que $r(s,\psi)$ est d'ordre inf\'erieur ou \'egal \`a z\'ero au voisinage de tout $s$ r\'eel strictement positif. Cela donne le r\'esultat.

On traitera aussi facilement le cas o\`u $Jord(\psi)$ a de la multiplicit\'e cf. \ref{casinduit}.

La suite de la d\'emonstration du th\'eor\`eme se fait par r\'ecurrence sur les entiers suivants; soit $\zeta$ un signe. On pose $\ell(\psi; \zeta):=\sum_{inf(a,b); \zeta (a-b) \geq  0}(inf(a,b)-1)$. On note $\zeta_{0}$ le signe de $a_{0}-b_{0}$ en prenant $+$ si $a_{0}=b_{0}$. Si $\zeta_{0}=-$ on fait  une r\'ecurrence d'abord  sur $n(\psi,-):=\vert \{(\rho,a,b)\in Jord(\psi); a<b\}$ puis sur $\ell(\psi,-\zeta_{0})$. Le d\'ebut de la r\'ecurrence se fait donc quand $\ell(\psi;\zeta_{0})=n(\psi,-)=0$;  on a donc pour tout $(\rho,a,b)\in Jord(\psi)$, $a\leq b$ par $n(\psi,-)=0$ et $b=1$ par $\ell(\psi;\zeta_{0})=0$. Ainsi, dans ce cas, $\psi$ est temp\'er\'e et le d\'ebut de la r\'ecurrence  a d\'ej\`a \'et\'e montr\'ee. Supposons maintenant que $\zeta_{0}=+$; ici  la r\'ecurrence se fait sur $\ell(\psi,+)+\ell(\psi,-)$. Le d\'ebut de la r\'ecurrence est donc quand pour tout $(\rho,a,b)\in Jord(\psi)$, $inf(a,b)=1$ et ce cas sera trait\'e en \ref{holomorphiecaselementaire} on le fera en utilisant la description des repr\'esentations obtenue dans  \cite{elementaire} et rappel\'ee ici.

Un cas facile est le cas o\`u l'on peut construire un morphisme 
$\psi'$ tel que $\pi\in \Pi(\psi)\cap \Pi(\psi')$ et tel que l'on puisse appliquer la r\'ecurrence \`a $\psi'$; le seul point alors est de v\'erifier que $r(s,\psi')/r(s,\psi)$ est holomorphe. C'est d'ailleurs dans ce cas que l'on trouve facilement des exemples o\`u $N(s,\psi)$ est identiquement 0 en $s=s_{0}$.

Le cas le plus fr\'equent est celui o\`u on construit un morphisme $\psi'$ est relatif \`a un groupe de rang plus petit que le groupe de d\'epart et auquel on peut appliquer l'hypoth\`ese de r\'ecurrence. Ce que l'on montre alors est qu'il existe une repr\'esentation $\pi'\in \Pi(\psi')$ et une repr\'esentation compl\`etement explicite $\sigma$ d'un groupe lin\'eaire avec une inclusion
$$
\pi\hookrightarrow \sigma\times \pi'.
$$
Ici il faut donc aussi consid\'erer les op\'erateurs d'entrelacements $$St(\rho,a_{0})\vert\,\vert^{s_{0}}\times \sigma \rightarrow \sigma \times St(\rho,a_{0}\vert\,\vert^{s_{0}}$$
$$
\sigma\times St(\rho,a_{0})\vert\,\vert^{-s_{0}}\rightarrow St(\rho,a_{0})\vert\,\vert^{-s_{0}}\times \sigma.
$$
Donc il faut aussi \'etudier ces op\'erateurs d'entrelacements (ce qui a \'et\'e fait en \cite{mw}) et on est presque oblig\'e de ne consid\'erer que des cas o\`u apr\`es une normalisation explicite ils sont holomorphes; ensuite on compare les normalisations ce qui est facile. Comme les induites \'ecrites ne sont pas irr\'eductibles en g\'en\'eral, dans cette \'etape on perd des informations. Le plus ennuyeux est qu'avec cette m\'ethode, on ne contr\^ole pas la non nullit\'e des op\'erateurs. Toutefois, une fois l'holomorphie d\'emontr\'ee, on pourra utiliser des op\'erateurs non holomorphes en contr\^olant l'ordre des p\^oles, ce qui, je l'esp\`ere, permettra d'aller plus loin.

On explique  les \'etapes de la r\'ecurrence et la construction des $\psi'$; on suppose d'abord que $\ell(\psi;-\zeta_{0})>0$ et on  baisse $\ell(\psi;-\zeta_{0})$ sans modifier $n(\psi,\zeta_{0})$ (d\'efini si $\zeta_{0}=-$) ni $\ell(\psi,\zeta_{0})$. On consid\`ere alors un ordre sur $Jord(\psi)$ qui est tel que $(\rho,a,b)>_{Jord(\psi)}(\rho,a',b')$ si $\zeta_{0}(a-b)\leq 0$ et $\zeta_{0}(a'-b')>0$. On remplace chaque $(\rho,a,b)$ tel que $\zeta_{0}(a-b)\leq 0$ par $(\rho, A=(a+b)/2-1, B=\vert (a-b)/2\vert,-\zeta_{0})$ et on fixe un tel \'el\'ement tel que $A-B>0$ et $B$ est maximal avec cette propri\'et\'e; on met alors un ordre sur $Jord(\psi)$ tel que $(\rho,A',B',\zeta')>_{Jord(\psi)} $ entra\^{\i}ne que $\zeta'=-\zeta_{0}$ $B'>B$ et $A'>A$. Ainsi pour $(\rho,A',B',\zeta')>_{Jord(\psi)}(\rho,A,B,\zeta)$, on aura n\'ecessairement $A'=B'$. Un cas est alors facile, celui o\`u pour tout $(\rho,A',B',\zeta)>_{Jord(\psi)}(\rho,A,B,-\zeta_{0})$, on a$B'\geq A$; dans ce cas, on montre (cf. \ref{isole}) que l'on peut remplacer $\psi$ par un morphisme $\psi'$ qui se d\'eduit de $\psi$ soit en rempla\c{c}ant $(\rho,A,B,-\zeta_{0})$ par $\cup_{C\in [B,A]}(\rho,C,C,-\zeta_{0})$ soit par $(\rho,A-1,B+1,-\zeta_{0})$ ce qui  a remplac\'e $\ell(\psi;-\zeta_{0})$ par $\ell(\psi;-\zeta_{0})-(A-B)$ dans le premier cas et  par $\ell(\psi;-\zeta_{0})-2$ dans le deuxi\`eme cas.  Dans le cas contraire trait\'e en \ref{casnonisole}, on ordonne les demi-entiers $B'$ tels que $B'\in ]B,A[$ et $(\rho,B',B',-\zeta_{0})\in Jord(\psi)$; on note $B_{inf+1}$ le plus petit et $B_{max}$ le plus grand, en posant $B_{inf}=B$. On se ram\`ene d'abord au cas o\`u pour $B'> B''$ cons\'ecutifs dans cet ensemble auquel on a ajout\'e $B$, on a $Jac_{\zeta B', \cdots,\zeta B'' }\pi=0$; cette r\'eduction n'am\'eliore par la r\'ecurrence. Mais avec cette hypoth\`ese, on montre que l'on peut remplacer $\psi$ par un morphisme $\psi'$ qui se d\'eduit de $\psi$ en rempla\c{c}ant $(\rho,A,B,-\zeta_{0})$ soit par $(\rho,A-1,B+1,-\zeta_{0})$ soit par $(\rho,B_{max},B,-\zeta_{0})\cup_{C\in [B_{max},A]}(\rho,C,C,-\zeta_{0})$; dans le premier cas, on a remplac\'e $\ell(\psi;\zeta_{0})$ par $\ell(\psi;-\zeta_{0})-2$ et dans le deuxi\`eme cas par $\ell(\psi;-\zeta_{0})-(A-B_{max})$. Les cons\'equences sur l'holomorphie cherch\'ee sont  en \ref{holomorphiecasnonisolesignecontraire}. On remarque pour la suite que dans toutes ces \'etapes on n'a pas modifi\'e $n(\psi,\zeta_{0})$.

Ainsi, on est ramen\'e au cas o\`u $\ell(\psi;-\zeta_{0})=0$ sans avoir modifi\'e $n(\psi,\zeta_{0})$ ni $\ell(\psi,\zeta_{0})$. 

Dans le cas o\`u $\zeta_{0}=+$, on baisse ensuite $\ell(\psi,\zeta_{0})$ en \ref{holomorphiecasnonisolememesigne}

Il faut donc encore expliquer comment on baisse $n(\psi,\zeta_{0})$ ce qui sera utile pour $\zeta_{0}=-$  mais dans ce cas on modifie $\ell(\psi;-\zeta_{0})$ qui peut redevenir $>0$ d'o\`u l'ordre de la r\'ecurrence.  On commence par partir d'un $\psi$ tel que $\ell(\psi,-\zeta_{0})=0$ ce qui est loisible. On fixe un \'el\'ement $(\rho,A,B,\zeta_{0})$ dans $Jord(\psi)$ tel que $B$ soit minimal; on suppose \'evidemment qu'il en existe. On met un ordre sur $Jord(\psi)$ tel que cet \'el\'ement soit le plus petit \'el\'ement de $Jord(\psi)$; on montre d'abord que l'on peut se ramener au cas o\`u $B=0$ ou $1/2$ au prix \'eventuellement d'une modification des $(\rho,A',B',-\zeta_{0})$; puis presque par d\'efinition que l'on peut alors changer $-\zeta_{0}$ en $\zeta_{0}$. Dans cet \'etape on a donc diminu\'e $n(\psi,\zeta_{0})$ de $1$ \'eventuellement en augmentant $\ell(\psi,-\zeta_{0})$. Il faut remarquer que pour mener \`a bien cette preuve, on a pris un ordre sur $Jord(\psi)$ tel que $(\rho,A,B,\zeta_{0})$ est le plus petit \'el\'ement alors que tout autre \'el\'ement de la forme $(\rho,A',B',\zeta_{0})$ est plus grand que tout \'el\'ement de la forme $(\rho,A'',B'',-\zeta_{0})$; cf \ref{reductionparlebas} et \ref{finpreuve-}.
\subsection{Propri\'et\'es des facteurs de normalisation\label{calculdelordre}}
On fixe $\rho$ une repr\'esentation cuspidale autoduale et des entiers $a_{0},a,b_{0},b$. On pose $s_{0}=(b_{0}-1)/2$. Chaque \'el\'ement $(\rho,A,B,\zeta)$ de $Jord(\psi)$ contribue par un facteur dans $r(s,\psi)$; le tableau ci-dessous indique s'il contribue aux p\^oles de $r(s,\psi)$ ou non; s'il remplit la condition \'ecrite dans la case, il contribue et sinon il n'y contribue pas.

\begin{center}
\begin{tabular}{|c|c|c|}
\hline
$\zeta\backslash\zeta_{0} $&+&-
\\
\hline
+ &$ B\leq B_{0}\leq A_{0}\leq A$ &
\\
\hline
-& $ B\leq A_{0}\leq A$ & $B_{0}\leq B\leq A_{0}\leq A$\\
\hline
\end{tabular}
\end{center}

En particulier si $\zeta_{0}=-$ aucun \'el\'ement de la forme $(\rho,A,B,\zeta)$ avec $\zeta=+$ contribue aux p\^oles de $r(s,\psi)$. On peut maintenant expliquer pourquoi la d\'emonstration est diff\'erente suivant les valeurs de $\zeta_{0}$; si on baisse $n(\psi,\zeta_{0})$ en augmentant \'eventuellement $\ell(\psi,-\zeta_{0})$, $r(s,\psi')/r(s,\psi)$ a de bonne chance d'\^etre holomorphe en $s=s_{0}$ si $\zeta_{0}=-$ mais pas si $\zeta_{0}=+$.

A ces contributions, s'ajoute par d\'efinition $L(St(\rho,a_{0})\times St(\rho,a_{0}), r_{G},2s)/L(St(\rho,a_{0}\times St(\rho,a_{0}),2s+1)$. Ce facteur n'a pas de p\^ole au voisinage d'un r\'eel strictement positif. Ce terme ne joue donc de r\^ole qu'en $s=0$ o\`u il peut avoir un p\^ole.

\section{D\'emonstration}
\subsection{Un lemme connu\label{lemmeconnu}}
Soient $A'_{0},B'_{0}, A,B$ des demi-entiers tels que $A'_{0}+B'_{0}$ et $A+B$ soient des entiers. On suppose que $B'_{0}\leq A'_{0}$ et que $B\leq A$. Ainsi les repr\'esentations $<A'_{0}, \cdots, B'_{0}>_{\rho}$ et $<A, \cdots, B>_{\rho}$ sont des s\'eries discr\`etes tordues par un caract\`ere.
\begin{lem}Dans le GL convenable, l'op\'erateur d'entrelacement standard $$
<A'_{0}, \cdots, B'_{0}>_{\rho}\vert\,\vert^{s}\times <A, \cdots, B>_{\rho}\rightarrow <A, \cdots, B>_{\rho}\times <A'_{0}, \cdots, B'_{0}>_{\rho}\vert\,\vert^{s}$$
est holomorphe en $s=0$ si $B'_{0}> B$ ou si $A'_{0}> A$ et a un p\^ole d'ordre 1 exactement si $B\geq B'_{0}$ et $A\geq A'_{0}$ avec l'une des in\'egalit\'e \'etant une \'egalit\'e. L'op\'erateur d'entrelacement normalis\'e \`a la Shahidi est holomorphe en $s=0$ si $B'_{0}\geq B$ ou $A'_{0}\geq A$.
\end{lem}
La deuxi\`eme partie du lemme est d\'emontr\'ee en \cite{mw} (c'est un cas particulier de \cite{mw} 1.6.3)  et est une cons\'equence des travaux d'Harish-Chandra. Le facteur de normalisation vaut
$$
\frac{L(St(\rho,A'_{0}-B'_{0}+1)\times St(\rho,A-B+1),s+(A'_{0}+B'_{0})/2-(A+B)/2)}{L(St(\rho,A'_{0}-B'_{0}+1)\times St(\rho,A-B+1),s+(A'_{0}+B'_{0})/2-(A+B)/2)}=$$
$$\frac{L(\rho\times \rho,s+\vert (A'_{0}-A+B-B'_{0})/2\vert+(A'_{0}-A+B'_{0}-B)/2)}{L(\rho\times \rho,s+(A'_{0}-B'_{0}+A-B)/2+1+(A'_{0}-A+B'_{0}-B)/2}=
$$
$$
\frac{L(\rho\times \rho,s+sup((A'_{0}-A); (B'_{0}-B)))}{L(\rho\times \rho,s+A'_{0}-B+1)}.
$$
Le facteur de normalisation a un p\^ole si $A'_{0}\leq A$ et $B'_{0}\leq B$ avec l'une des in\'egalit\'es \'etant une \'egalit\'e. Cela termine la preuve.

\subsection{R\'eduction dans le cas induit\label{casinduit}}
Ici on suppose que $Jord(\psi)$ a de la multiplicit\'e; on fixe $(\rho,a,b)\in Jord(\psi)$ et on suppose que $(\rho,a,b)$ a multiplicit\'e au moins 2 dans $Jord(\psi)$. On sait alors (cf. \ref{induction}) qu'il existe une repr\'esentation $\pi'$ dans le paquet associ\'e \`a $\psi'$ et une inclusion
$$
\pi \hookrightarrow Speh(St(\rho,a),b)\times \pi'.
$$
\begin{lem} L'holomorphie pour $N_{\psi}(s,\pi)$ r\'esulte de celle de $N_{\psi'}(s,\pi')$ en tout $s\in {\mathbb R}_{\geq 0}.$
\end{lem}
On consid\`ere la suite d'op\'erateurs d'entrelacement d'abord standard et que l'on normalisera ensuite: 
$St(\rho,a_{0})^s\times \pi \hookrightarrow$
$$
 St(\rho,a_{0})^s\times Speh (St(\rho,a),b)\times \pi' \rightarrow Speh(St(\rho,a),b)\times St(\rho,a_{0})^s\times \pi' \eqno(1)
$$
$$Speh(St(\rho,a),b)\times St(\rho,a_{0})^s\times \pi'
\rightarrow Speh(St(\rho,a),b)\times St(\rho,a_{0})^{-s}\times \pi' \eqno(2)$$
$$
Speh(St(\rho,a),b)\times St(\rho,a_{0})^{-s}\times \pi'
\rightarrow St(\rho,a_{0})^{-s}\times Speh(St(\rho,a),b)\times \pi'.\eqno(3)
$$
On a d\'efini la normalisation pour l'op\'erateur (2) et pour les op\'erateurs (1) et (3), on prend la normalisation de Shahidi ( \cite{mw}), c'est-\`a-dire pour (1) $$L(St(\rho,a_{0})\times St(\rho,a),s-(b-1)/2)/L(St(\rho,a_{0}) \times St(\rho,a),s+(b+1)/2)$$ et pour (3) $L(St(\rho,a)\times St(\rho,a_{0}),-(b-1)/2+s)/L(St(\rho,a)\times St(\rho,a_{0}),(b+1)/2+s).$

 On remarque alors que $r(\psi,s)$ est par d\'efinition le produit de ces 2 facteurs avec $r(\psi',s)$. Or les op\'erateurs d'entrelacement normalis\'es (1) et (3) sont holomorphes en tout $s\in {\mathbb R}_{\geq 0}$. Donc l'holomorphie pour $N(s,\pi)$ r\'esulte de l'holomorphie de $N(s,\pi')$. On se ram\`ene ainsi au cas o\`u $Jord(\psi)$ n'a pas de multiplicit\'e.
\subsection{R\'eduction, le cas isol\'e\label{isole}}
Ici on suppose que $Jord(\psi)$ contient un \'el\'ement $(\rho,A,B,\zeta)$ tel que $A>B$ et tel qu'il n'existe pas d'\'el\'ement $(\rho,A',B',\zeta')\in Jord(\psi)(\rho,A,B,\zeta)$ avec $\zeta'=\zeta$ et $B'\in ]B,A[$. On note $\psi'$ le morphisme qui se d\'eduit de $\psi$ en rempla\c{c}ant $(\rho,A,B,\zeta)$ par $\cup_{C\in [B,A]}(\rho,C,C,\zeta)$ et $\psi''$ le morphisme qui se d\'eduit de $\psi$ en rempla\c{c}ant $(\rho,A,B,\zeta)$ par $(\rho,A-1,B+1,\zeta)$ (ou en supprimant $(\rho,A,B,\zeta)$ si $A=B+1$).
\begin{lem}Sous les hypoth\`eses faites, 2 cas sont possibles. Soit $\pi$ appartient aussi au paquet associ\'e \`a $\psi'$, soit il existe $\pi''$ dans le paquet associ\'e \`a $\psi''$ tel que l'on ait une inclusion
$$
\pi\hookrightarrow <\rho\vert\,\vert^{\zeta B}, \cdots, \rho\vert\,\vert^{-\zeta A}>\times \pi'.
$$
\end{lem}
On fixe un ordre sur $Jord(\psi)$ pour pouvoir d\'efinir $\pi$; on impose \`a l'ordre de v\'erifier $$(\rho,A',B',\zeta')<_{Jord(\psi)}(\rho,A,B,\zeta)\hbox{ si soit $\zeta'=-\zeta$ soit $B'\leq B$.}$$ Avec l'hypoth\`ese faite, on a donc $(\rho,A',B',\zeta')>_{Jord(\psi)}(\rho,A,B,\zeta)$ exactement quand $\zeta'=\zeta$ et  $B'\geq A$. Avec ce choix d'ordre, on fixe $\psi_{>>}$ dominant tous les \'el\'ements de $Jord(\psi)$ sup\'erieur ou \'egaux \`a $(\rho,A,B,\zeta)$ et on note $\pi_{>>}$ la repr\'esentation dans le paquet associ\'e \`a $\psi_{>>}$ permettant de d\'efinir $\pi$ par module de Jacquet. On associe $\underline{t}$ et $\underline{\eta}$ les param\`etres de $\pi$ et le premier cas se produit exactement quand $\underline{t}(\rho,A,B,\zeta)=0$. En effet, supposons que $\underline{t}(\rho,A,B,\zeta)=0$; alors par d\'efinition $\pi_{>>}$ appartient au paquet associ\'e au morphisme $\psi'_{>>}$ qui se d\'eduit de $\psi'$ en rempla\c{c}ant $(\rho,A+T,B+T,\zeta)$ (o\`u $T$ est un entier grand) par $\cup_{C\in [B+T,A+T]}(\rho,C,C,\zeta)$. Les modules de Jacquet qui permettent de passer de $\psi_{>>}$ \`a $\psi$ sont les m\^emes que ceux qui permettent de passer de $\psi'_{>>}$ \`a $\psi'$ car $$\circ_{j\in [1,T]}Jac_{\zeta (B+T-j+1), \cdots \zeta (A+T-j+1)}=\circ_{C\in [A,B]}\circ_{j\in [1,T]} Jac_{\zeta (C+T-j+1)}.$$
On suppose maintenant que $\underline{t}(\rho,A,B,\zeta)\geq 1$. On note $\psi''_{>>}$ le morphisme qui se d\'eduit de $\psi_{>>}$ en rempla\c{c}ant $(\rho,A+T,B+T,\zeta)$ par $(\rho,A+T-1,B+T+1,\zeta)$ et  on sait qu'il existe $\pi''_{>>}$ une repr\'esentation dans le paquet associ\'e \`a $\psi''_{>>}$ avec une inclusion:
$$
\pi_{>>}\hookrightarrow <\rho\vert\,\vert^{\zeta (B+T)}, \cdots, \rho\vert\,\vert^{-\zeta (A+T)}>\times \pi''_{>>}.$$
On sait que $Jac_{\zeta x}\pi''_{>>}=0$ pour tout $x\in [B+1,B+T]$ car $(\rho,A',B',\zeta)\in Jord(\psi''_{>>})$ entra\^{\i}ne soit $B'\leq B$ soit $B'= B+T+1$ soit $B'>>B+T$; les formules standard donne donc l'inclusion
$$
\circ_{j\in [1,T]}Jac_{\zeta (B+T-j+1}\pi_{>>}\hookrightarrow <\rho\vert\,\vert^{\zeta B}, \cdots, \rho\vert\,\vert^{-\zeta (A+T}>\times \pi''_{>>}
$$
Ensuite on applique $\circ_{j\in [1,T]}Jac_{\zeta (B+T+1-j+1), \cdots, \zeta (A+T-1-j+1)}$ aux deux membres de l'inclusion ci-dessus et sur le membre de droite, ce module de Jacquet ne peut que s'appliquer \`a $\pi''_{>>}$; le r\'esultat est une repr\'esentation $\pi''_{>}$ dans le paquet associ\'e au morphisme $\psi''_{>}$ qui se d\'eduit de $\psi''_{>>}$ en rempla\c{c}ant $(\rho,A+T-1,B+T+1,\zeta)$ par $(\rho,A-1,B+1,\zeta)$ et cette repr\'esentation est irr\'eductible. De plus elle v\'erifie $Jac_{\zeta x}\pi''_{>}=0$ pour tout $x\in [A+1,A+T]$ car tout \'el\'ement $(\rho,A',B',\zeta)$ de $Jord(\psi''_{>})$ v\'erifie soit $B'\leq B+1$ soit $B'>>A+T$. Il reste \`a appliquer $Jac_{\zeta (A+T), \cdots, \zeta (A+1)}$ qui sur le membre de droite n'affecte pas $\pi''_{>}$. On obtient donc
$$
\circ_{j\in [1,T]}Jac_{\zeta (B+T-j+1), \cdots, \zeta (A+T-j+1)}\pi_{>>}\hookrightarrow 
<\rho\vert\,\vert^{\zeta B}, \cdots, \rho\vert\,\vert^{-\zeta A}>\times \pi''_{>}.
$$
Puis ensuite on redescend \`a $\pi$ en appliquand des $Jac_{\zeta x}$ pour $x$ parcourant un ensemble d'\'el\'ements tous strictement plus grand que $A$  (on a utilis\'e ici l'hypoth\`ese); quand on applique ces modules de Jacquet au membre de droite de l'inclusion ci-dessus, ils ne s'appliquent qu'\`a $\pi''_{>}$ pour donner une repr\'esentation $\pi''$ ayant les propri\'et\'es de l'\'enonc\'e.

\subsection{R\'eduction dans le cas de 2 blocs \'el\'ementaires cons\'ecutifs\label{consecutif}}
Ici on suppose que $Jord(\psi)$ contient 2 \'el\'ements de la forme $(\rho,A=B,\zeta), (\rho,A'=B',\zeta)$ avec $B'>B$, tels que pour tout $(\rho,A'',B'',\zeta)\in Jord(\psi)$, on ait soit $B''>B'$ soit $B''\leq B$. On note $\psi'$ le morphisme qui se d\'eduit de $\psi$ en enlevant ces 2 blocs.
\begin{lem}Avec les hypoth\`eses ci-dessus, 2 cas sont possibles; soit $Jac_{\zeta B', \cdots, \zeta (B+1)}\pi=0$ soit il existe une repr\'esentation irr\'eductible $\pi'$ dans le paquet associ\'e \`a $\psi'$ et une inclusion
$$
\pi\hookrightarrow <\rho\vert\,\vert^{\zeta B'}, \cdots, \rho\vert\,\vert^{-\zeta B}> \times \pi'.
$$
\end{lem}
Dans la preuve on interpr\`ete cette dichotomie en termes de param\`etres et ceci est important:
on met un ordre sur $Jord(\psi)$ tel que $(\rho,A'',B'',\zeta'')>_{Jord(\psi)}(\rho,A',B',\zeta)$ exactement quand $\zeta''=\zeta$ et $B''>B'$ et $(\rho,A'',B'',\zeta'')<_{Jord(\psi)}(\rho,A,B,\zeta)$ soit si $\zeta''=-\zeta$ soit si $B''\leq B$. Cela \'epuise tous les \'el\'ements de $Jord(\psi)$ d'apr\`es l'hypoth\`ese. On peut alors associer \`a $\psi$ ses param\`etres $\underline{t},\underline{\eta}$. 

\

On va montrer que le premier cas, $Jac_{\zeta B', \cdots, \zeta (B+1)}\pi=0$, se produit exactement quand $\underline{\eta}(\rho,A,B,\zeta)=-\underline{\eta}(\rho,A',B',\zeta)$. 

\

On fixe un morphisme $\psi_{>>}$ qui domine tous les \'el\'ements de $Jord(\psi)$ sup\'erieurs ou \'egaux \`a $(\rho,A,B,\zeta)$. On pose $\pi_{>>}:=\pi(\psi_{>>},\underline{t},\underline{\eta})$.
 Par d\'efinition, on sait que $\pi$ s'obtient \`a partir de  $\tilde{\pi}$ en prenant des modules de Jacquet de la forme suivante:
 $$
 \pi=Jac_{ \zeta x; x\in {\cal E}}Jac_{\zeta (B'+T'), \cdots, \zeta (B'+1)}Jac_{\zeta (B+T),\cdots, \zeta (B+1)}\pi_{>>}, \eqno(1)
 $$
 o\`u ${\cal E}$ est un ensemble totalement ordonn\'e et o\`u pour tout $x\in {\cal E}, x>B'+1$ et o\`u $T$ et $T'$ sont des entiers ''grands'' avec $T'>>T$. Quand on calcule $Jac_{\zeta B', \cdots, \zeta (B+1)}\pi$, on peut faire commuter cette op\'eration avec $Jac_{x\in {\cal E}}$ (sur le terme de droite) \`a cause de l'hypoth\`ese $x>B'+1$. Ainsi $$
 Jac_{\zeta B', \cdots, \zeta (B+1)}\pi\neq 0 \Rightarrow Jac_{\zeta B', \cdots, \zeta (B+1)}Jac_{\zeta (B'+T'), \cdots, \zeta (B'+1)}Jac_{\zeta (B+T), \cdots, \zeta (B+1)}\pi_{>>}\neq 0.
 $$
 Le terme de droite s'\'ecrit aussi $Jac_{\zeta (B'+T'), \cdots, \zeta (B+1)}Jac_{\zeta (B+T), \cdots, \zeta (B+1)}\pi_{>>}$. La r\'eciprocit\'e de Frobenius donne donc l'existence d'une repr\'esentation $\sigma$ et d'une inclusion
 $$
 \pi_{>>}\hookrightarrow \rho\vert\,\vert^{\zeta (B+T)}\times \cdots \times \rho\vert\,\vert^{\zeta (B+1)}\times  \rho\vert\,\vert^{\zeta (B'+T')}\times \cdots, \times \rho\vert\,\vert^{\zeta (B+1)}\times \sigma.
 $$
 On sait que $Jac_{\zeta x}\pi_{>>}=0$ pour tout $x\in ]B'+T',B+T[\cup ]B+T,B+1[$ et que $Jac_{\zeta (B+T), \zeta (B+T)}\pi_{>>}=0$. Ceci montre que l'inclusion ci-dessus se factorise par l'unique sous-module irr\'eductible de l'induite (pour un GL convenable) $ \rho\vert\,\vert^{\zeta (B+T)}\times \cdots \times \rho\vert\,\vert^{\zeta (B+1)}\times  \rho\vert\,\vert^{\zeta (B'+T')}\times \cdots, \times \rho\vert\,\vert^{\zeta (B+1)}$; d'o\`u une inclusion:
 $$
 \pi_{>>}\hookrightarrow <\rho\vert\,\vert^{\zeta (B+T)}, \cdots, \rho\vert\,\vert^{(B+1)}>\times <\rho\vert\,\vert^{\zeta (B'+T')}, \cdots, \rho\vert\,\vert^{\zeta (B+1)}>\times \sigma.
 $$Par irr\'eductibilit\'e, on peut \'echanger les 2 premi\`eres repr\'esentations et on obtient donc, a fortiori, $Jac_{\zeta (B'+T'), \cdots, \zeta (B+T+1)}\pi_{>>}\neq 0$.  Par d\'efinition, ceci est \'equivalent \`a $\underline{\eta}(\rho,A',B',\zeta)=\underline{\eta}(\rho,A,B,\zeta)$. Supposons donc que l'on ait cette \'egalit\'e et on note $\psi'_{>>}$ le morphisme qui se d\'eduit de $\psi_{>>}$ en enlevant les blocs $(\rho,A+T,B+T,\zeta), (\rho,A'+T',B'+T',\zeta)$. On sait qu'il existe une repr\'esentation $\pi'_{>>}$ dans le paquet associ\'e \`a $\psi'_{>>}$ et une inclusion 
 $$
 \pi_{>>}\hookrightarrow <\rho\vert\,\vert^{\zeta (B'+T')}, \cdots, \rho\vert\,\vert^{-\zeta (B+T)}>\times \pi'_{>>}.
 $$
 Pour calculer $\pi$, on revient \`a (1) ci-dessus et on commence par calculer $Jac_{\zeta (B+T), \cdots, \zeta (B+1)}\pi_{>>}$. Comme $Jac_{\zeta x}\pi'_{>>}=0$ pour tout $x\in [B+T,B+1]$ puisque $Jord(\psi'_{>>})$ ne contient aucun \'el\'ement de la forme $(\rho,\tilde{A},x,\zeta)$ pour ces valeurs de $x$, on obtient:
 $$
 Jac_{\zeta (B+T), \cdots, \zeta (B+1)}\pi_{>>}\hookrightarrow <\rho\vert\,\vert^{\zeta (B'+T')}, \cdots, \rho\vert\,\vert^{-\zeta B}>\times \pi'_{>>}.
 $$
 De m\^eme on obtient
 $$
 Jac_{\zeta (B'+T'), \cdots, \zeta (B'+1)} Jac_{\zeta (B+T), \cdots, \zeta (B+1)}\pi_{>>}\hookrightarrow <\rho\vert\,\vert^{\zeta B'}, \cdots, \rho\vert\,\vert^{-\zeta B}>\times \pi'_{>>}.
 $$
 Ensuite on doit encore appliquer $Jac_{\zeta x\in {\cal E}}$ avec $x>B'$; on obtient donc
 $$
 \pi\hookrightarrow <\rho\vert\,\vert^{\zeta B'}, \cdots, \rho\vert\,\vert^{-\zeta B}>\times Jac_{\zeta x; x\in {\cal E}}\pi'_{>>}.
 $$
 Or par d\'efinition, $\pi':=Jac_{\zeta x; x\in {\cal E}}\pi'_{>>}$ est dans le paquet associ\'e \`a $\psi'$ et on a donc d\'emontr\'e l'alternative du lemme.
\subsection{Cons\'equence sur l'holomorphie des op\'erateurs d'entrelacements}
 On fixe encore $a_{0}$ un entier et $s_{0}$ un r\'eel strictement positif et on d\'efinit $b_{0}$ par $s_{0}=(b_{0}-1)/2$. On note $\zeta_{0}$ le signe de $a_{0}-b_{0}$ si ce nombre est non nul et on pose $\zeta_{0}=+$ sinon.
 
\subsubsection{Le cas de \ref{consecutif}\label{holomorphieconsecutif}}
Ici on suppose que les hypoth\`eses de \ref{consecutif} (dont on adopte les notations) sont satisfaites; on suppose de plus que $\zeta=-\zeta_{0}$. Et on suppose que $Jac_{\zeta B', \cdots, \zeta (B+1)}\pi\neq 0$. D'o\`u l'existence de $\pi'$ 
\begin{lem} L'op\'erateur $N_{\psi}(s,\pi)$ est holomorphe en $s=s_{0}$ si l'op\'erateur $N_{\psi'}(s,\pi')$ est holomorphe en $s=s_{0}$.
\end{lem}
 
 On d\'ecompose \'evidemment l'op\'erateur d'entrelacement en produit:
 $
 St(\rho,a_{0})\vert\,\vert^{s}\times \pi\hookrightarrow $
 $$St(\rho,a_{0})\vert\,\vert^{s}\times <\rho\vert\,\vert^{\zeta B'}, \cdots, \rho\vert\,\vert^{-\zeta B}>\times \pi'
 \rightarrow 
 <\rho\vert\,\vert^{\zeta B'}, \cdots, \rho\vert\,\vert^{-\zeta B}>\times St(\rho,a_{0})\vert\,\vert^{s}\times \pi' \eqno(1)
 $$
 $$
 \rightarrow <\rho\vert\,\vert^{\zeta B'}, \cdots, \rho\vert\,\vert^{-\zeta B}>\times St(\rho,a_{0})\vert\,\vert^{-s}\times \pi'\eqno(2)
 $$
 $$\rightarrow
 St(\rho,a_{0})\vert\,\vert^{-s}\times <\rho\vert\,\vert^{\zeta B'}, \cdots, \rho\vert\,\vert^{-\zeta B}>\times \pi' .
 \eqno(3)
 $$
Supposons d'abord que $\zeta_{0}=-$; alors $\zeta=+$ et $(\rho,A,B,\zeta)$ est de la forme $(\rho,a,1)$ tandis que $(\rho,A',B',\zeta)$ et de la forme $(\rho,a',1)$ avec $a'>a$. La repr\'esentation   $<\rho\vert\,\vert^{\zeta B'}, \cdots, \rho\vert\,\vert^{-\zeta B}>$ est une repr\'esentation de Steinberg tordue, $St(\rho, (a+a')/2)\vert\,\vert^{(a'-a)/4}$. 
En $s_{0}$ l'ordre de $r(s,\psi)$ est le m\^eme que celui de $r(s,\psi')$ d'apr\`es \ref{calculdelordre}.
L'op\'erateur d'entrelacement standard (3) est holomorphe par positivit\'e stricte, d'apr\`es Harish-Chandra; l'op\'erateur d'entrelacement normalis\'e (1) est holomorphe en $s_{0}$ d'apr\`es par exemple \cite{mw} 1.6.3 parceque $-(a_{0}-1)/2+s_{0}=(b_{0}-a_{0})/2>0\geq -\zeta B=-B$. Le facteur de normalisation est $$\frac{L(St(\rho,a_{0})\times St(\rho,(a+a')/2),s_{0}-(a'-a)/4)}{L(St(\rho,a_{0})\times St(\rho,(a+a')/2),s_{0}-(a'-a)/4+1)}=$$
$$\frac{L(\rho\times \rho, \vert (a_{0}-1)/2-((a+a')/2-1)/2\vert+s_{0}-(a'-a)/4)}{L(\rho\times \rho,(a_{0}+((a+a')/2)/2+s_{0}-(a'-a)/4)}.
$$
L'ordre de ce facteur de normalisation en $s=s_{0}$ est l'ordre du num\'erateur; or ce num\'erateur n'a pas de z\'ero en $s=s_{0}$ et \`a un p\^ole d'ordre 1 quand
$$
sup((a_{0}-1)/2-((a+a')/2-1)/2-(a'-a)/4; -(a_{0}-1)/2+((a+a')/2-1)/2-(a'-a)/4)+s_{0})=0.
$$
Ou encore, $sup((a_{0}-1)/2-(a'-1)/2;-(a_{0}-1)/2+(a-1)/2)+s_{0}= sup((a_{0}-a')/2; (a-a_{0})/2)+s_{0}=0$.  Or $(a-a_{0})/2+(b_{0}-1)/2=(a-1)/2+ (b_{0}-a_{0})/2$. Or on a suppos\'e que $b_{0}>a_{0}$ et la nullit\'e ne peut donc \^etre obtenue. Ainsi l'holomorphie de $N(s,\pi)$ en $s=s_{0}$ r\'esulte de la m\^eme propri\'et\'e pour $N(s,\pi')$.

On suppose maintenant que $\zeta_{0}=+$ d'o\`u $\zeta=-$. Ici, $(\rho,A,B,\zeta)$ est de la forme $(\rho,1,b)$ et $(\rho,A',B',\zeta)$ est de la forme $(\rho,1,b')$. Et la repr\'esentation $<\rho\vert\,\vert^{\zeta B'}, \cdots, \rho\vert\,\vert^{-\zeta B}>$ est le module de Speh $<\rho\vert\,\vert^{-(b'-1)/2}, \cdots, \rho\vert\,\vert^{(b-1)/2}>$. Les op\'erateurs d'entrelacements normalis\'es (1) et (3) sont holomorphes d'apr\`es \cite{mw} 1.6.3. Il faut donc comparer les normalisations.

 Le facteur de normalisation pour l'entrelacement (1) est $L(St(\rho,a_{0})\times \rho, s_{0}-(b-1)/2)/L(St(\rho,a_{0})\times \rho, s_{0}+(b'+1)/2)$ et celui pour (2) vaut $L(St(\rho,a_{0})\times \rho,-(b'-1)/2+s_{0})/L(St(\rho,a_{0})\times \rho,(b-1)/2+s_{0})$. Or $r(s,\psi)$ est exactement le produit de ces 2 facteurs avec $r(s,\psi')$. Ainsi l'holomorphie pour $N(s,\pi)$ en $s=s_{0}$ r\'esulte de l'holomorphie pour $N(s,\pi')$ en $s=s_{0}$.
\subsubsection{Le cas de \ref{isole}\label{holomorphiecasisole}} 
Ici on suppose que les hypoth\`eses de \ref{isole} (dont on adopte les notations) sont satisfaites. On a donc d\'efini $\psi'$ et $\psi''$ et on a montr\'e que soit $\pi$ est dans le paquet associ\'e \`a $\psi'$ soit il existe une repr\'esentation $\pi''$ dans le paquet associ\'e \`a $\psi''$ avec une inclusion
$$
\pi\hookrightarrow <\rho\vert\,\vert^{\zeta B}, \cdots, \rho\vert\,\vert^{-\zeta A}>\times \pi''.\eqno(1)
$$On suppose encore que $\zeta=-\zeta_{0}$
\begin{lem}(i) Si $\pi$ est dans le paquet associ\'e \`a $\psi'$, l'holomorphie de $N_{\psi}(s,\pi)$ en $s=s_{0}$ r\'esulte de l'holomorphie de $N_{\psi'}(s,\pi)$

(ii) On suppose que l'inclusion (1) est satisfaite pour un bon choix de $\pi''$; l'holomorphie de $N_{\psi}(s,\pi)$ en $s=s_{0}$ r\'esulte de l'holomorphie de $N_{\psi''}(s,\pi'')$ en $s=s_{0}$. 
\end{lem}
Pour prouver (i), il suffit de v\'erifier que $r(s,\psi')/r(s,\psi)$ est holomorphe en $s=s_{0}$. Or on passe de $\psi$ \`a $\psi'$ en rempla\c{c}ant $(\rho,A,B,\zeta)$ par $\cup_{C\in [B,A]}(\rho,C,C,\zeta)$. On revient \`a \ref{calculdelordre}; on suppose que $\zeta_{0}=+$, par hypoth\`ese $\zeta=-$ et la contribution des \'el\'ements $(\rho,C,C,\zeta)$ aux p\^oles de $r(s,\psi')$ en $s=s_{0}$ ne peut que donner un p\^ole d'ordre 1 et ceci se produit exactement si $A_{0}\in [B,A]$; dans ce cas, $(\rho,A,B,\zeta)$ fournit aussi un p\^ole d'ordre 1 \`a $r(s,\psi)$ en $s=s_{0}$. On suppose maintenant que $\zeta_{0}=-$ et on a donc $\zeta=+$ par hypoth\`ese. Ici les \'el\'ements $(\rho,C,C,\zeta)$ ne contribue pas aux p\^oles de $r(s,\psi')$ en $s=s_{0}$. D'o\`u (i).

La preuve de (ii) est analogue \`a celle du paragraphe pr\'ec\'edent. On commence par regarder le cas o\`u $\zeta_{0}=-$ d'o\`u $\zeta_{0}=+$; ici l'op\'erateur d'entrelacement standard:
$$
St(\rho,a_{0})\vert\,\vert^{s}\times <\rho\vert\,\vert^{\zeta B}, \cdots, \rho\vert\,\vert^{-\zeta A}>\rightarrow
<\rho\vert\,\vert^{\zeta B}, \cdots, \rho\vert\,\vert^{-\zeta A}>\times St(\rho,a_{0})\vert\,\vert^{s},\eqno(2)
$$
est holomorphe en $s=s_{0}$ par positivit\'e stricte. On obtient l'holomorphie de l'op\'erateur d'entrelacement standard
$$<\rho\vert\,\vert^{\zeta B}, \cdots, \rho\vert\,\vert^{-\zeta A}>\times St(\rho,a_{0})\vert\,\vert^{-s}\rightarrow St(\rho,a_{0})\vert\,\vert^{s}\times <\rho\vert\,\vert^{\zeta B}, \cdots, \rho\vert\,\vert^{-\zeta A}>\eqno(3)
$$par \cite{mw} 1.6.3 et le r\'esultat est alors clair car $r(s,\psi'')/r(s,\psi)$ est d'ordre exactement 0 en $s=s_{0}$ d'apr\`es \ref{calculdelordre}. On suppose maintenant que $\zeta_{0}=+$ et $\zeta=-$, les op\'erateurs d'entrelacement normalis\'es \`a la Shahidi (2) et (3) sont holomorphes en $s=s_{0}$ d'apr\`es \cite{mw} et les produits de leur facteur de normalisation avec $r(s,\psi'')$ est exactement $r(s,\psi)$, c'est le calcul fait pr\'ec\'edemment.

\subsection{R\'eduction dans le cas d'un bloc non isol\'e\label{casnonisole}}
\subsubsection{Une propri\'et\'e des modules de Jacquet\label{2eproprietedujac}}
 On suppose ici que $Jord(\psi)$ contient un \'el\'ement $(\rho,A,B,\zeta)$ tel que $B>1/2$ et pour tout $(\rho,A',B',\zeta')$ $\in Jord(\psi)$ (diff\'erent de $(\rho,A,B,\zeta)$) v\'erifiant $\zeta'=\zeta$, $B'\geq B$ et $A'\geq A$, on a $B'>A$; en d'autres termes $Jord(\psi)$ ne contient pas d'\'el\'ements $(\rho,A',B',\zeta')$ avec $\zeta'=\zeta$, $B'\in [B,A]$ et $A'\geq A$.  On note $\psi'$ le morphisme qui se d\'eduit de $\psi$ en rempla\c{c}ant $(\rho,A,B,\zeta)$ par $(\rho,A-1,B-1,\zeta)$. 
 \begin{lem} On suppose aussi que $Jac_{\zeta B, \cdots, \zeta A}\pi\neq 0$. Alors il existe $\pi'$ dans le paquet associ\'e \`a $\psi'$ et une inclusion
 $$
 \pi \hookrightarrow <\rho\vert\,\vert^{\zeta B}, \cdots, \rho\vert\,\vert^{\zeta A}>\times \pi'.
 $$
 \end{lem}
On fixe un ordre sur $Jord(\psi)$ tel que $(\rho,A',B',\zeta')>_{Jord(\psi)} (\rho,A,B,\zeta)$ exactement quand $B'>A$; par hypoth\`ese si $(\rho,A',B',\zeta')<_{Jord(\psi)}(\rho,A,B,\zeta)$ on a soit $\zeta'=-\zeta$ soit $B'<B$ soit $A'<A$. Ainsi l'ordre sur $Jord(\psi)$ induit un ordre sur $Jord(\psi')$ qui a encore les bonnes propri\'et\'es. 
On choisit  $\psi_{>>}$ dominant tous les \'el\'ements de $Jord(\psi)$ sup\'erieurs strictement \`a $(\rho,A,B,\zeta)$. On note $\psi'_{>>}$ le morphisme qui se d\'eduit de $\psi_{>>}$ en rempla\c{c}ant simplement $(\rho,A,B,\zeta)$ par $(\rho,A-1,B-1,\zeta)$ et d'apr\`es ce que l'on vient de voir $\psi'_{>>}$ domine tous les blocs de Jordan de $\psi'$ strictement sup\'erieurs \`a $(\rho,A-1,B-1,\zeta)$. On note $\pi_{>>}$ la repr\'esentation du paquet associ\'e \`a $\psi_{>>}$ servant \`a d\'efinir $\pi$. On v\'erifie comme dans la preuve de \ref{consecutif} que $Jac_{\zeta B, \cdots, \zeta A}\pi\neq 0$, entra\^{\i}ne  que $Jac_{\zeta B, \cdots, \zeta A}\pi_{>>}\neq 0$. On note $\pi'_{>>}$ cette repr\'esentation dont on sait qu'elle est irr\'eductible et dans le paquet associ\'e \`a $\psi'_{>>}$. On en d\'eduit l'inclusion
$$
\pi_{>>}\hookrightarrow <\rho\vert\,\vert^{\zeta B}, \cdots, \rho\vert\,\vert^{\zeta A}> \times \pi'_{>>}.
$$
Ensuite on applique les modules de Jacquet qui calcule $\pi$ en fonction de $\pi_{>>}$; on les note $Jac_{\zeta x; x\in {\cal E}}\pi_{>>}$ $=\pi$, o\`u ${\cal E}$ est un ensemble convenable totalement ordonn\'e. Pour tout $x\in {\cal E}$, on a $x>A$ et on a donc
$$
\pi=Jac_{\zeta x; x\in {\cal E}}\pi_{>>}\hookrightarrow <\rho\vert\,\vert^{\zeta B}, \cdots, \rho\vert\,\vert^{\zeta A}>\times Jac_{\zeta x; x\in {\cal E}}\pi'_{>>}.
$$
On pose $\pi':= Jac_{\zeta x; x\in {\cal E}}\pi'_{>>}$; on sait que $\pi'$ est non nulle et c'est donc une repr\'esentation irr\'eductible dans le paquet associ\'e \`a $\psi'$.
\subsubsection{R\'eduction de l'amplitude dans le cas non isol\'e\label{reductiondelamplitude}}
Ici, on suppose que $Jord(\psi)$ contient un \'el\'ement $(\rho,A,B,\zeta)$ avec $A>B$ et que pour tout $(\rho,A',B',\zeta')$ $\in Jord(\psi)$ tel que $\zeta'=\zeta$ et $B'\in ]B,A[$, on a $A'=B'$. On note $\ell$ le cardinal des \'el\'ements $(\rho,A',B',\zeta)\in Jord(\psi)$ avec $B'\in ]B,A[$ et on les note $(\rho,A_{i}=B_{i},\zeta)$ pour $i\in [1,\ell]$ avec $B_{1}< \cdots <B_{\ell}$.  On note $\psi'_{\ell}$ le morphisme qui se d\'eduit de $\psi$ en rempla\c{c}ant les blocs  $(\rho,A_{\ell}=B_{\ell},\zeta)$ et $(\rho,A,B,\zeta)$ par $(\rho,A,A,\zeta)$ et $(\rho,A-1,B+1,\zeta)$; on note $\psi'$ le morphisme qui se d\'eduit de $\psi$ en rempla\c{c}ant simplement $(\rho,A,B,\zeta)$ par $(\rho,A-1,B+1,\zeta)$ et on note $\psi_{0}$ le morphisme qui se d\'eduit de $\psi$ en enlevant tous le blocs $(\rho,A_{i}=B_{i},\zeta)$ pour $i\in [1,\ell]$ et en rempla\c{c}ant $(\rho,A,B,\zeta)$ par $\cup_{C\in [B,A-\ell]}(\rho,C,C,\zeta)$. On suppose que l'on est pas dans le cas de \ref{isole} et on a donc certainement $A>B+1$.
\begin{lem}Sous les hypoth\`eses ci-dessus, 3 cas sont possibles: 
(i) soit il existe $\pi'_{\ell}$ dans le paquet associ\'e \`a $\psi'_{\ell}$ et une inclusion
$$
\pi\hookrightarrow <\rho\vert\,\vert^{\zeta B}, \cdots, \rho\vert\,\vert^{-\zeta B_{\ell}}>\times \pi'_{\ell}.
$$
(ii) Soit il existe $\pi'$ dans le paquet associ\'e \`a $\psi'$ et une inclusion
$$
\pi\hookrightarrow <\rho\vert\,\vert^{\zeta B}, \cdots, \rho\vert\,\vert^{-\zeta A}>\times \pi'.
$$
(iii) Soit finalement il existe $\pi'_{0}$ dans le paquet associ\'e \`a $\psi'_{0}$ et une inclusion 
$$
\pi\hookrightarrow \times_{i\in [1,\ell]}<\rho\vert\,\vert^{\zeta B_{i}}, \cdots, \rho\vert\,\vert^{-\zeta (A-i+1)}>\times \pi'_{0}.
$$
\end{lem}
On fixe un ordre sur $Jord(\psi)$ tel que pour tout  $(\rho,A',B',\zeta')>_{Jord(\psi)} (\rho,A,B,\zeta)$ on a soit $\zeta'=\zeta$, $B'>B$ soit $\zeta'=-\zeta$. Ainsi un \'el\'ement $(\rho,A',B',\zeta')$ est plus grand que $(\rho,A,B,\zeta)$ soit s'il est l'un des $(\rho,A_{i}=B_{i},\zeta)$ pour $i\in [1,\ell]$ soit si $\zeta'=\zeta$ et $B'>A$. On fixe un morphisme $\psi_{>>}$ qui domine tous les \'el\'ements de $Jord(\psi)$ strictement sup\'erieurs \`a $(\rho,A,B,\zeta)$ et une repr\'esentation $\pi_{>>}$ qui permet de calculer $\pi$. On applique \`a $\pi_{>>}$ le lemme de \ref{isole} puisque les hypoth\`eses sont satisfaites. On distingue donc suivant les 2 cas de ce lemme. 

On note $\psi'_{>>}$ le morphisme qui se d\'eduit de $\psi_{>>}$ en rempla\c{c}ant $(\rho,A,B,\zeta)$ par $(\rho,A-1,B+1,\zeta)$. On suppose d'abord qu'il existe une repr\'esentation $\pi'_{>>}$ dans le paquet associ\'e \`a $\psi'_{>>}$ et une inclusion
$$
\pi_{>>}\hookrightarrow <\rho\vert\,\vert^{\zeta B}, \cdots, \rho\vert\,\vert^{-\zeta A}>\times \pi'_{>>}.\eqno(1)
$$
Ainsi $\pi$ est de la forme $Jac_{-\zeta y; y\in {\cal F}}Jac_{\zeta x; x\in {\cal E}}\circ_{i\in [\ell,1]}Jac_{\zeta B_{i}+T_{i}, \cdots, \zeta (B_{i}+1)} \pi_{>>}$, o\`u ${\cal E}$ est un ensemble totalement ordonn\'e convenable de demi-entiers strictement sup\'erieur \`a $A$, tandis que ${\cal F}$ est un ensemble totalement ordonn\'e de demi-entiers strictement positifs et o\`u les $0<<T_{1}<< \cdots << T_{\ell}$ sont d\'etermin\'es par $\psi_{>>}$. On applique ces modules de Jacquet au membre de droite de (1) ce qui donne (malheureusement) plusieurs termes. L'un des termes consiste \`a appliquer le module de Jacquet \`a $\pi'_{>>}$ et on obtient alors l'induite $$
<\rho\vert\,\vert^{\zeta B}, \cdots, \rho\vert\,\vert^{-\zeta A}>\times \pi',$$
o\`u $\pi'$ est une repr\'esentation dans le paquet associ\'e \`a $\psi'$. Si l'inclusion (1) donne une inclusion de $\pi$ dans cette induite, on est dans le 2e cas de l'\'enonc\'e. Sinon, on remarque d'abord le point suivant; soit $i\in [1,\ell]$ et soit $x\in [B_{i}+1, A[$ alors
$$
Jac_{\zeta x}Jac_{\zeta (B_{i}+T_{i}), \cdots, \zeta (A+1)} \circ_{j\in ]i,1]}Jac_{\zeta (B_{j}+T_{j}), \cdots, \zeta (B_{j}+1)}\pi'_{>>}=0;
$$
en effet on peut \'echanger $Jac_{\zeta x}$ et $Jac_{\zeta (B_{i}+T_{i}), \cdots, \zeta (A+1)}$ car $A+1-x>1$; de plus $$\pi_{i,>>}:=Jac_{\zeta (B_{j}+T_{j}), \cdots, \zeta (B_{j}+1)}\pi'_{>>}$$ se trouve dans le paquet associ\'e au morphisme $ \psi'_{i,>>}$ qui se d\'eduit de $\psi'_{>>}$ en rempla\c{c}ant  les $(\rho,B_{j}+T_{j},B_{j}+T_{j},\zeta)$ par $(\rho,B_{j},B_{j},\zeta)$ pour tout $j\in ]i,1]$. Ainsi $Jac_{\zeta x}\pi_{i,>>}=0$ parce qu'il n'existe pas $(\rho,\tilde{A},\tilde{B},\zeta)$ dans $Jord(\psi'_{i,>>})$ avec $\tilde{B}=x$, un tel \'el\'ement devrait d\'ej\`a \^etre dans $Jord(\psi)$ et  ceci a \'et\'e exclu puisque $x\in ]B_{i},A[$.
Donc si l'on n'est pas dans le 2e cas de l'\'enonc\'e,
il existe $i\in [1,\ell]$ tel que l'inclusion (1) donne une inclusion (on oublie $\cal {F}$ qui ne joue pas de r\^ole ici)
$$
\pi \hookrightarrow Jac_{\zeta x; x\in {\cal E}}\circ_{j\in [\ell, i+1]}Jac_{\zeta (B_{j}+T_{j}), \cdots, \zeta (B_{j}+1)}\biggl( <\rho\vert\,\vert^{\zeta B}, \cdots, \rho\vert\,\vert^{-\zeta (B_{i})}> \times Jac_{\zeta (B_{i}+T_{i}), \cdots, \zeta (A+1)} \pi'_{i,>>}\biggr),
$$
avec la notation $\pi'_{i,>>}$ d\'ej\`a introduite. On veut maintenant d\'emontrer que n\'ecessairement $i=\ell$. D'abord, on remarque que le terme de droite est l'induite:
$$
<\rho\vert\,\vert^{\zeta B}, \cdots, \rho\vert\,\vert^{-\zeta (B_{i})}> \times Jac_{\zeta x; x\in {\cal E}}\circ_{j\in [\ell, i+1]}Jac_{\zeta (B_{j}+T_{j}), \cdots, \zeta (B_{j}+1)} Jac_{\zeta (B_{i}+T_{i}), \cdots, \zeta (A+1)}\pi'_{i,>>}.
$$
On remarque maintenant que la non nullit\'e du terme de droite entra\^{\i}ne une inclusion
$$
\pi'_{i,>>}\hookrightarrow \rho\vert\,\vert^{\zeta (B_{i}+T_{i})}\times \cdots \times \rho\vert\,\vert^{\zeta (A+1)}\times \rho\vert\,\vert^{\zeta (B_{i+1}+T_{i+1}}\times \cdots \times \rho\vert\,\vert^{\zeta (B_{i+1}+1)}\times \sigma,
$$
o\`u $\sigma$ est une repr\'esentation convenable. Avec les m\'ethodes habituelles, on montre que cette inclusion se factorise par
$$
\pi'_{i,>>}\hookrightarrow  <\rho\vert\,\vert^{\zeta (B_{i}+T_{i})}, \cdots, \rho\vert\,\vert^{\zeta (A+1)}>\times <\rho\vert\,\vert^{\zeta (B_{i+1}+T_{i+1})}, \cdots, \rho\vert\,\vert^{(B_{i+1}+1)}>\times\sigma.
$$
Dans le GL convenable l'induite $<\rho\vert\,\vert^{\zeta (B_{i}+T_{i})}, \cdots, \rho\vert\,\vert^{\zeta (A+1)}>\times <\rho\vert\,\vert^{\zeta (B_{i+1}+T_{i+1})}, \cdots, \rho\vert\,\vert^{(B_{i+1}+1)}>$ est irr\'eductible car les segments sont emboit\'es; on peut donc \'echanger les 2 induites. D'o\`u en particulier par r\'eciprocit\'e de Frobenius
$$
Jac_{\zeta (B_{i+1}+T_{i+1}), \cdots, (B_{i}+T_{i}+1)}\pi'_{i,>>}\neq 0.
$$
Cela s'interpr\`ete avec \ref{consecutif} en terme de param\`etres et donne en particulier l'assertion analogue pour $\pi$ \`a savoir $Jac_{\zeta B_{i+1}, \cdots, \zeta (B_{i}+1)}\pi\neq 0$. Ceci a \'et\'e exclu par hypoth\`ese et on a donc $i=\ell$. On est donc dans le 1e cas de l'\'enonc\'e. Ainsi il nous reste \`a voir le 1e cas de \ref{isole}.

On suppose donc maintenant que $\pi_{>>}$ est aussi dans le paquet associ\'e au morphisme qui se d\'eduit de $\psi_{>>}$ en rempla\c{c}ant $(\rho,A,B,\zeta)$ par $\cup_{C\in [B,A]}(\rho,C,C,\zeta)$; on note $\underline{t}$ et $\underline{\eta}$ les param\`etres pour ce morphisme permettant de d\'efinir $\pi$. Puisque $\pi$ est non nul, il faut certainement $Jac_{\zeta (B_{1}+T_{1}), \cdots, \zeta (A+1)}\pi_{>>}\neq 0$. On applique alors \ref{consecutif} et on sait que le param\`etres $\underline{\eta}$ prend la m\^eme valeur sur $(\rho,A,A,\zeta)$ et $(\rho,B_{1}+T_{1},B_{1}+T_{1},\zeta)$. On sait aussi que $\underline{\eta}$ alterne sur les blocs $(\rho,B_{i}+T_{i},B_{i}+T_{i},\zeta)$ pour $i\in [1,\ell]$ et sur les blocs $(\rho,C,C,\zeta)$ pour $C\in [B,A]$. On note $\psi'_{0,>>}$ le morphisme qui se d\'eduit de $\psi_{>>}$ en enlevant tous les blocs $(\rho,B_{i}+T_{i},B_{i}+T_{i},\zeta)$ pour $i\in [1,\ell]$ et les blocs $(\rho,C,C,\zeta)$ pour $C\in [A, A-\ell+1]$ (on est s\^ur que $\ell < A-B$ par d\'efinition). On peut donc appliquer \ref{consecutif} $\ell$ fois et obtenir l'existence d'une repr\'esentation irr\'eductible dans le paquet associ\'e \`a $\psi'_{0,>>}$ avec une inclusion
$$
\pi_{>>}\hookrightarrow \times_{i\in [1,\ell]}<\rho\vert\,\vert^{\zeta (B_{i}+T_{i}), \cdots, \rho\vert\,\vert^{-\zeta (A-i+1)}}> \times \pi'_{0,>>}. \eqno(2)
$$
Comme dans ce qui suit (1), on a 
$$
\pi=Jac_{\zeta x; x\in {\cal E}}\circ_{i\in [\ell,1]}Jac_{\zeta (B_{i}+T_{i}), \cdots, \zeta (B_{i}+1)} \pi_{>>}.
$$
On applique ce module de Jacquet au membre de droite de (2). Calculons d'abord $Jac_{\zeta (B_{1}+T_{1}, \cdots, \zeta (B_{1}+1)}$ du membre de droite. Pour tout $x\in [B_{1}+1,B_{1}+T_{1}]$, $Jac_{\zeta x}\pi'_{0,>>}=0$. On montre par r\'ecurrence descendante sur $x$ que
$$
Jac_{\zeta (B_{1}+T_{1}), \cdots, \zeta x}\pi_{>>}\hookrightarrow <\rho\vert\,\vert^{\zeta (x-1)}, \cdots, \rho\vert\,\vert^{-\zeta A}>\times_{i\in ]1,\ell]}<\rho\vert\,\vert^{\zeta (B_{i}+T_{i}), \cdots, \rho\vert\,\vert^{-\zeta (A-i+1)}}> \times \pi'_{0,>>}.
$$
S'il n'en est pas ainsi, on aurait une inclusion pour une repr\'esentation $\tau$ convenable
$$
Jac_{\zeta (B_{1}+T_{1}), \cdots, \zeta x}\pi_{>>}\hookrightarrow <\rho\vert\,\vert^{\zeta x}, \cdots, \rho\vert\,\vert^{-\zeta (A-1)}>\times \tau.
$$
Ceci entra\^{\i}ne $Jac_{\zeta (B_{1}+T_{1}), \cdots, \zeta x, \zeta x}\pi_{>>}\neq 0$ et encore $Jac_{\zeta  x}\pi_{0,>>}\neq 0$ ce qui est impossible si $x\in ]B_{1},B_{1}+T_{1}[$ et $Jac_{\zeta x, \zeta x}\pi_{0,>>}\neq 0$ si $x=B_{1}+T_{1}$, ce qui est aussi impossible. Donc finalement on trouve une inclusion
$$
Jac_{\zeta (B_{1}+T_{1}), \cdots, \zeta (B_{1}+1)}\pi_{>>}\hookrightarrow <\rho\vert\,\vert^ {\zeta B_{1}},\cdots, \rho\vert\,\vert^{-\zeta A}> \times_{i\in ]1,\ell]}<\rho\vert\,\vert^{\zeta (B_{i}+T_{i})}, \cdots, \rho\vert\,\vert^{\zeta (A-i+1)}> \times \pi'_{0,>>}.
$$
Ensuite on continue pour $i\in ]1,\ell]$ en utilisant les m\^emes propri\'et\'es de nullit\'e pour certains modules de Jacquet; ici on utilise le fait que le membre de gauche ci-dessus est une repr\'esentation irr\'eductible dans le paquet associ\'e au morphisme qui se d\'eduit de $\psi_{>>}$ en rempla\c{c}ant $(\rho,B_{1}+T_{1},B_{1}+T_{1},\zeta)$ par $(\rho,B_{1},B_{1},\zeta)$. Et finalement on d\'emontre que $\pi$ satisfait au 3e cas du lemme.
\subsubsection{Suppression d'un bloc dans le cas non isol\'e\label{casnonisolememesigne}}
On fixe $B_{0}$ un demi-entier positif ou nul et on note $\ell'$ le cardinal des $i\in [1,\ell]$ tel que $B_{i}\leq B_{0}$ (\'eventuellement $\ell'=0$).

On fait les hypoth\`eses suivantes un peu diff\'erentes de celles de \ref{reductiondelamplitude}

pour tout $i\in ]1,\ell']\cup ]\ell'+1,\ell]$, $Jac_{\zeta B_{i}, \cdots, \zeta (B_{i-1}+1)}\pi=0$

pour tout $x\in ]B,B_{0}]$, $Jac_{\zeta x}\pi=0$.

\begin{lem} Sous les hypoth\`eses ci-dessus, il existe

un morphisme $\psi'$ qui se d\'eduit de $\psi$ en rempla\c{c}ant $(\rho,A,B,\zeta)$ et $\cup_{i\in [1,\ell]}(\rho,B_{i},B_{i},\zeta)$ par des \'el\'ements de la forme $(\rho,C,C,\zeta)$ avec $C\in [B,A]$ et une repr\'esentation $\pi'$ dans le paquet associ\'e \`a $\psi'$

et deux ensembles de multisegments (dont l'un peut \^etre vide) de la forme ci-dessous avec $ t'\geq t''\in {\mathbb Z}_{\geq 0}$ et des demi-entiers $B'_{i_{t'-t''}}<\cdots < B'_{i_{1}}\leq A$ tous strictement sup\'erieurs \`a $B+t'-1$ (si $t'=t''-1$ les derni\`eres lignes n'apparaissent pas)
$$
\begin{matrix} \zeta B &\cdots&\cdots &-\zeta A
\\
\vdots &\cdots&\cdots &\vdots\\
\zeta (B+t''-1) &\cdots&\cdots &-\zeta (A-t''+1)\\
\zeta (B+t'') &\cdots &-\zeta B'_{i_{1}}\\
\vdots &\cdots& &\\
\zeta(B+t'-1)  &\cdots&-\zeta B'_{i_{t'-t''}}
\end{matrix}
$$
et le 2e multisegment est l'union de $v$ segments avec $v\in [0, \ell-\ell']$ de la forme $\zeta D_{i}, \cdots -\zeta F_{i}$ pour $i\in [1,v]$   avec $B_{0}<D_{1}< \cdots <D_{v} $ et $F_{1} > \cdots > F_{v}>B_{0}$; on note $\sigma_{1}$ et $\sigma_{2}$ les repr\'esentations de $GL$ convenable associ\'ees par Zelevinsky \`a ces multisegments et \`a la cuspidale $\rho$

tels que l'on ait $\pi\hookrightarrow \sigma_{1}\times \sigma_{2}\times \pi'$.
\end{lem}
On consid\`ere ici un ordre qui est tel que $(\rho,B_{i},B_{i},\zeta')<_{Jord(\psi)}(\rho,A,B,\zeta)$ si $i\leq \ell'$ (cette condition est vide si $\ell'=0$) et $(\rho,B_{i},B_{i},\zeta)>_{Jord(\psi)}(\rho,A,B,\zeta)$ si $i\in ]\ell',\ell]$ (cette condition est vide si $\ell'=\ell$). On fixe un morphisme $\psi_{>>}$ qui domine tous les \'el\'ements de $Jord(\psi)$ sup\'erieurs ou \'egaux \`a $(\rho,A,B,\zeta)$. On consid\`ere les param\`etres $\underline{t}$ et $\underline{\eta}$ et la repr\'esentation $\pi_{>>}$ dans le paquet associ\'e \`a $\psi_{>>}$ qui permettent de d\'efinir $\pi$. On pose $t:=\underline{t}(\rho,A,B,\zeta)$ et on pose $\psi_{t,>>}$ le morphisme qui se d\'eduit de $\psi_{>>}$ en rempla\c{c}ant $(\rho,A+T,B+T,\zeta)$ par $\cup_{C\in [B+t,A-t}(\rho,C+T,C+T,\zeta)$. Par d\'efinition, il existe une repr\'esentation dans le paquet associ\'e \`a $\psi_{t,>>}$, not\'ee $\pi_{t,>>}$ avec une inclusion
$$
pi_{>>}\hookrightarrow \bigg< \begin{matrix}
\zeta (B+T) &\cdots &\cdots &-\zeta (A+T)\\
\vdots &\vdots &\vdots &\vdots\\
\zeta (B+T+t-1)&\cdots &\cdots &-\zeta (A+T-t+1)
\end{matrix}\bigg>_{\rho}\times \pi'_{t,>>}.\eqno(1)
$$
On applique d'abord $\circ_{j\in [1,T]}Jac_{\zeta (B+T-j+1), \cdots, \zeta (B+T+-1-j+1)}$ \`a chacun des 
2 membres de (1). Quand on applique ces modules de Jacquet au terme de droite, il y a plusieurs solution car a priori on peut avoir $Jac_{\zeta x}\pi'_{t,>>}\neq 0$ pour une de ces valeurs de $x$. Toutefois si on n'utilise pas le d\'ebut des lignes pour calculer ces modules de Jacquet, il existera $x_{0}\leq B_{\ell}$ et $y_{0}\in ]B,B+t]$ avec $y_{0}\leq x_{0}$ et $$Jac_{\zeta x_{0}, \cdots, \zeta y_{0}}
\circ_{j\in [1,T]}Jac_{\zeta (B+T-j+1), \cdots, \zeta (B+T+t-1-j+1)}\pi_{>>}\neq 0.
$$
Pour repasser de $\circ_{j\in [1,T]}Jac_{\zeta (B+T-j+1), \cdots, \zeta (B+T+t-1-j+1)}\pi_{>>}$ \`a $\pi$ il faut encore appliquer des $Jac_{\zeta x}$ convenable mais on est s\^ur qu'il restera $x'_{0}\in [x_{0},y_{0}]$ tel que $Jac_{\zeta x'_{0}, \cdots, \zeta y_{0}}\pi\neq 0$. Ceci est contradictoire avec les hypoth\`eses.

On note $t'$ le plus grand entier strictement sup\'erieur ou \'egal \`a $t$, s'il existe tel que $B+t'\leq B_{\ell'-t'+t+1}$ et $B+t'\leq A-t'$ (sinon on prend $t'=t-1$). 
Si $t'>t-1$, on applique encore $\circ_{j\in [1,T]}Jac_{\zeta (B+T+t-j+1)}$; pour ne pas avoir 0, il faut certainement que le param\`etre $\underline{\eta}$ prenne la m\^eme valeur sur $(\rho,B+t+T,B+t+T,\zeta)$ et $(\rho,B_{\ell'},B_{\ell'},\zeta)$ (c'est l'argument que l'on a d\'ej\`a vu dans la d\'emonstration de \ref{reductiondelamplitude}). On note $\psi'_{>>}$ le morphisme qui se d\'eduit de $\psi_{>>}$ en rempla\c{c}ant $(\rho,A+T,B+t,\zeta)$ par $\cup_{C\in ]B+t',A-t]}(\rho,C+T,C+T,\zeta)$ et  il faut qu'il existe une repr\'esentation $\pi'_{>>}$ dans le paquet associ\'e \`a $\psi'_{>>}$ et une inclusion
$$
\pi_{>>}\hookrightarrow \bigg< \begin{matrix}
\zeta (B+T) &\cdots &\cdots &-\zeta (A+T)\\
\vdots &\vdots &\vdots &\vdots\\
\zeta (B+T+t-1)&\cdots &\cdots &-\zeta (A+T-t+1)\\
\zeta (B+T+t) &\cdots &-\zeta B_{\ell'}\\
\vdots &\vdots&\\
\zeta (B+T+t')& \cdots &-\zeta B_{\ell'-t'+t}
\end{matrix}\bigg>_{\rho}\times \pi'_{>>}.\eqno(2)
$$
On note $\psi'_{>}$ le morphisme qui se d\'eduit de $\psi_{>>}$ en rempla\c{c}ant $(\rho,A+T,B+T,\zeta)$ par $\cup_{C\in ]B+t',A-t]}$ $(\rho,C,C,\zeta)$.
Avec l'argument que l'on a d\'ej\`a donn\'e, on montre alors qu'il existe une repr\'esentation irr\'eductible $\pi'_{>}$ dans le paquet associ\'e \`a $\psi'_{>}$ et une inclusion
$$
\pi_{>}:=\circ_{j\in [1,T]}Jac_{\zeta (B+T-j+1), \cdots, \zeta (A+T-j+1)}\pi_{>>}\hookrightarrow
$$
$$
\bigg<\begin{matrix} \zeta B &\cdots& &-\zeta A
\\
\vdots &\cdots& &\vdots\\
\zeta (B+t-1) &\cdots& &-\zeta (A-t+1)\\
\zeta(B+t)  &\cdots&-\zeta B_{\ell'}\\
\vdots &\cdots &\\
\zeta(B+t')  &\cdots&-\zeta B_{\ell'-t'+t}
\end{matrix}\bigg>_{\rho}\times \pi'_{>}.\eqno(3)
$$
La d\'emonstration se termine tr\`es rapidement si $A-t< B_{\ell'+1}$. On doit encore calculer $\circ_{i\in [\ell,\ell']}\circ_{j\in [1,T_{i}}Jac_{\zeta (B_{i}+T_{i}-j+1)}$ aux 2 membres de (3) et sur le membre de droite on ne peut l'appliquer qu'\`a $\pi'_{>}$. Finalement on obtient une inclusion
$$
\pi\hookrightarrow
\bigg<\begin{matrix} \zeta B &\cdots& &-\zeta A
\\
\vdots &\cdots& &\vdots\\
\zeta (B+t-1) &\cdots& &-\zeta (A-t+1)\\
\zeta(B+t)  &\cdots&-\zeta B_{\ell'}\\
\vdots &\cdots &\\
\zeta(B+t')  &\cdots&-\zeta B_{\ell'-t'+t}
\end{matrix}\bigg>_{\rho}\times \pi'\eqno(4)
$$
avec $\pi'$ une repr\'esentation irr\'eductible dans le paquet associ\'e \`a un morphisme qui r\'epond aux conditions du lemme.

On suppose donc que $A-t\geq B_{\ell'+1}$ et on note $t''$ le plus grand entier tel que $A-t''\geq B+t'+1$ et $A-t''> B_{\ell'+t''-t+1}$. On distingue encore suivant que le param\`etre signe d\'efinissant $\pi'_{>}$ prend la m\^eme valeur sur $(\rho,A-t,A-t,\zeta)$ et sur $(\rho,B_{\ell'+1}+T_{\ell'+1},B_{\ell'+1}+T_{\ell'+1},\zeta)$ ou non. Dans le cas o\`u ce param\`etre prend la m\^eme valeur, on note $\psi''_{>>}$ le morphisme qui se d\'eduit de $\psi'_{>}$ en enlevant les blocs $(\rho,C,C,\zeta)$ pour $C\in [A-t'',A-t]$ et les blocs $(\rho,B_{\ell'+i}+T_{\ell'+i},B_{\ell'+i}+T_{\ell'+i},\zeta)$ pour $i\in [1,t-t''+1]$ et il existe d'apr\`es \ref{consecutif} appliqu\'e $t''-t+1$ fois, une repr\'esentation $\pi''_{>>}$ dans le paquet associ\'e \`a $\psi''_{>>}$ avec une inclusion:
$
\pi'_{>}\hookrightarrow $
$$<
\zeta (B_{\ell'+t''-t+1}+T_{\ell'+t''-t+1}),\cdots, -\zeta (A-t'') >_{\rho}\times \cdots \times
<
\zeta (B_{\ell'+1}+T_{\ell'+1}),\cdots, -\zeta (A-t)
>_{\rho}\times \pi''_{>>}.
$$Mais on sait aussi que le signe alterne sur chacun des blocs $(\rho,B_{\ell'+i}+T_{\ell'+i},B_{\ell'+i}+T_{\ell'+i},\zeta)$ pour $i\in [1,t-t''+1]$ (cela fait partie des hypoth\`eses) et on en d\'eduit que cette inclusion se factorise par l'unique sous-module c'est-\`a-dire:
$$
\pi'_{>}\hookrightarrow \bigg< \begin{matrix}
\zeta (B_{\ell'+t''-t+1}+T_{\ell'+t''-t+1})&\cdots &-\zeta (A-t'') \\
\vdots &\vdots&\vdots\\
\zeta (B_{\ell'+1}+T_{\ell'+1})&\cdots &-\zeta (A-t)
\end{matrix}\bigg>_{\rho}\times \pi''_{>>};\eqno(5)
$$
En mettant ensemble (3) et (5) on obtient une inclusion de $\pi_{>}$ dans une certaine induite de la forme $\tilde{\sigma}_{1}\times \tilde{\sigma}_{2}\times \pi''_{>>}$. On applique $\circ_{i\in [t''-t+1,1]}\circ_{j\in [1,T_{\ell'+i}}Jac_{\zeta (B_{\ell'+i}+T_{\ell'+i}-j+1)}$ \`a $\pi_{>}$ c'est-\`a-dire une succession de $Jac_{\zeta x}$ avec $x>B_{0}+1$. On l'applique au membre de droite de l'inclusion obtenue; cette op\'eration laisse certainement inchang\'e $\pi''_{>>}$ et modifie les d\'ebuts de lignes de la matrice de (5) et/ou les fins de lignes des matrices de (3) et (5), en gardant, dans tous les cas, la liaison des segments \`a l'int\'erieur de chaque matrice. Mais on remarque que le nouveau d\'ebut de ligne ou la nouvelle fin de colonne est certainement sup\'erieure strictement \`a $B_{0}$. Ensuite il y a encore \`a appliquer $\circ_{i\in [\ell,\ell'+t''-t+1]}\circ_{j\in [1,T_{i}}Jac_{\zeta (B_{i}+T_{i}-j+1}$; ceci peuvent s'appliquer comme ci-dessus aux d\'ebuts de lignes, \`a la fin des lignes ou \`a $\pi''_{>>}$. Ce qui est s\^ur est que $\pi''_{>>}$ devient une repr\'esentation dans un paquet associ\'e \`a un morphisme qui s'obtient en transformant les \'el\'ements $(\rho,B_{i}+T_{i},B_{i}+T_{i},\zeta)$, pour $i\in [\ell,\ell'+t''-t+1]$ en des $(\rho,C_{i},C_{i},\zeta)$ avec $C_{i}\leq A$ et les d\'ebuts ou fin de lignes modifi\'es \'etaient et restent strictement sup\'erieures \`a $B_{0}$; de plus les liaisons entre les lignes de chacune des matrices restent inchang\'ees. On obtient alors le lemme. Dans le cas o\`u le param\`etre signe ne prend pas la m\^eme valeur sur $(\rho,A-t,A-t,\zeta)$ et $(\rho,B_{\ell'+1},B_{\ell'+1},\zeta)$ on proc\`ede comme ci-dessus simplement en n'ayant pas la premi\`ere \'etape qui m\`ene \`a (5); on n'a qu'une seule matrice celle qui vient de (3). Cela termine la d\'emonstration.


\subsection{Cons\'equence pour l'holomorphie des op\'erateurs d'entrelacement}
\subsubsection{Le cas des signes contraires \label{holomorphiecasnonisolesignecontraire}}
On se place en un point $s_{0}$ de la forme $(b_{0}-1)/2$. On a d\'efini $\zeta_{0}$ et on suppose que $\zeta=-\zeta_{0}$; on suppose que les hypoth\`eses de \ref{reductiondelamplitude} sont satisfaites.

Dans le cas 2 de \ref{reductiondelamplitude}, on a vu en \ref{holomorphiecasisole} que l'holomorphie cherch\'ee se d\'eduit de la m\^eme propri\'et\'e pour $\pi'$ et $\psi'$. 

Dans le premier cas de \ref{reductiondelamplitude}, on distingue suivant les valeurs de $\zeta_{0}$. On factorise l'op\'erateur d'entrelacement en tenant compte de l'inclusion \'ecrite. Supposons que $\zeta_{0}=-$; ainsi $\zeta=+$ et la repr\'esentation $<\rho\vert\,\vert^{\zeta B}, \cdots, \rho\vert\,\vert^{-\zeta B_{\ell}}>$ est une s\'erie discr\`ete $St(\rho, (B+B_{\ell})+1)$ tensoris\'ee par le caract\`ere $\vert\,\vert^{(B-B_{\ell})/2}$.  Ici il suffit de v\'erifier que les op\'erateurs d'entrelacement standard
$$
St(\rho,a_{0})\vert\,\vert^s\times St(\rho,(B+B_{\ell})+1)\vert\,\vert^{(B-B_{\ell})/2}\rightarrow 
St(\rho,(B+B_{\ell})+1)\vert\,\vert^{(B-B_{\ell})/2}\times St(\rho,a_{0})\vert\,\vert^s
$$
est holomorphe en $s=s_{0}$ et que l'op\'erateur d'entrelacement standard
$$
St(\rho,(B+B_{\ell})+1)\vert\,\vert^{(B-B_{\ell})/2}\times St(\rho,a_{0})\vert\,\vert^{-s} \rightarrow St(\rho,a_{0})\vert\,\vert^{-s}\times St(\rho,(B+B_{\ell})+1)\vert\,\vert^{(B-B_{\ell})/2}$$
est holormophe en $s=s_{0}$. 
La premi\`ere assertion est vraie par positivit\'e car $B-B_{\ell}<0$. Pour le deuxi\`eme op\'erateur en $s=s_{0}$, $St(\rho,a_{0})\vert\,\vert^s
$ est la repr\'esentation $<-B_{0}, \cdots, -A_{0}>_{\rho}$; clairement $B>-B_{0}$ et l'holomorphie r\'esulte de \ref{lemmeconnu}.

On suppose maintenant que $\zeta_{0}=+$ et donc $\zeta=-$. Ici la repr\'esentation $<\rho\vert\,\vert^{-B}, \cdots, \rho\vert\,\vert^{B_{\ell}}>$ est un module de Speh que l'on note $J(\rho,B_{\ell},-B)$. Les op\'erateurs d'entrelacement normalis\'es
$$
St(\rho,a_{0})\vert\,\vert^{s}\times J(\rho,B_{\ell},-B)  \rightarrow J(\rho,B_{\ell},-B)\times St(\rho,a_{0})\vert\,\vert^{s}
$$
$$
J(\rho,B_{\ell},-B)\times St(\rho,a_{0})\vert\,\vert^{-s}\rightarrow St(\rho,a_{0})\vert\,\vert^{-s}\times J(\rho,B_{\ell},-B)
$$
sont holomorphes en $s=s_{0}$ d'apr\`es \cite{mw} 1.6.3. On note $r_{1}(s)$ et $r_{2}(s)$ les facteurs de normalisation et on va v\'erifier que l'ordre de la fonction
$
r_{1}(s)r_{2}(s) r(s,\psi'_{\ell})/r(s,\psi)$ en $s=s_{0}$ est exactement 0: le bloc de Jordan $(\rho,A,B,\zeta)$ donne un p\^ole d'ordre 1 \`a $r(s,\psi)$ exactement quand $B\leq A_{0}\leq A$ et le bloc $(\rho,A-1,B+1,\zeta)$ donne un p\^ole d'ordre 1 \`a $r(s,\psi'_{\ell})$ exactement quand $A-1\geq A_{0}\geq B+1$ et le nouveau bloc $(\rho,A,A,\zeta)$ donne un p\^ole d'ordre $1$ exactement quand $A_{0}=A$. Le bloc $(\rho,B_{\ell},B_{\ell},\zeta)$ donne un p\^ole d'ordre 1 \`a $r(s,\psi)$ quand $A_{0}=B_{\ell}$. La fonction $r_{1}(s)$ vaut $L(St(\rho,a_{0})\times \rho, s-B_{\ell})/L(St(\rho,a_{0}\times \rho,s+B+1)$; donc en $s=s_{0}$ sont ordre est $0$ sauf si $L(St(\rho,a_{0})\times \rho, s-B_{\ell})$ a un p\^ole ce qui se produit exactement quand $A_{0}=B_{\ell}$. La fonction $r_{2}(s)$ vaut $L(\rho\times St(\rho,a_{0}),s-B)/L(\rho\times St(\rho,a_{0}),s+B_{\ell}+1)$; cette fonction est donc d'ordre $0$ en $s=s_{0}$ sauf exactement quand la fonction $L(\rho\times St(\rho,a_{0}),s-B)$ a un p\^ole c'est-\`a-dire quand $A_{0}=B$. L'ordre de $r(s,\psi)$ en $s=s_{0}$ est donc celui de $r_{1}(s)r_{2}(s)r(s,\psi'_{\ell})$ comme annonc\'e. Cela termine l'\'etude du cas 2.

L'\'etude du 3e cas de \ref{reductiondelamplitude} est exactement analogue.
\subsubsection{Le cas de m\^eme signe\label{holomorphiecasnonisolememesigne}}
On fixe encore $s_{0}=(b_{0}-1)/2$ avec $b_{0}\in {\mathbb N}_{>1}$ et $a_{0}\in {\mathbb N}_{\geq 1}$ et on note $\zeta_{0}$ le signe de $(a_{0}-b_{0})/2$ en prenant le signe $+$ si $a_{0}=b_{0}$. On remarque que $A_{0}:= (a_{0}+b_{0})/2-1> (a_{0}-b_{0})/2=B_{0}$ car $b_{0}>1$.

Ici on suppose que $\zeta=\zeta_{0}=+$. On fixe $(\rho,A,B,+)\in Jord(\psi)$ tel que $A>B$ avec $B$ maximal pour cette propri\'et\'e. Ainsi pour tout $(\rho,A',B',+)\in Jord(\psi)$ avec $B'>B$, on a $A'=B'$. On veut montrer que l'on peut supposer que les hypoth\`eses de \ref{casnonisolememesigne} sont satisfaites. On consid\`ere donc d'abord des demi-entiers $B'',B'$ v\'erifiant $B''>B'$ tels que $(\rho,B',B',+)$ et $(\rho,B'',B'',+)$ sont des \'el\'ements de $Jord(\psi)$ cons\'ecutifs au sens de \ref{consecutif}. On suppose que soit $B_{0}\geq B''$ soit que $B_{0}<B'$.  On reprend les notations de \ref{consecutif} d'o\`u une inclusion
$$
\pi\hookrightarrow <\zeta B'', \cdots, -\zeta B'>_{\rho}\times \pi'.
$$
\begin{lem}Avec les hypoth\`eses pr\'ec\'edentes, l'holomorphie en $s=s_{0}$ de $N_{\psi}(s,\pi)$ r\'esulte de celle de $N_{\psi'}(s,\pi')$ au m\^eme point.\end{lem}
C'est \'evidemment la d\'emonstration de \ref{holomorphieconsecutif} que l'on reprend. On \'etudie les op\'erateurs d'entrelacement standard en $s=0$ (on a d\'ecal\'e par $(b_{0}-1)/2$ pour que ce soit plus joli)
$$
<A_{0}, \cdots, -B_{0}>_{\rho}\vert\,\vert^{s}\times <B'', \cdots, -B'>_{\rho}\rightarrow
<B'', \cdots, -B'>_{\rho} \times <A_{0}, \cdots, -B_{0}>_{\rho}\vert\,\vert^{s}
$$
$$
<B'', \cdots, -B'>_{\rho} \times <B_{0}, \cdots, -A_{0}>_{\rho}\vert\,\vert^{-s}\rightarrow <B_{0}, \cdots, -A_{0}>_{\rho}\vert\,\vert^{-s}\times <B'', \cdots, -B'>_{\rho}.
$$On suppose que $B_{0}\geq B''$. On applique \ref{lemmeconnu}
Le premier op\'erateur est holomorphe car $A_{0}>B_{0}\geq B''$ et le deuxi\`eme l'est car $-B'>-B''\geq -B_{0}>-A_{0}$. On suppose maintenant que $B_{0}<B'$ et avec la m\^eme r\'ef\'erence, on a l'holomorphie du premier op\'erateur car $-B_{0}>-B'$ et celle du deuxi\`eme car $B''>B'> B_{0}$. En utilisant \ref{calculdelordre}, on v\'erifie que $r(s,\psi)$ et $r(s,\psi')$ (on utilise le fait que $A_{0}>B_{0}$ vu ci-dessus) ont le m\^eme ordre en $s=s_{0}$. Cela prouve le lemme

\

On veut se ramener aussi au cas o\`u $Jac_{\zeta x}\pi=0$ pour tout $x\in ]B,B_{0}]$. On montre d'abord que si $Jac_{\zeta x}\pi\neq 0$ avec $x>B$ alors $Jord(\psi)$ contient $(\rho,x,x,\zeta)$ et en notant $\psi'$ le morphisme qui se d\'eduit de $\psi$ en rempla\c{c}ant $(\rho,x,x,\zeta)$ par $(\rho,x-1,x-1,\zeta)$, il existe une repr\'esentation $\pi'$ dans le paquet associ\'e \`a $\psi'$ avec une inclusion
$$
\pi\hookrightarrow \rho\vert\,\vert^{\zeta x}\times \pi'.
$$
On sait d\'ej\`a que $Jac_{\zeta x}\pi\neq 0$ entra\^{\i}ne l'existe d'un \'el\'ement $(\rho,A',B',\zeta)\in Jord(\psi)$ avec $B'=x$. Avec les hypoth\`eses mises ici et le fait que $x>B$, cela entra\^{\i}ne que $A'=B'$. L'assertion r\'esulte alors du lemme plus g\'en\'eral d\'emontr\'e en \ref{2eproprietedujac}.

On v\'erifie encore que l'holomorphie de $N(s,\pi)$ en $s=s_{0}$ r\'esulte de l'holomorphie de $N(s,\pi')$ en $s=s_{0}$; ce sont les m\'ethodes ci-dessus, les facteurs de normalisations ont m\^eme ordre en $s=s_{0}$. L'op\'erateur d'entrelacement standard, en $s'=0$:
$$
<A_{0}, \cdots, -B_{0}>_{\rho}\vert\,\vert^{s'}\times \rho\vert\,\vert^{\zeta x} \rightarrow \rho\vert\,\vert^{\zeta x}\times <A_{0}, \cdots, -B_{0}>_{\rho}\vert\,\vert^{s'}$$
est holomorphe d'apr\`es \ref{lemmeconnu} car $A_{0}>B_{0}\geq x$ et l'op\'erateur d'entrelacement standard
$$
\rho\vert\,\vert^{\zeta x}\times <B_{0}, \cdots, -A_{0}>_{\rho}\vert\,\vert^{-s'}\rightarrow 
<B_{0}, \cdots, -A_{0}>_{\rho}\vert\,\vert^{-s'}\times \rho\vert\,\vert^{\zeta x}
$$
est holomorphe en $s'=0$ car $x>-A_{0}$.

On reprend maintenant les notations de \ref{casnonisolememesigne} et on \'ecrit
$$
\pi\hookrightarrow \sigma_{1}\times \sigma_{2}\times \pi'.
$$
et on va montrer que $N(s,\pi)$ est holomorphe en $s=s_{0}$ si $N(s,\pi')$ est holomorphe en $s=s_{0}$. On remarque d'abord avec \ref{calculdelordre} que $r(s,\psi')/r(s,\psi)$ est d'ordre $\geq 0$ en $s=s_{0}$ et a un z\'ero d'ordre 1 si $B\leq B_{0}\leq A_{0}\leq A$; on rappelle que l'on a suppos\'e que $\zeta=\zeta_{0}=+$.
L'op\'erateur d'entrelacement standard:
$$
<A_{0}, \cdots, -B_{0}>_{\rho}\vert\,\vert^{s'}\times \sigma_{2}\rightarrow \sigma_{2}\times <A_{0}, \cdots, -B_{0}>_{\rho}\vert\,\vert^{s'}$$
est holomorphe en $s'=0$ d'apr\`es \ref{lemmeconnu} car la fin de chaque ligne d\'efinissant $\sigma_{2}$ est strictement inf\'erieur \`a $-B_{0}$. L'op\'erateur d'entrelacement standard
$$
\sigma_{2}\times <B_{0}, \cdots, -A_{0}>_{\rho}\vert\,\vert^{-s'}\rightarrow 
<B_{0}, \cdots, -A_{0}>_{\rho}\vert\,\vert^{-s'}\times \sigma_{2}
$$
est holomorphe en $s'=0$ d'apr\`es \ref{lemmeconnu} car $B_{0}$ est strictement inf\'erieur \`a tout d\'ebut de ligne de la matrice d\'efinissant $\sigma_{2}$.

Il faut faire la m\^eme chose avec $\sigma_{1}$; pour le premi\`ere op\'erateur on a l'holomorphie parce que chaque ligne d\'efinissant $\sigma_{1}$ est une s\'erie discr\`ete tensoris\'ee par un caract\`ere strictement n\'egatif. Le 2e op\'erateur n'est pas toujours holomorphe en $s'=0$. On r\'ecrit $\sigma_{1}$ comme la repr\'esentation associ\'ee aux multisegments $[B+i,-A'_{i}]$ o\`u $i\in [0,t']$ et o\`u les $A'_{i}$ sont des demi-entiers v\'erifiant $A\geq A'_{0}> \cdots >A'_{t'}$. Fixons $i\in [1,t']$; d'apr\`es \ref{lemmeconnu} l'op\'erateur d'entrelacement standard
$$
<B+i, \cdots, -A'_{i}>_{\rho}\times <B_{0},\cdots, -A_{0}>_{\rho}\vert\,\vert^{-s'} \rightarrow
<B_{0},\cdots, -A_{0}>_{\rho}\vert\,\vert^{-s'}\times <B+i, \cdots, -A'_{i}>_{\rho}
$$
est holomorphe en $s'=0$ sauf si $A'_{i}\geq A_{0}>B_{0}\geq B+i$. On note $i_{0}$ le plus grand entier s'il existe dans $[0,t']$ tel que ces in\'egalit\'es soient satisfaites. On peut alors remplacer $\sigma_{1}$ par la repr\'esentation associ\'ee aux $i_{0}+1$ premi\`eres lignes puisque les suivantes ne donne pas de p\^ole. On peut encore couper les lignes de fa\c{c}on \`a avoir une matrice rectangulaire; cela revient \`a remplacer $A'_{i}$ par $A'_{i_{0}}+i-i_{0}$ car, en notant $\sigma'$ la nouvelle repr\'esentation, on a clairement 
$$
\sigma_{1}\hookrightarrow \sigma'\times_{x\in {\cal E}}\rho\vert\,\vert^{-x},
$$
o\`u ${\cal E}$ est un ensemble de demi-entiers tous strictement sup\'erieur \`a $A_{0}+1$; les entrelacement de la repr\'esentation $<B_{0}, \cdots, -A_{0}>_{\rho}\vert\,\vert^{-s'}$ avec ces repr\'esentations $\rho\vert\,\vert^{-x}$ n'ont donc pas de p\^oles.  On peut aussi s'arranger (quitte \`a raccourcir encore les lignes) pour que dans les in\'egalit\'es $A'_{i_{0}}\geq A_{0}>B_{0}\geq B+i_{0}$ l'une soit une \'egalit\'e (c'est uniquement pour pouvoir appliquer tel quel \cite{mw} 1.6.3). On peut maintenant appliquer \cite{mw} 1.6.3 qui montre que l'op\'erateur d'entrelacement normalis\'e
$$
\sigma'\times <B_{0}, \cdots, -A_{0}>_{\rho}\vert\,\vert^{-s'}\rightarrow <B_{0}, \cdots, -A_{0}>_{\rho}\vert\,\vert^{-s'} \times \sigma'
$$est holomorphe en $s'=0$. Un calcul facile de facteur de normalisation montre que l'op\'erateur d'entrelacement standard a au plus un p\^ole d'ordre 1. Mais les in\'egalit\'es $A'_{i_{0}}\geq A_{0}>B_{0}\geq B+i_{0}$ forcent $A\geq A_{0}>B_{0}\geq B$ car $A'_{i_{0}}\leq A$ et le p\^ole que l'on a ci-dessus est compens\'e par le z\'ero de la fonction $r(s,\psi')/r(s,\psi)$ en $s=s_{0}$. Cela termine la preuve. Donc dans le cas o\`u $\zeta_{0}=+$, on a ramen\'e le th\'eor\`eme d'holomorphie pour $N(s,\pi)$ au cas particulier o\`u $\psi$ est \'el\'ementaire c'est-\`a-dire o\`u pour tout $(\rho,A,B,\zeta)\in Jord(\psi)$, $A=B$.

\subsection{R\'eduction par le bas}
\subsubsection{R\'eduction par le bas, premi\`ere \'etape\label{reductionparlebas1etape}}
Toutes nos d\'efinitions reposent  sur une m\'ethode constructive qui commence pr\'ecis\'ement par r\'eduire le plus petit bloc de Jordan; c'est ce que l'on exploite ici. On fixe $\zeta$ et on note $(\rho,A,B,\zeta)$ un \'el\'ement de $Jord(\psi)$ tel que $B$ soit minimal; on fixe un ordre sur $Jord(\psi)$ tel que $(\rho,A,B,\zeta)$ en soit le plus petit \'el\'ement. On fixe $\psi_{>>}$ un morphisme dominant $\psi$ et $\pi_{>>}$ une repr\'esentation dans le paquet associ\'e \`a $\psi_{>>}$ permettant de d\'efinir $\pi$. On note $\underline{t}$ et $\underline{\eta}$ les param\`etres. Le lemme ci-dessous est, certes, effrayant mais il donne une description tout \`a fait pr\'ecise des repr\'esentations dans le cas consid\'er\'e.
\begin{lem}
(i) On suppose que $B$ est 1/2 entier non entier et que soit $\underline{t}(\rho,A,B,\zeta)\neq 0$ soit $\underline{\eta}(\rho,A,B,\zeta)=+$. On note $\psi'_{>>}$ le morphisme qui se d\'eduit de $\psi$ en rempla\c{c}ant $(\rho,A,B,\zeta)$ par $(\rho,A-B-1/2,1/2,-\zeta)$ (ce terme n'appara\^{\i}t pas si $A=B$). Il existe une repr\'esentation $\pi'_{>>}$ dans le paquet associ\'e \`a $\psi'_{>>}$ et une inclusion
$$
\pi_{>>}\hookrightarrow <\begin{matrix} \zeta B &\cdots &\zeta A\\ \vdots &\cdots &\vdots \\ \zeta 1/2 &\cdots &\zeta (A-B+1/2)\end{matrix}>_{\rho}\times \pi'_{>>}.
$$
(ii) On suppose que $B$ est 1/2 entier non entier et que $\underline{t}(\rho,A,B,\zeta)=0$ et que $\underline{\eta}(\rho,A,B,\zeta)=-$. On note ici $\psi''_{>>}$ le morphisme qui se d\'eduit de $\psi$ en rempla\c{c}ant $(\rho,A,B,\zeta)$ par $(\rho, A-B+1/2,1/2,-\zeta)$. Il existe une repr\'esentation $\pi''_{>>}$ dans le paquet associ\'e \`a $\psi''_{>>}$ et une inclusion
$$
\pi_{>>}\hookrightarrow <\begin{matrix} \zeta B &\cdots &\zeta A\\ \vdots &\cdots &\vdots \\ \zeta 3/2 &\cdots &\zeta (A-B+3/2)\end{matrix}>_{\rho}\times \pi''_{>>}.
$$
(iii) On suppose que $B$ est entier; on note $\psi'_{0,>>}$ le morphisme qui se d\'eduit de $\psi_{>>}$ en rempla\c{c}ant $(\rho,A,B,\zeta)$ par $(\rho,A-B,0,-\zeta)$. Il existe une repr\'esentation $\pi'_{0,>>}$ dans le paquet associ\'e \`a $\psi'_{0,>>}$ et une inclusion
$$
\pi_{>>}\hookrightarrow <\begin{matrix} \zeta B &\cdots &\zeta A\\ \vdots &\cdots &\vdots \\ \zeta 1 &\cdots &\zeta (A-B+1)\end{matrix}>_{\rho}\times \pi'_{0,>>}.
$$
\end{lem}
Cela traduit pr\'ecis\'ement la construction des repr\'esentations: on commence par se ramener au cas o\`u $B=0$ ou $1/2$. On applique simplement \cite{paquetdiscret} 3.1: on note $\psi''_{>>}$ le morphisme qui se d\'eduit de $\psi_{>>}$ en rempla\c{c}ant $(\rho,A,B,\zeta)$ par $(\rho,A-[B],B-[B],\zeta)$ et on pose $\pi''_{>>}:=\circ_{i\in [1,[B]]}Jac_{\zeta B-i+1, \cdots, \zeta A-i+1}\pi_{>>}$. On sait que $\pi''_{>>}$ est non nul et est dans le paquet de repr\'esentations associ\'ees \`a $\psi''_{>>}$ et on a une inclusion:
$$
\pi_{>>}\hookrightarrow <\begin{matrix} \zeta B &\cdots, &\zeta A\\ \vdots &\cdots &\vdots \\ \zeta B-[B]+1 &\cdots &\zeta (A-[B]+1)\end{matrix}>_{\rho}\times \pi''_{>>}.
$$
Maintenant on distingue suivant que $B=[B]$ ou non. Dans le premier cas, on sait que le paquet associ\'e \`a $\psi''$ se d\'efinit aussi en rempla\c{c}ant $(\rho,A-B,0,\zeta)$ par $(\rho,A-B,0,-\zeta)$ et on a donc directement (iii). On suppose maintenant que $B-[B]=1/2$; on a directement (ii), comme pr\'ec\'edemment, quand les param\`etres $\underline{t}$ et $\underline{\eta}$ qui d\'efinissent $\pi$ v\'erifient $\underline{t}(\rho,A,B,\zeta)=0$ et $\underline{\eta}(\rho,A,B,\zeta)=-$; c'est la d\'efinition. 
Supposons que $\underline{t}(\rho,A,B,\zeta)=0$ mais que $\underline{\eta}(\rho,A,B,\zeta)=+$. On peut alors remplacer $(\rho,A-[B],1/2,\zeta)$ par $\cup_{C\in [B,A]}(\rho,C,C,\zeta)$ avec un nouveau param\`etre $\underline{\eta}$ qui alterne sur ces blocs en commen\c{c}ant par $+$ sur $(\rho,1/2,1/2,\zeta)$. On a alors par d\'efinition l'existence d'une repr\'esentation $\pi'_{1}$ dans le morphisme qui se d\'eduit de $\psi_{>>}$ en rempla\c{c}ant $(\rho,A-[B],1/2,\zeta)$ par $\cup_{C\in [1/2,A-1]}(\rho,C,C,\zeta)$ et une inclusion
$$
\pi_{>>}\hookrightarrow <\rho\vert\,\vert^{\zeta 1/2}, \cdots, \rho\vert\,\vert^{\zeta A-[B]}>_{\rho}\times \pi'_{1}.
$$
Le param\`etre d\'efinissant $\pi'_{1}$ contient en particulier un signe qui alterne sur les blocs $(\rho,C,C,\zeta)$ pour $C\in [1/2,A-1]$, en commen\c{c}ant par $-$ sur $(\rho,1/2,1/2,\zeta)$; on peut donc remplacer $\zeta$ par $-\zeta$ sur ces blocs et finalement remplacer l'ensemble de ces blocs par $(\rho,A-[B]-1,1/2,-\zeta)$. Cela donne (i) dans ce cas. On remarque ici que les param\`etres de $\pi'_{>>}$ se d\'eduisent naturellement de ceux de $\pi_{>>}$ avec toutefois $\underline{t}_{\pi'_{>>}}(\rho,A-B-1/2,1/2,-\zeta)=0$ et $\underline{\eta}_{\pi'_{>>}}(\rho,A-B-1/2,1/2,-\zeta)=-=-\underline{\eta}(\rho,A,B,\zeta)$. 
 
On consid\`ere maintenant les cas restants. On note d'abord $\pi'_{1,>>}:= \circ_{i\in [1, [B]]}Jac_{\zeta (B-i+1), \cdots, \zeta (A-i+1)}\pi_{>>}$ et on sait que cette repr\'esentation est irr\'eductible dans le paquet associ\'e au morphisme $\psi_{1,>>}$ qui se d\'eduit de $\psi_{>>}$ en rempla\c{c}ant $(\rho,A,B,\zeta)$ par $(\rho,A-[B],1/2,\zeta)$. De plus on a une inclusion
$$
\pi_{>>}\hookrightarrow \bigg< \begin{matrix} \zeta B &\cdots & \zeta A\\
\vdots&\vdots&\vdots\\
\zeta 3/2 &\cdots & \zeta (A-[B]+1)
\end{matrix}
\bigg>_{\rho}\times \pi_{1,>>}.\eqno(1)
$$
On pose $t_{0}:=\underline{t}(\rho,A,B,\zeta)$ et on suppose que $t_{0}\neq 0$. Il faut aussi distinguer suivant la valeur de $\underline{\eta}(\rho,A,B,\zeta)$; on va montrer que les param\`etres pour la repr\'esentation obtenue sont 
$$
\underline{t}_{\pi'_{>>}}(\rho,A-B-1/2,1/2,-\zeta)=\begin{cases}t_{0} \hbox { si }\underline{\eta}(\rho,A,B,\zeta)=+\\
t_{0}-1 \hbox { si }\underline{\eta}(\rho,A,B,\zeta)=-.\end{cases}
$$
$$
\underline{\eta}_{\pi'_{>>}}(\rho,A-B-1/2,1/2,-\zeta)=\begin{cases}- \hbox { si }\underline{\eta}(\rho,A,B,\zeta)=+\\
+ \hbox { si }\underline{\eta}(\rho,A,B,\zeta)=-.\end{cases}
$$
On pose $\pi'_{2,>>}:=\circ_{i\in [1,t_{0}]}Jac_{\zeta 1/2+i-1, \cdots, -\zeta (A-[B]-i+1}\pi'_{1}$ puis $$\pi'_{3,>>}:= \circ_{j\in [1,t'_{0}]}Jac_{\zeta (1/2+t_{0}-j+1), \cdots, \zeta (A-[B]-t_{0}-j+1)} \pi'_{2,>>},$$ o\`u $t'_{0}=t_{0}$ si $\underline{\eta}(\rho,A,B,\zeta)=-$ et $t'_{0}+1$ si $\underline{\eta}(\rho,A,B,\zeta)=+$. On sait que $\pi'_{3,>>}$ est une repr\'esentation irr\'eductible dans le paquet associ\'e au morphisme qui se d\'eduit de $\psi_{>>}$ en rempla\c{c}ant $(\rho,A,B,\zeta)$ par $\cup_{C\in [1/2, A-[B]-t_{0}-t'_{0}}(\rho,C,C,\zeta)$ et avec un param\`etre signe qui alterne sur tous ces blocs en commen\c{c}ant par $-$; on peut donc encore remplacer $\zeta$ en $-\zeta$. De plus, $\pi'_{1,>>}$ est l'unique sous-module irr\'eductible de l'induite
$$
\bigg<\begin{matrix} \zeta 1/2 &\cdots&\zeta(1/2-t'_{0}-1) &\cdots &-\zeta (A-[B])\\
\vdots &\vdots &\vdots &\vdots\\
\zeta (1/2+t_{0}-1 &\cdots& \zeta(1/2+t_{0}-t'_{0}) &\cdots &-\zeta (A-[B]-t_{0}+1)\\
\zeta (1/2+t_{0}) &\cdots & \zeta (1/2+t_{0}-t'_{0}+1)&\\
\vdots &\vdots&\vdots\\
\zeta (A-[B]-t_{0}+1) &\cdots &\zeta (A-[B]-t_{0}-t'_{0})&
\end{matrix}
\bigg>_{\rho}
 \times \pi'_{3,>>}.
 $$
 De plus l'induite $\rho\vert\,\vert^{x}\times \pi'_{3}$ est irr\'eductible pour tout $x\in ]A-[B]-t_{0}-t'_{0}+1, A-[B]]$. On peut donc retourner les $t'_{0}$ derni\`eres colonnes pour les mettre \`a la fin des $t'_{0}$ premi\`eres colonnes et l'induite ci-dessus est isomorphe \`a l'induite
 $$
 \bigg<\begin{matrix} \zeta 1/2 &\cdots&\zeta(1/2-t'_{0}-1) &\cdots &-\zeta (A-[B]-t'_{0}+1)\\
\vdots &\vdots &\vdots &\vdots\\
\zeta (1/2+t_{0}-1 &\cdots& \zeta(1/2+t_{0}-t'_{0}) &\cdots &-\zeta (A-[B]-t_{0}+1-t'_{0})\\
\zeta (1/2+t_{0}) &\cdots & \zeta (1/2+t_{0}-t'_{0}+1)&\\
\vdots &\vdots&\vdots\\
\zeta (A-[B]-t_{0}+1) &\cdots &\zeta (A-[B]-t_{0}-t'_{0})\\
\vdots &\vdots &\vdots\\
\zeta (A-[B]) &\cdots &\zeta (A-[B]-t'_{0}+1)
\end{matrix}
\bigg>_{\rho}
 \times \pi'_{3,>>}.
 $$
 On r\'eserve la premi\`ere colonne de la matrice, on r\'ecrit les $t'_{0}-1$ colonnes suivantes
 $$
\mathcal{A}_{t'_{0}-1}:= \begin{matrix} -\zeta (1/2) &\cdots &-\zeta (1/2+(t'_{0}-1)+1)\\
 \vdots&\vdots&\vdots\\
 \zeta (A-[B]-1) &\cdots&\zeta (A-[B]-1-(t'_{0}-1)+1)
 \end{matrix}
 $$
 Et les colonnes restantes s'\'ecrivent sous la forme
$$
\mathcal{B}:= \begin{matrix} -\zeta (1/2+t'_{0}-1 )&\cdots &-\zeta (A-[B]-t'_{0}+1)\\
\vdots&\vdots&\vdots\\
-\zeta (1/2+t'_{0}-t_{0}) &\cdots &-\zeta (A-[B]-t_{0}-t'_{0}+1)
\end{matrix}
$$
L'induite $<{\mathcal{B}}>_{\rho}\times \pi'_{3,>>}$ a un unique sous-module irr\'eductible qui se trouve dans le paquet associ\'e au morphisme qui se d\'eduit de $\psi_{>>}$ en rempla\c{c}ant $(\rho,A,B,\zeta)$ par $\cup_{C\in [1/2+t'_{0}-1, A-[B]-t'_{0}+1]}(\rho,C,C,-\zeta)$; on note $\pi'_{4,>>}$ ce sous-module irr\'eductible, son param\`etre signe alterne sur les blocs \'ecrits en commen\c{c}ant par $+$ si $t'_{0}=t_{0}$ et par $-$ si $t'_{0}=t_{0}+1$. On remarque ensuite que $<{\mathcal{A}_{t'_{0}-1}}>_{\rho}\times \pi'_{4,>>}$ a un unique sous-module irr\'eductible qui est exactement le $\pi'_{>>}$ d\'ecrit. La premi\`ere colonne se retourne pour s'ajouter comme derni\`ere ligne \`a la matrice de (1) et on obtient le (i) de l'\'enonc\'e.
\subsubsection{R\'eduction par le bas, fin\label{reductionparlebas}}
On reprend les hypoth\`eses de \ref{reductionparlebas1etape} en les renfor\c{c}ant ici. On a fix\'e $(\rho,A,B,\zeta)$ dans $Jord(\psi)$ tel que pour tout $(\rho,A',B',\zeta')\in Jord(\psi)$ avec $\zeta'=\zeta$ on a $B'\geq B$ et on suppose en plus ici que pour tout $(\rho,A',B',\zeta')\in Jord(\psi)$ si $\zeta'=-\zeta$, on a $A'=B'$ et que pour $(\rho,A'=B',-\zeta), (\rho,A''=B'',-\zeta)$ des \'el\'ements de $Jord(\psi)$ avec $B'>B''$, $Jac_{-\zeta B', \cdots, -\zeta(B''+1)}\pi=0$. On note $\psi'$ le morphisme qui se d\'eduit de $\psi$ en enlevant $(\rho,A,B,\zeta)$ et tous les \'el\'ements de la forme $(\rho,A'=B',-\zeta)$ avec $B'< A$ et en ajoutant $(\rho,A-[B]-1,1/2,-\zeta)$ resp. $(\rho,A-[B],1/2,-\zeta)$ et resp. $(\rho,A-B,0,-\zeta)$ suivant les cas de (i) \`a (iii) de \ref{reductionparlebas1etape}. En suivant les 3 cas de cette r\'ef\'erence, on pose $\delta_{B}=1/2,3/2,1$.
\begin{lem}Il existe:

 un ensemble de demi-entiers (ensemble \'eventuellement vide) $\mathcal{C}$ tous inf\'erieurs ou \'egaux \`a $A$ qui permet de d\'efinir le morphisme $\psi''$ obtenu en rajoutant \`a $Jord(\psi')$ les \'el\'ements $(\rho,C,C,-\zeta)$ pour $C\in {\mathcal{C}}$
 
 une repr\'esentation irr\'eductible $\pi''$ dans le paquet associ\'e \`a $\psi''$ 
 
 une suite de demi-entier $ B''_{1}<\cdots <B''_{\ell}\leq A$ avec $\ell=B-\delta_{B}+1$

tels que $\pi$ soit un sous-module de l'induite
$$
<\begin{matrix}\zeta B &\cdots& & \zeta B''_{\ell}\\
\vdots &\vdots&&\\ 
 \zeta \delta_{B}&\cdots &\zeta B''_{1}&\end{matrix}>_{\rho}\times \pi'',
 $$en prenant comme convention que les derni\`eres lignes peuvent ne pas appara\^{\i}tre si $\zeta$ fois l'extr\'emit\'e de la ligne est strictement inf\'erieure \`a $\zeta$ fois le d\'ebut de la ligne \end{lem}
 On reprend le lemme de \ref{reductionparlebas1etape}; pour passer de $\pi_{>>}$ \`a $\pi$ il faut d'abord appliquer des $Jac_{\zeta x; x\in {\cal E}}$ o\`u ${\cal E}$ est un ensemble totalement ordonn\'e de demi-entiers tous sup\'erieurs strictement \`a $B$; cette op\'eration appliqu\'ee au membre de droite des inclusions de loc.cite. ne s'applique qu'aux repr\'esentations $\pi'_{>>},\pi''_{>>},\pi_{0,>>}$. Ensuite il faut faire redescendre les \'el\'ements de $Jord(\psi_{>>})$ qui domine les \'el\'ements de $Jord(\psi)$ de la forme $(\rho,A'=B',-\zeta)$. On va du plus petit au plus grand. Pour cela il peut y avoir un choix quand on applique le module de Jacquet aux membres de droite des inclusions de \ref{reductionparlebas1etape} tant que $A'=B'<A$, on peut appliquer ces modules de Jacquet \`a la fois aux analogues de $\pi'_{>>}\cdots$ soit aux bouts des lignes de la matrice d\'ecrire en loc. cite. Comme expliqu\'e au d\'ebut de \ref{casnonisolememesigne}, on obtient alors les donn\'ees un peu impr\'ecises de l'\'enonc\'e et cela termine la d\'emonstration.
\subsection{Fin de la preuve de l'holomorphie des op\'erateurs d'entrelacement dans le cas $\zeta_{0}=-$\label{finpreuve-}}
On fixe ici $a_{0},b_{0}, s_{0}=(b_{0}-1)/2$ et on suppose que $a_{0}<b_{0}$. 
\begin{prop}Sous l'hypoth\`ese faite $a_{0}<b_{0}$, l'op\'erateur d'entrelacement $N_{\psi}(s,\pi)$ est holomorphe en $s=s_{0}$. \end{prop}
Gr\^ace \`a \ref{casinduit} on suppose que $Jord(\psi)$ est sans multiplicit\'e.
En \ref{holomorphiecasnonisolesignecontraire}, on a ramen\'e la propri\'et\'e d'holomorphie de $N(s,\pi)$ en $s=s_{0}$ au cas o\`u pour tout $(\rho,A,B,\zeta)\in Jord(\psi)$ avec $\zeta=+$, $A=B$. Si pour tout $(\rho,A,B,\zeta)\in Jord(\psi)$, on a $\zeta=+$, alors $\psi$ est temp\'er\'e et il est facile de v\'erifier que $r(s,\psi)$ est holomorphe non nul en $s=s_{0}$. L'holomorphie cherch\'ee est donc \'equivalente \`a l'holomorphie de l'op\'erateur d'entrelacement standard:
$$
St(\rho,a_{0})\vert\,\vert^{s}\times \pi \rightarrow St(\rho,a_{0})\vert\,\vert^{-s}\times \pi.
$$
Mais ici $\pi$ est une s\'erie discr\`ete et le r\'esultat a \'et\'e d\'emontr\'e par Harish-Chandra. On raisonne donc par r\'ecurrence sur le nombre d'\'el\'ements de $Jord(\psi)$ de la forme $(\rho,A,B,-)$. On suppose qu'il existe de tels \'el\'ements et on en fixe 1 avec $B$ minimum. On lui applique \ref{reductionparlebas}. On reprend les notations de loc. cite et on pose
$$
\sigma:=<\begin{matrix}- B &\cdots& &- B''_{\ell}\\
\vdots &\vdots&\\ 
- \delta_{B}&\cdots &- B''_{1}\end{matrix}>_{\rho}
 $$
 vue comme une repr\'esentation d'un GL convenable. La repr\'esentation $\pi''$ est associ\'ee \`a un morphisme $\psi''$ auquel on peut appliquer l'hypoth\`ese de r\'ecurrence; de plus d'apr\`es \ref{calculdelordre} $r(s,\psi'')/r(s,\psi)$ en $s=s_{0}$ est holomorphe avec un z\'ero d'ordre 1 exactement quand
 $$
 B_{0}\leq B\leq A_{0}\leq A.\eqno(1)
 $$
 L'op\'erateur d'entrelacement standard $St(\rho,a_{0})\vert\,\vert^{s}\times \sigma \rightarrow \sigma \times St(\rho,a_{0})\vert\,\vert^{s}$ est  holomorphe en $s=s_{0}$ par positivit\'e. Le point est donc de d\'emontrer que l'op\'erateur d'entrelacement standard $\sigma \times St(\rho,a_{0})\vert\,\vert^{-s}\rightarrow St(\rho,a_{0})\vert\,\vert^{-s}\times \sigma$ a au plus un p\^ole d'ordre 1 en $s=s_{0}$ et que ce p\^ole n'existe que si (1) est satisfait. D'abord on remarque que pour $i\in [1,\ell]$ l'op\'erateur d'entrelacement standard:
 $$
 <\rho\vert\,\vert^{-\delta_{B}-i+1}, \cdots, \rho\vert\,\vert^{-B''_{i}}>\times St(\rho,a_{0})\vert\,\vert^{-s}\rightarrow St(\rho,a_{0})\vert\,\vert^{-s}\times  <\rho\vert\,\vert^{-\delta_{B}-i+1}, \cdots, \rho\vert\,\vert^{-B''_{i}}>
 $$
 est holomorphe  en $s=s_{0}$ (cf. \ref{lemmeconnu}) sauf si $B_{0}\leq \delta_{B}+i-1 \leq A_{0}\leq B''_{i}$ cas o\`u il a un p\^ole d'ordre au plus 1. Donc si pour tout $i\in [1,\ell]$ ces in\'egalit\'es ne sont pas satisfaites, on a l'holomorphie cherch\'ee. Sinon, on note $j$ le plus petit entier tel que $B_{0}\leq \delta_{B}+j-1 \leq A_{0}\leq B''_{j}$ et on note $\sigma_{j}$ la repr\'esentation associ\'ee \`a la matrice
 $$
 \begin{matrix}-B &\cdots & -B''_{j}-\ell+j\\
 \vdots &\vdots &\vdots\\
 -\delta_{B}-j+1 &\cdots & -B''_{j}
 \end{matrix}.
 $$
 Ainsi on peut \'ecrire $$\sigma\hookrightarrow \sigma_{j}\times_{x\in {\cal E}}\rho\vert\,\vert^{-x} \times_{i\in ]j,1]}<\rho\vert\,\vert^{-\delta_{B}-i+1}, \cdots, \rho\vert\,\vert^{-B''_{i}}>,$$o\`u ${\cal E}$ est un ensemble totalement ordonn\'e d'\'el\'ements tous strictement sup\'erieurs \`a $A_{0}+1$. Les propri\'et\'es cherch\'ees pour l'entrelacement avec $\sigma$ d\'ecoulent donc des m\^emes propri\'et\'es pour l'entrelacement avec $\sigma_{j}$ rempla\c{c}ant $\sigma$; pour la m\^eme raison et quitte \`a  augmenter $B_{0}$ on peut supposer que $B_{0}=\delta_{B}+j-1$ pour pouvoir ais\'ement appliquer \cite{mw} 1.6.3. Dans cette r\'ef\'erence, on a alors montr\'e que l'op\'erateur d'entrelacement normalis\'e  \`a la Shahidi
 $$
 \sigma_{j}\times St(\rho,a_{0})\vert\,\vert^{-s}\rightarrow St(\rho,a_{0})\vert\,\vert^{-s}\times \sigma_{j}
 $$
 est holomorphe en $s=s_{0}$; il faut donc prendre en compte les p\^oles \'eventuels du facteur de normalisation. Le facteur de normalisation est, \`a une fonction holomorphe inversible pr\`es,
$$L(St(\rho,\alpha)\times St(\rho,a_{0}), s-x_{f})/L(St(\rho,\alpha)\times St(\rho,a_{0}), s-x_{f}+1)$$
o\`u $\alpha=B''_{j}-(\delta_{B}+j-1)+1$, $x_{f}$ et $x_{d}$ sont tels que $(\alpha-1)/2+x_{f}=-B$ et $-(\alpha-1)/2+x_{d}=-B''_{j}$. Le num\'erateur a un p\^ole d'ordre 1 exactement quand $B_{0}\leq B \leq A_{0}\leq B''_{j}+\ell-j$; ces derni\`eres \'egalit\'es entrainent en particulier que $B_{0}\leq B\leq A_{0}\leq A$ car $A\geq B''_{\ell}\geq B''_{j}+\ell-j$. On en d\'eduit donc l'assertion cherch\'ee.
\subsection{Fin de la preuve de l'holomorphie des op\'erateurs d'entrelacement dans le cas $\zeta_{0}=+$ \label{holomorphiecaselementaire}}
Ici on fixe $a_{0},b_{0}, s_{0}=(b_{0}-1)/2$ et on suppose que $a_{0}\geq b_{0}$. En \ref{holomorphiecasnonisolememesigne} et \ref{holomorphiecasnonisolesignecontraire}, on a r\'eduit l'holomorphie de $N(s,\pi)$ en $s=s_{0}$ au cas o\`u le morphisme $\psi$ est \'el\'ementaire. On suppose donc ici que $\psi$ est \'el\'ementaire, c'est-\`a-dire que pour tout $(\rho,A,B,\zeta)\in Jord(\psi)$, on a $A=B$ et que $Jord(\psi)$ n'a pas de multiplicit\'e. On suppose que $s_{0}>0$ ce qui veut dire que $b_{0}>1$ et on a donc $A_{0}>B_{0}$. Ainsi si pour tout $(\rho,A,B,\zeta)\in Jord(\psi)$, $\zeta=+$, $\pi$ est une s\'erie discr\`ete, le facteur $r(s,\psi)$ est holomorphe non nul en $s=s_{0}$ et l'holomorphie de $N(s,\pi)$ en $s=s_{0}$ est \'equivalente \`a celle de l'op\'erateur d'entrelacement standard, holomorphie qui r\'esulte des travaux d'Harish-Chandra. Pour d\'emontrer l'holomorphie de $N(s,\pi)$ en $s=s_{0}$, on raisonne donc par r\'ecurrence sur le nombre d'\'el\'ements de $Jord(\psi)$ de la forme $(\rho,A=B,-)$. On suppose qu'il existe un tel \'el\'ement et on le fixe avec $B$ minimum. On fixe un ordre sur $Jord(\psi)$ tel que cet \'el\'ement soit minimal et on note $\psi_{>>}$ un morphisme qui domine tous les \'el\'ements de $Jord(\psi)$ qui lui sont strictement sup\'erieurs et on note $\pi_{>>}$ la repr\'esentation dans le paquet associ\'e \`a $\psi_{>>}$ permettant de d\'efinir $\pi$. On peut utiliser \ref{reductionparlebas1etape}: l'une des situations suivantes est r\'ealis\'ee

(i) $B$ est un demi-entier non entier, on note $\psi'_{>>}$ le morphisme qui se d\'eduit de $\psi$ en enlevant $(\rho,B,B,-)$ et il existe une repr\'esentation $\pi'_{>>}$ dans le paquet associ\'e \`a $\psi'_{>>}$ avec une inclusion
$$\pi_{>>}\hookrightarrow
<\rho\vert\,\vert^{-B}, \cdots, \rho\vert\,\vert^{-1/2}>\times \pi'_{>>};
$$
(ii) $B$ est un demi-entier non entier, on note $\psi'_{>>}$ le morphisme qui se d\'eduit de $\psi$ en rempla\c{c}ant $(\rho,B,B,-)$ par $(\rho,1/2,1/2,+)$ et il existe une repr\'esentation $\pi'_{>>}$ dans le paquet associ\'e \`a $\psi'_{>>}$ avec une inclusion
$$\pi_{>>}\hookrightarrow
<\rho\vert\,\vert^{-B}, \cdots, \rho\vert\,\vert^{-3/2}>\times \pi'_{>>};
$$
(iii) $B$ est entier et on note $\psi'_{>>}$ le morphisme qui se d\'eduit de $\psi_{>>}$ en rempla\c{c}ant $(\rho,B,B,-)$ par $(\rho,0,0,+)$ et il existe $\pi'_{>>}$ dans le paquet associ\'e \`a ce morphisme avec une inclusion
$$\pi_{>>}\hookrightarrow
<\rho\vert\,\vert^{-B}, \cdots, \rho\vert\,\vert^{-1}>\times \pi'_{>>}.
$$
On passe maintenant de $\pi_{>>}$ \`a $\pi$. Pour cela on fait d'abord redescendre les \'el\'ements de $Jord(\psi_{>>})$ dominant les \'el\'ements de $Jord(\psi)$ de la forme $(\rho,B',B',-)$; n\'ecessairement $B'>B$ et on applique donc des $Jac_{-x}$ pour $x>B+1$; on note $\pi_{+,>>}$ le r\'esultat. On applique ces modules de Jacquet aux membres de droite des inclusions ci-dessus et ces modules de Jacquet ne s'appliquent  qu'\`a la repr\'esentation  $\pi'_{>>}$; on note $\pi'_{+,>>}$ le r\'esultat. On a donc des inclusions analogues avec $_{>>}$ remplac\'e par $_{+,>>}$.

Il faut maintenant faire redescendre les \'el\'ements de $Jord(\psi_{>>})$ qui dominent les \'el\'ements de $Jord(\psi)$ de la forme $(\rho,B',B',+)$; on note $v$ le nombre de ces \'el\'ements et on les ordonne par $0\leq B'_{1}< \cdots < B'_{v}$. Pour tout $i\in [1,v]$ il existe un entier $T'_{i}$ tel que $(\rho,B'_{i}+T'_{i},B'_{i}+T'_{i},+)\in Jord(\psi'_{+,>>})$. On doit calculer $\circ_{i\in [v,1]}Jac_{B'_{i}+T'_{i}, \cdots, B'_{i}+1}\pi_{+,>>}$ et on obtient $\pi$. On applique ces modules de Jacquet aux membres de droite des inclusions obtenues; dans le cas (i), ces modules de Jacquet ne s'appliquent qu'\`a $\pi'_{+,>>}$ et on se trouve donc dans le cas

(i)bis on note $\psi'$ le morphisme qui se d\'eduit de $\psi$ en supprimant$(\rho,B,B,-)$ et il existe une repr\'esenta\-tion $\pi'$ dans le paquet associ\'e \`a ce morphisme avec une inclusion
$$\pi\hookrightarrow
<\rho\vert\,\vert^{-B}, \cdots, \rho\vert\,\vert^{-1/2}>\times \pi';
$$
 Dans les cas (ii); si $B'_{1}>1/2$, c'est le m\^eme ph\'enom\`ene que ci-dessus en prenant pour $\psi'$ le morphisme qui se d\'eduit de $\psi$ en rempla\c{c}ant $(\rho,B,B,-)$ par $(\rho,1/2,1/2,+)$. Si $B'_{1}=1/2$, on a 2 possibilit\'e soit
 $$
 Jac_{B'_{1}+T'_{1}, \cdots, 3/2}\pi_{+,>>}\hookrightarrow <\rho\vert\,\vert^{-B}, \cdots, \rho\vert\,\vert^{-3/2}>\times  Jac_{B'_{1}+T'_{1}, \cdots, 3/2}\pi_{+,>>}$$
 soit
 $$
  Jac_{B'_{1}+T'_{1}, \cdots, 3/2}\pi_{+,>>}\hookrightarrow <\rho\vert\,\vert^{-B}, \cdots, \rho\vert\,\vert^{-5/2}>\times  Jac_{B'_{1}+T'_{1}, \cdots, 5/2}\pi_{+,>>},
  $$
  le premier terme dans l'induite n'intervenant pas si $B=3/2$. Le premier cas donne le m\^eme r\'esultat que si $B'_{1}>1/2$ mais dans le deuxi\`eme cas la repr\'esentation $Jac_{B'_{1}+T'_{1}, \cdots, 5/2}\pi_{+,>>}$ se trouve dans le paquet associ\'e au morphisme qui se d\'eduit de $\psi'_{+,>>}$ en remplac\c{c}ant $(\rho,B'_{1}+T'_{1},B'_{1}+T'_{1},+)$ par $(\rho,B'_{1}+1,B'_{1}+1,+)$. Ensuite il faut appliquer $Jac_{B'_{2}+T'_{2}, \cdots, B'_{2}+1}$ et en l'appliquant \`a l'induite de droite on peut encore avoir 2 possibilit\'es si l'on est dans le 2e cas; en tenant compte du fait que $B'_{2}\geq B'_{1}+1$, on contr\^ole les 2 possibilit\'e en faisant appara\^{\i}tre une repr\'esentation qui est dans le paquet associ\'e au morphisme qui se d\'eduit de $\psi'_{+,>>}$ en rempla\c{c}ant $(\rho,B'_{1}+T'_{1},B'_{1}+T'_{1},+)$ par $(\rho,B'_{1}+1,B'_{1}+1,+)$ et $(\rho,B'_{2}+T'_{2},B'_{2}+T'_{2},+)$ soit par $(\rho,B'_{2},B'_{2},+)$ soit par $(\rho,B'_{2}+1,B'_{2}+1,+)$. Et on continue jusqu\`a $i=v$. Donc finalement on trouve un demi-entier $\tilde{B}< B$ avec $\tilde{B}\geq 1/2$ et un morphisme $\tilde{\psi}$ qui se d\'eduit de $\psi$ en rempla\c{c}ant $(\rho,B,B,-)$ par $(\rho,1/2,1/2,+)$ et les $(\rho,B'_{i},B'_{i},+)$ qui v\'erifie $B'_{i}<\tilde{B}$ par $(\rho,B'_{i}+1,B'_{i}+1,+)$ et une repr\'esentation $\tilde{\pi}$ dans le paquet associ\'e \`a $\tilde{\psi}$ avec une inclusion
$$
 \pi\hookrightarrow <\rho\vert\,\vert^{-B}, \cdots, \rho\vert\,\vert^{-\tilde{B}-1}>\times \tilde{\pi},\eqno(1)
$$
o\`u $<\rho\vert\,\vert^{-B}, \cdots, \rho\vert\,\vert^{-\tilde{B}-2}>$ n'intervient pas si $\tilde{B}=B-1$.
Dans le cas (iii) on a une inclusion de m\^eme type que (1) simplement $(\rho,1/2,1/2,+)$ devient $(\rho,0,0,+)$.
On fait rentrer le cas (i)bis dans le cas (1) ci-dessus en posant $\tilde{B}=-1/2$ dans ce cas.

On d\'emontre l'holomorphie de $N(s,\pi)$ en $s=s_{0}$ en utilisant l'inclusion (1); par r\'ecurrence, on sait que $N(s,\tilde{\pi})$ est holomorphe en $s=s_{0}$; on v\'erifie ais\'ement que $r(s,\tilde{\psi})/r(s,\psi)$ est holomorphe avec un z\'ero d'ordre 1 quand $A_{0}=B$. L'op\'erateur d'entrelacement standard $$St(\rho,a_{0})\vert\,\vert^{s}\times <\rho\vert\,\vert^{-B}, \cdots, \rho\vert\,\vert^{-\tilde{B}-1}> \rightarrow
<\rho\vert\,\vert^{-B}, \cdots, \rho\vert\,\vert^{-\tilde{B}-1}>\times St(\rho,a_{0})\vert\,\vert^{s},
$$
est holomorphe en $s=s_{0}$ par positivit\'e.  D'apr\`es \cite{mw} 1.6.3, l'op\'erateur d'entrelacement normalis\'e \`a la Shahidi
$$
<\rho\vert\,\vert^{-B}, \cdots, \rho\vert\,\vert^{-\tilde{B}-1}>\times St(\rho,a_{0})\vert\,\vert^{-s}\rightarrow 
St(\rho,a_{0})\vert\,\vert^{-s}\times <\rho\vert\,\vert^{-B}, \cdots, \rho\vert\,\vert^{-\tilde{B}-1}>,
$$
est holomorphe en $s=s_{0}$ sauf si $A_{0}=\tilde{B}$; dans ce cas on v\'erifie que l'op\'erateur d'entrelacement standard est holomorphe, ce qui r\`egle ce cas. On calcule ais\'ement (ici) le facteur de normalisation \`a une fonction holomorphe inversible pr\`es, il vaut
$$L(St(\rho,a_{0})\times \rho,s-B)/L(St(\rho,a_{0})\times \rho,s-\tilde{B}).$$
En posant $s'=s-s_{0}$, cela vaut encore $L(\rho\times \rho,s'+A_{0}-B)/L(\rho\times \rho,s'+A_{0}-\tilde{B})$. Donc  l'op\'erateur d'entrelacement standard  ne peut avoir de p\^ole que si $A_{0}=B$ c'est-\`a-dire que si $r(s,\tilde{\psi})/r(s,\psi)$ a un z\'ero. D'o\`u l'holomorphie cherch\'ee.

 \end{document}